\documentclass[10pt]{article}
\usepackage[T1]{fontenc}
\usepackage{amsthm}
\usepackage{graphicx}%
\usepackage{float}
\usepackage{placeins}
\usepackage{needspace}
\usepackage{multirow}%
\usepackage{amsmath,amssymb,amsfonts,mathtools}%
\usepackage{mathrsfs}%
\usepackage[title]{appendix}%
\usepackage[table]{xcolor}%
\usepackage{textcomp}%
\usepackage{manyfoot}%
\usepackage{booktabs}%
\usepackage{algorithm}%
\usepackage{algorithmicx}%
\usepackage[noEnd=true,indLines=true]{algpseudocodex}%
\tikzset{algpxIndentLine/.style={draw=black,very thin}}%
\usepackage{arydshln}
\usepackage{amssymb}
\usepackage{xspace}
\usepackage{listings}%
\usepackage{arydshln}
\usepackage{enumitem}
\usepackage[numbers]{natbib}
\usepackage[a4paper]{geometry}
\usepackage[colorlinks=true,
            hypertexnames=false]{hyperref}
\newcommand{\email}[1]{\protect\href{mailto:#1}{#1}}

\usepackage{soul}

\usepackage{aliascnt}
\usepackage[capitalise,noabbrev]{cleveref}
\newcommand{\GITHUBREPO}{\url{https://github.com/ForeverHaibara/strongly_convex_stepsize_concatenation}~\cite{foreverhaibara2026stepsizeconcatenation}}

\newcommand{\x}{{\mathbf x}}

\def\g{{\mathbf g}}
\def\N{{\mathbb N}}
\def\R{{\mathbb R}}
\DeclareMathOperator*{\argmin}{arg\,min}
\theoremstyle{plain}%
\newtheorem{theorem}{Theorem}[section]% meant for sectionwise numbers
\newaliascnt{proposition}{theorem}
\newtheorem{proposition}[proposition]{Proposition}
\aliascntresetthe{proposition}
\crefname{proposition}{Proposition}{Propositions}
\Crefname{proposition}{Proposition}{Propositions}
\newaliascnt{lemma}{theorem}
\newtheorem{lemma}[lemma]{Lemma}
\aliascntresetthe{lemma}
\crefname{lemma}{Lemma}{Lemmas}
\Crefname{lemma}{Lemma}{Lemmas}
\newaliascnt{corollary}{theorem}
\newtheorem{corollary}[corollary]{Corollary}
\aliascntresetthe{corollary}
\crefname{corollary}{Corollary}{Corollaries}
\Crefname{corollary}{Corollary}{Corollaries}

\theoremstyle{plain}%
\newtheorem{remark}{Remark}%
\newtheorem{definition}{Definition} %
\newtheorem{conjecture}{Conjecture}
\crefname{conjecture}{Conjecture}{Conjectures}
\Crefname{conjecture}{Conjecture}{Conjectures}

\numberwithin{equation}{section}

\newcommand{\algA}{{\circ}}
\newcommand{\algB}{{\bullet}}
\newcommand{\algC}{{\blacktriangle}}
\newcommand{\dom}{$f$-dominant\xspace}
\newcommand{\dual}{$g$-bounded\xspace}

\renewcommand{\geqslant}{\geq}

\renewcommand{\a}{{\alpha}}
\renewcommand{\b}{{\beta}}

\newcommand{\rk}[1]{\mathcal{A}_{q}^{(#1)}}
\newcommand{\onet}{{\mathbf{1}^\top}}

\newcommand{\hubertype}{Huber-type function\xspace} 
\DeclareMathOperator*{\argmax}{arg\,max}
\DeclareMathOperator*{\rev}{rev}
\DeclareRobustCommand{\concatname}[1]{%
  \ifmmode
    \mathop{\textnormal{\textsc{#1}}}\nolimits
  \else
    \textnormal{\textsc{#1}}\xspace
  \fi
}
\DeclareRobustCommand{\ConPF}{\concatname{ConPF}}
\DeclareRobustCommand{\ConPP}{\concatname{ConPP}}
\DeclareRobustCommand{\ConGP}{\concatname{ConGP}}
\DeclareRobustCommand{\OConPP}{\textnormal{\texttt{O-ConPP}}\xspace}
\DeclareRobustCommand{\OConPPr}{\textnormal{\texttt{O-ConPPr}}\xspace}
\DeclareRobustCommand{\OConPPl}{\textnormal{\texttt{O-ConPPl}}\xspace}
\DeclareRobustCommand{\OConPF}{\textnormal{\texttt{O-ConPF}}\xspace}
\DeclareRobustCommand{\OConGP}{\textnormal{\texttt{O-ConGP}}\xspace}

\begin{document}

\title{Accelerated Concatenation of Stepsize Schedules for Smooth Strongly Convex Gradient Descent}
\author{
Zehao Zhang\thanks{School of Data Science, Fudan University.}
\and Rujun Jiang\thanks{Corresponding author. School of Data Science, Fudan University. (\email{rjjiang@fudan.edu.cn}).}
}
\date{\today}

\maketitle
\begin{abstract}
\unboldmath 
This work considers stepsize schedules for gradient descent on smooth strongly convex objectives. We develop a curvature-dependent concatenation framework for constructing stepsize schedules of arbitrary length with sharp worst-case guarantees. Within this framework, we introduce two families of schedules for the function-value gap and gradient norm of the last iterate, respectively. For each fixed condition number \(\kappa=L/\mu>1\), both families attain iteration complexity \(\Theta(\kappa^\rho\log(1/\epsilon))\) for sufficiently small \(\epsilon\), where $\rho=\log_{\sqrt2+1}2\approx0.7864$. This complexity is optimal up to a subpolynomial factor in $\kappa$ among predetermined nonnegative stepsize schedules. The framework also recovers the strongly convex silver schedule at power-of-two iteration budgets and applies to smooth convex and smooth weakly convex functions.

\end{abstract}



\section{Introduction}

We consider minimizing an $L$-smooth, $\mu$-strongly convex function $f:\mathbb R^d\to\R$ by gradient descent, where $0<\mu<L$.
% Our certificates and concatenation rules also cover the smooth convex case $\mu=0$ and the weakly convex case $\mu<0$.
% \todo{This sentence is weird.}
We study fixed-step memoryless gradient descent in the form
\begin{equation}\label{eq:gd}
    \x_{i+1}=\x_i-\frac{h_i}{L}\g_i,
    \qquad \g_i=\nabla f(\x_i),\qquad i=0,1,\dotsc,(n-1).
\end{equation}
The iteration budget $n$ and the stepsize schedule $h=[h_0,\dotsc,h_{n-1}]^\top$ are fixed in advance, and the method uses no momentum. Our goal is to construct stepsize schedules of arbitrary length and establish sharp worst-case guarantees.
 
For smooth convex functions, it is known that gradient descent with constant stepsizes has complexity $\Theta (\epsilon^{-1}) $ to reach $\epsilon$-accuracy, while momentum-based methods, e.g., Nesterov's AGM~\cite{nesterov1983method,nesterov2013introductory}, OGM~\cite{kim2016optimized,kim2017ogm} and OGM-G~\cite{kim2021ogm-g} have complexity $\Theta(\epsilon^{-1/2})$.  It was also established that OGM and OGM-G are optimal momentum-based methods  for smooth convex functions in the corresponding metrics~\cite{kim2017ogm,kim2021ogm-g}. 
Adaptive stepsize methods, including Polyak's stepsize~\cite{polyak1987introduction}, the Barzilai--Borwein method~\cite{barzilai1988two}, and recent adaptive gradient methods~\cite{malitsky2020adaptive,zhou2026adabb}, select stepsizes based on observed function values or gradients.

In their seminal 2014 work, Drori and Teboulle~\cite{drori2014performance} proposed a framework, the performance estimation problem (PEP), to analyze the worst-case performance of fixed-step first-order methods.  A series of  
works~\cite{taylor2017thesis,daccache2019performance,eloi2022worst,teboulle2023elementary,grimmer2023longsteps,grimmer2024provably,silver1,silver2,grimmer2024acceleratedobjective,zhang2024accelerated,grimmer2024composing}
focused on designing or analyzing efficient memoryless schedules.
For smooth convex functions, Altschuler and Parrilo~\cite{silver2} proposed the silver stepsize schedule. This schedule and subsequent stepsize constructions~\cite{grimmer2024acceleratedobjective,zhang2024accelerated,grimmer2024composing} attain the best-known iteration complexity $\Theta(\epsilon^{-\rho})$ for memoryless gradient descent, where $\rho=\log_{\sqrt2+1}2\approx0.7864$. Recent studies~\cite{ma2026lower,tsai2026lower,ye2026improved,vernimmen2026can,jung2026stronger,ye2026silver} have proved that the rate is almost optimal  for smooth convex memoryless gradient descent.

For smooth strongly convex problems with a fixed condition number $\kappa=L/\mu>1$, gradient descent with constant stepsizes requires $\Theta(\kappa\log(1/\epsilon))$ iterations to reach $\epsilon$-accuracy.
Momentum-based methods, including Polyak's Heavy Ball method~\cite{polyak1964some} (for quadratics),
 Nesterov's AGM~\cite{nesterov2013introductory}, SC-OGM~\cite{park2023sc-ogm}, ITEM~\cite{taylor2023ITEM}, and ITEM-f~\cite{kim2026lemniscate},
 reduce this complexity to approximately $\Theta(\sqrt{\kappa}\log\frac{1}{\epsilon})$.  The memoryless silver schedule~\cite{silver1}  attains iteration complexity $\Theta(\kappa^{\rho}\log(1/\epsilon))$.

 The original recursive constructions of the silver schedule~\cite{silver1,silver2} and its variants~\cite{grimmer2024acceleratedobjective} are organized around power-of-two iteration budgets.
In our earlier work~\cite{zhang2024accelerated}, we developed a concatenation approach for constructing memoryless fixed-step schedules of arbitrary length for smooth convex functions. The basic idea is to join two shorter stepsize schedules $h_a$ and $h_b$ by inserting a connecting step $\alpha$:
\[
    h_c=[h_a^\top,\alpha,h_b^\top]^\top,
\]
where $\alpha$ is chosen based on $h_a$ and $h_b$. A related composition framework was developed concurrently by Grimmer, Shu, and Wang~\cite{grimmer2024composing}. The approach is also useful in establishing  anytime convergence~\cite{kornowski2024open,zhang2024anytimeaccelerationofgradientdescent} and analyzing constant stepsizes~\cite{grimmer2024strengthened,zhang2026grimmer,chikhi2026last,zhao2026constant}. The schedules proposed by~\cite{zhang2024accelerated,grimmer2024composing} are also conjectured to be the optimal fixed-step memoryless stepsizes for smooth convex optimization. %For smooth strongly convex problems, we study the following design question: how can one construct a stepsize schedule for a prescribed curvature ratio and iteration budget, and certify its worst-case performance in the function-value or gradient-norm metric?

Extending this approach to strongly convex problems requires accounting for the curvature ratio $q=\mu/L=1/\kappa$ in both the concatenation rules and the trajectory certificates, even after normalizing $L$ to one. The main challenge is to construct certificates that yield sharp performance bounds and remain valid under concatenation. The connecting steps must allow the input certificates and interpolation inequalities to combine with nonnegative weights into a certificate for the full schedule.

We develop a curvature-dependent concatenation framework for designing fixed-step memoryless gradient descent schedules on smooth strongly convex functions. We introduce a cumulative stepsize quantity that incorporates curvature and three classes of schedules---primitive, \dom{}, and \dual{}---characterized by inequalities along gradient descent trajectories. We establish three concatenation rules, \ConPP{}, \ConPF{}, and \ConGP{}, for constructing schedules in these classes.
% The trajectory certificates, connecting-step constructions, and growth analysis of the nonlinear cumulative-stepsize recursion differ substantially\todo{not substantial} from those in our smooth convex work~\cite{zhang2024accelerated}.
Using dynamic programming, we construct the \OConPF{} and \OConGP{} schedules of arbitrary length for the function-value gap and gradient norm, respectively.
We prove worst-case bounds for these schedules and construct examples that attain them, which means that these bounds are sharp.
The certificates and concatenation rules also apply to smooth convex and smooth weakly convex functions.

For the function-value gap and gradient-norm criteria, we derive explicit two-sided iteration bounds for \OConPF{} and \OConGP{}, respectively. For each fixed $\kappa=L/\mu>1$ and $0<\epsilon\leq1/(2\kappa^2)$, both families attain iteration complexity $\Theta(\kappa^\rho\log(1/\epsilon))$.
For sufficiently large $\kappa$ and $0<\epsilon\leq1/(2\kappa)$, comparison with the lower bound of Liu and Ye~\cite{liu2026silverstronglyconvex} shows that these complexities are optimal up to a subpolynomial factor in $\kappa$ among all predetermined nonnegative stepsize schedules.
For squared distance and gradient norm, we also construct the \OConPPr{} and \OConPPl{} schedules, respectively. At every power-of-two iteration budget, \OConPPr{} coincides with the strongly convex silver schedule~\cite{silver1}.
The framework also recovers the strongly or weakly convex specialization of the 
dynamic schedule of Rotaru et al.~\cite{rotaru2024exact}.

\textbf{Organization.} \Cref{sect:preliminaries} introduces the curvature-dependent cumulative stepsize and interpolation inequalities.
\Cref{sect:stepsize-shcedules} defines three schedule classes and establishes their worst-case guarantees and concatenation rules.
\Cref{sect:convergence} develops optimized schedules and analyzes their finite-step guarantees and iteration complexity, with remarks on the smooth convex and weakly convex cases.
\Cref{sect:connections-extensions} discusses connections with silver and dynamic schedules.
\Cref{sect:compare} presents numerical comparisons and optimality conjectures.
\Cref{sect:conclusion} concludes the paper.
%Some proofs in the paper involve  cumbersome symbolic computations, and we provide symbolic verification scripts implemented in SymPy at \GITHUBREPO.
Some proofs involve symbolic computations, and the corresponding SymPy verification scripts are available at \GITHUBREPO.

\section{Preliminaries and Notation}
\label{sect:preliminaries}

Let $L>0$ and $\mu<L$. We write $\mathcal F_{\mu,L}(\R^d)$ for the class of differentiable functions such that both $f(\x)-\frac\mu2\|\x\|^2$ and $\frac L2\|\x\|^2-f(\x)$ are convex. This class covers strongly convex, convex, and weakly convex functions for $\mu>0$, $\mu=0$, and $\mu<0$, respectively. The curvature ratio is $q=\mu/L\in(-\infty,1)$. For $q>0$, the condition number is $\kappa=L/\mu=1/q$. 
Throughout the paper, we assume that $f$ admits at least one stationary point. This assumption is automatic when $\mu>0$. Unless stated otherwise, let $\x_\star$ be any stationary point, with $f_\star=f(\x_\star)$ and $\g_\star\triangleq\nabla f(\x_\star)=\mathbf 0$. When $\mu\geq0$, $\x_\star$ is a global minimizer. For $\mu<0$, any additional assumption of global optimality or boundedness below is stated where needed.
%In the remainder of the paper, $q$ is fixed. In iteration-complexity statements for $q>0$, $\kappa$ is fixed in advance and the precision $\epsilon$ varies.

The endpoint $\mu=L$, or equivalently $q=1$, can be handled separately. In this case, the two convexity conditions imply that $f(\x)-\frac L2\|\x\|^2$ is affine, so $f(\x)=f_\star+\frac L2\|\x-\x_\star\|^2$. Thus, one gradient step with stepsize $1/L$, corresponding to $h_0=1$ in~\eqref{eq:gd}, reaches $\x_\star$ from any initial point. We therefore restrict the subsequent analysis to $\mu<L$.

We write $f_i=f(\x_i)$ and $\g_i=\nabla f(\x_i)$ for the function value and gradient at the iterate $\x_i$. We use $h=[h_0,\dotsc,h_{n-1}]^\top\in\R^n$ to denote a schedule of $n$ iterations. The norm $\|\cdot \|$ is always the Euclidean norm. We use $\N=\{0,1,2,\dotsc\}$ and $\N_+=\{1,2,3,\dotsc\}$ to denote the sets of natural numbers and positive integers, respectively.
We use $\mathbb Q(a_1,...,a_m)$ to denote the rational function field with rational coefficients.
In~\cref{tab:notation} we record other notations that will be defined and used  throughout the paper.

\begin{table}[htbp]
    \centering
    \caption{Notation used throughout the paper.}
    \label{tab:notation}
    \small
    \begin{tabular}{@{}lll@{}}
        \toprule
        Notation & Definition & Description \\
        \midrule
        % $q$ & $ $ & Ratio of the lower and upper curvatures.\\
        % $\kappa $ &  & Condition number when $\mu>0$. \\ 
        $\sigma_q$ &\eqref{eq:sigma} & Cumulative quantity of a schedule.\\
        $\rk{n}$ & \Cref{def:rkn} & The domain of considered schedules.\\
      %  $Q_{i,j}$ & \Cref{def:Qij} & Co-coercivity inequalities $(Q_{i,j}\leq 0)$.\\
      %  $\ConPP$, $\ConPF$, $\ConGP$ & \Cref{theorem:conpp,theorem:conpf,theorem:congp} & Stepsize concatenation procedure. \\
        $\varphi_q,\psi_q,\Psi_{x,y,q}$ & \Cref{theorem:conpp,theorem:conpf,theorem:congp}  & Operators used in concatenation procedure.\\
        $h_\algA,h_\algB,h_\algC$ & \cref{def:dynamic-programming} & \OConPP, \OConPF, and \OConGP{} schedules. \\

        $\nu(q)$ &\cref{prop:prim-asymp} & Asymptotic rate of the cumulative stepsize.\\
        $\hat h,\check h$ & \cref{sect:silver-stepsizes} & \OConPPr and \OConPPl schedules. \\
        \bottomrule
    \end{tabular}
\end{table}

We first introduce a curvature-dependent cumulative quantity of a schedule $h=[h_0,\dotsc,h_{n-1}]^\top$.
\begin{equation*}
    \label{eq:sigma}\tag{Def-$\sigma$}
    \sigma_{q}(h) = \begin{cases}
        q^{-1} \left(\prod_{i=0}^{n-1}(1-q h_i)^{-1}-1\right), &\text{if } q \in (-\infty,0)\cup (0,1),
        \\
        \sum_{i=0}^{n-1}h_i, &\text{if } q=0.
    \end{cases}
\end{equation*}
As $q\to 0$, $(1-q h_i)^{-1}=1+q h_i+o(q)$ and hence $\lim_{q\to 0}\sigma_q(h)=\sigma_0(h)=\onet h$. Thus, $\sigma_q(h)$ extends the stepsize sum from the convex case $q=0$ to the strongly and weakly convex settings. In our earlier work~\cite{zhang2024accelerated}, the sum $\onet h$ plays a central role in the concatenation analysis for smooth convex functions. Here, $\sigma_q(h)$ takes its place in the strongly and weakly convex settings.
As $q$ is fixed, we write $\sigma_q(h)$ as $\sigma(h)$ whenever no ambiguity can arise.
We adopt the convention that an empty sum equals $0$ and an empty product equals $1$. Hence $\sigma_q ([~])=0$ for an empty schedule.

The factors in the definition~\eqref{eq:sigma} require $1-q h_i\neq 0$. In this paper we work with the slightly stronger condition $1-q h_i>0$, which guarantees that all factors remain positive. 
This leads to the following definition of \emph{admissible} schedules.

\begin{definition}\label{def:rkn}
    For $n\in\N$, we define the set of admissible schedules of length $n$ by
    \[
    \rk{n}=
    \{[h_0,\dotsc,h_{n-1}]^\top:\ h_i\geq 0,\ 1-q h_i> 0,\ (i=0,\dotsc,n-1)\}    .\]
\end{definition}
Thus, the admissible stepsizes satisfy $0\leq h_i<1/q$ when $q>0$ and $h_i\geq0$ when $q\leq0$.

We next establish useful properties of $\sigma_q$ in~\cref{lemma:telescope,lemma:s-bound,lemma:recursion-of-sigma}.

\begin{lemma}\label{lemma:telescope}
Let $h=[h_0,\dotsc,h_{n-1}]^\top\in \rk{n}$. Then we have
\begin{subequations}\begin{align}\label{eq:qs-add-1}
 \prod_{i=0}^{n-1}(1-q h_i)^{-1}&= q \sigma_q(h)+1,\\
    \label{eq:telescope}
\sum_{i=0}^{n-1}h_i\prod_{j=0}^i(1-q h_j)^{-1}&=\sigma_q (h).
\end{align}
\end{subequations}
\end{lemma}

\begin{proof} Both~\eqref{eq:qs-add-1} and~\eqref{eq:telescope} are trivial when $q=0$. When $q\neq 0$,
    ~\eqref{eq:qs-add-1} directly follows from the definition~\eqref{eq:sigma}, while~\eqref{eq:telescope} follows from
    \[
 \begin{aligned}
\sum_{i=0}^{n-1} h_i\prod_{j=0}^i(1-q h_j)^{-1}
    &= \frac{1}{q}\sum_{i=0}^{n-1} ({1-(1-q h_i)})\prod_{j=0}^i(1-q h_j)^{-1}=
   \frac{1}{q}\left( \prod_{i=0}^{n-1}(1-q h_i)^{-1}-1\right)
   \overset{\eqref{eq:sigma}}{=}\sigma_q(h),
 \end{aligned}
    \]
where the second equality follows from telescoping the sum.\end{proof}

\begin{lemma}[Positivity and monotonicity of $\sigma_q$] \label{lemma:s-bound}
    Let $h\in \rk{n}$. Then we have
    \[
    \sigma_q(h)\geq 0,\qquad 1+q\sigma_q(h)>0,\qquad 1+2\sigma_q(h)+q\sigma_q^2(h)>0.
    \]
    Moreover, $\sigma_q$ is monotonically increasing in each coordinate on $\rk{n}$.
\end{lemma}

\begin{proof}
    The first two inequalities follow from~\cref{lemma:telescope}, since $h_i\geq 0$ and $1-q h_i>0$ for every $i$. Consequently,
    \[
    1+2\sigma_q(h)+q\sigma_q^2(h)
    =\sigma_q(h)\bigl(1+q\sigma_q(h)\bigr)+\sigma_q(h)+1>0.
    \]
    The entrywise monotonicity follows from the definition~\eqref{eq:sigma}.\end{proof}

\begin{lemma}[Concatenation identity for $\sigma_q$]\label{lemma:recursion-of-sigma} If $h_c=[h_a^\top,\alpha,h_b^\top]^\top$ and $h_c\in\rk{n}$, then
\begin{equation}
    \label{eq:recursion-of-sigma}
\sigma_q(h_c)=\frac{\sigma_q(h_a)+\alpha+\sigma_q(h_b)+q\sigma_q(h_a)\sigma_q(h_b)}{1-q\alpha}.
\end{equation}
When $q\neq 0$, the identity can also be written as
\begin{equation}
    \label{eq:recursion-of-sigma2}
\sigma_q(h_c)=\frac{(q \sigma_q(h_a)+1)(q \sigma_q(h_b)+1)(1-q\alpha)^{-1}-1}{q}.
\end{equation}
\end{lemma}

\begin{proof}
Splitting the sum in~\eqref{eq:telescope} into the prefix $h_a$, the connecting step $\alpha$, and the suffix $h_b$, and using~\eqref{eq:qs-add-1}, gives
\[
\begin{aligned}
\sigma_q(h_c)
&=\sigma_q(h_a)+\frac{1+q\sigma_q(h_a)}{1-q\alpha}\bigl(\alpha+\sigma_q(h_b)\bigr)=\frac{\sigma_q(h_a)+\alpha+\sigma_q(h_b)+q\sigma_q(h_a)\sigma_q(h_b)}{1-q\alpha}.
\end{aligned}
\]
This proves~\eqref{eq:recursion-of-sigma}. For $q\neq 0$, the right-hand sides of~\eqref{eq:recursion-of-sigma2}
and~\eqref{eq:recursion-of-sigma} agree by direct expansion.
\end{proof}

%The first lemma records the curvature-dependent telescoping identity, and the remaining results provide the weighted inequalities needed in the subsequent sections.
%Most of them are proved by telescoping summation over inequalities $Q_{i,j}\leq 0$ defined in~\cref{def:Qij}.
We next recall the interpolation inequalities in the following lemma. The co-coercivity inequalities in~\eqref{eq:Qij} for functions with upper and lower curvature bounds follow from~\cite[Theorem 3.2]{rotaru2024exact}, using $\g_\star=\mathbf 0$ when one of the points is $\x_\star$. The additional bound in~\eqref{eq:fi*} follows from~\cite[Lemma 3.5]{rotaru2024exact}.
\begin{lemma}\label{def:Qij} Suppose $L>0$, $\mu<L$, and $f\in\mathcal F_{\mu,L}(\mathbb R^d)$. For every $i,j\in\{0,\dotsc,n\}$, the quantities $Q_{i,j}$, $Q_{i,\star}$, and $Q_{\star,j}$ defined below satisfy
\[
Q_{i,j}\leq0,\qquad Q_{i,\star}\leq0,\qquad Q_{\star,j}\leq0,
\]
where
\begin{equation}\label{eq:Qij}
\begin{aligned}
Q_{i,j}&\triangleq f_i-f_j+\langle \g_i,\x_j-\x_i\rangle   +\frac{1}{2L}\|\g_i-\g_j\|^2 
%\\  &\qquad\qquad
+\frac{\mu}{2(1-\mu/L)}\|\x_i-\x_j-\frac{1}{L}(\g_i-\g_j)\|^2,\\
Q_{i,\star}&\triangleq f_i - f_\star+\langle \g_i,\x_\star-\x_i\rangle +\frac{1}{2L}\|\g_i\|^2 + \frac{\mu}{2(1-\mu/L)}\|\x_i-\x_\star-\frac{1}{L}\g_i\|^2,\\
Q_{\star,j}&\triangleq f_\star-f_j +\frac{1}{2L}\|\g_j\|^2 + \frac{\mu}{2(1-\mu/L)}\|\x_j-\x_\star-\frac{1}{L}\g_j\|^2.
\end{aligned}
\end{equation}
If $f_\star$ is the global infimum, then
\begin{equation}\label{eq:fi*}
    f_i -f_\star\geq \frac{1}{2L}\|\g_i\|^2.
\end{equation}
\end{lemma}

The following lemma relates weighted sums along a gradient descent trajectory to its endpoints and to the interpolation quantities $Q_{i,j}$ from~\cref{def:Qij}.
\begin{lemma}\label{lemma:multistep}
Consider gradient descent with a schedule $h\in\rk{n}$. Let $k,l,r\in\{0,\dotsc,n\}$ with $l\leq r$, and set $s=\sigma_q([h_l,\dotsc,h_{r-1}]^\top)$. Then the following statements hold.
\begin{enumerate}[label=(\roman*),ref=(\roman*),leftmargin=*]
\item\label{item:multistep-telescope}
We have
\begin{equation}
    \label{eq:sum-of-x}
    \sum_{i=l}^{r-1}h_i\prod_{j=l}^i (1-q h_j)^{-1}\left(\frac{\g_i}L-q(\x_i-\x_\star)\right)=(\x_l-\x_\star)-(q s+1)(\x_r-\x_\star).
\end{equation}
\item\label{item:multistep-weights}
For any $u=[u_l,\dotsc,u_{r-1}]^\top\in\R_+^{r-l}$,
\begin{equation}
\label{eq:A4}    
\begin{aligned}
% &\frac1L\left\langle \g_k-L(\x_k-\x_\star),\sum_{i=l}^{r-1}u_i\bigl(\g_i-\mu(\x_i-\x_\star)\bigr)\right\rangle
& \left\langle \g_k-L(\x_k-\x_\star),\sum_{i=l}^{r-1}u_i\bigl(\frac{\g_i}{L}-q(\x_i-\x_\star)\bigr)\right\rangle
-(1-q)\left[(\onet u)Q_{\star,k}+\sum_{i=l}^{r-1}u_iQ_{i,\star}\right]\\
&=-(1-q)\sum_{i=l}^{r-1}u_iQ_{i,k}\geq0.
\end{aligned}
\end{equation}
\end{enumerate}
\end{lemma}

\begin{proof}
For~\ref{item:multistep-telescope}, the gradient descent recurrence~\eqref{eq:gd} gives
\[
h_i\left(\frac{\g_i}L-q(\x_i-\x_\star)\right)=(1-q h_i)(\x_i-\x_\star)-(\x_{i+1}-\x_\star).
\]
Multiplying by $\prod_{j=l}^i(1-q h_j)^{-1}$ and telescoping over $i=l,\dotsc,r-1$, we obtain
\[
\sum_{i=l}^{r-1}h_i\prod_{j=l}^i(1-q h_j)^{-1}\left(\frac{\g_i}L-q(\x_i-\x_\star)\right)
=(\x_l-\x_\star)-\prod_{j=l}^{r-1}(1-q h_j)^{-1}(\x_r-\x_\star).
\]
Using $\prod_{j=l}^{r-1}(1-q h_j)^{-1}=1+q s$ from~\eqref{eq:qs-add-1}, we obtain~\eqref{eq:sum-of-x}.

For~\ref{item:multistep-weights}, expansion of the quantities in~\cref{def:Qij} yields
\[
Q_{i,\star}+Q_{\star,k}-\frac{1}{L(1-q)}\langle\g_k-L(\x_k-\x_\star),\bigl(\g_i-\mu(\x_i-\x_\star)\bigr)\rangle=Q_{i,k}\leq0.
\]
Multiplying by $-(1-q)u_i$ and summing over $i=l,\dotsc,r-1$ gives~\eqref{eq:A4}.\end{proof}

\section{Schedule Classes and Concatenation}
\label{sect:stepsize-shcedules}
In this section, we define three classes of stepsize schedules through curvature-dependent inequalities: primitive, \dom, and \dual.
Their smooth convex counterparts were introduced in our earlier work~\cite{zhang2024accelerated}.
For \dom and \dual schedules, we establish sharp function-value and gradient-norm bounds in terms of $\sigma_q(h)$, respectively.
We then develop concatenation rules that join two schedules by inserting a suitably chosen step between them.
The key is to construct connecting steps and weights so that the resulting schedules satisfy the defining inequalities of the required classes.

\subsection{Schedule Classes and Performance Guarantees}

For every $q<1$, we define three schedule classes through inequalities along gradient descent trajectories. In the following definitions, $f\in\mathcal F_{\mu,L}$ and the starting point $\x_0$ are arbitrary, and $\x_\star$ denotes any stationary point of $f$. The prefix $f$ in \dom refers to function-value guarantees.

\begin{definition}[Primitive schedule]\label{def:primitive} A schedule $h=[h_0,\dotsc,h_{n-1}]^T\in\rk{n}$ is primitive if 
\begin{equation}
\label{eq:primitive}
\begin{aligned}    
&\frac L2\|\x_0  -\x_\star\|^2+(1-q)\sum_{i=0}^{n-1}h_i\prod_{j=0}^{i}(1-q h_j)^{-1}Q_{i,\star}+q \sigma(h)Q_{n,\star}-\sigma(h)(f_n-f_\star)\\ &
\geq \frac L2\|\x_n-\x_\star\|^2+\frac{\sigma(h)(\sigma(h)+1)}{2L}\|\g_n\|^2+\frac{q \sigma(h)}{2L(1-q)}\|\g_n-L(\x_n-\x_\star)\|^2.
\end{aligned}
\end{equation}
\end{definition}

\begin{definition}[\dom schedule]\label{def:f-dominant} A schedule $h=[h_0,\dotsc,h_{n-1}]^\top\in\rk{n}$ is \dom if there exists $u = [u_0,\dotsc,u_{n-1}]^\top\in \R_+^{n}$, independent of $f$, $\x_0$, and $\x_\star$, such that $\eta-(1-q)(\onet u)\geq 0$ and
% \scriptsize 
% ...
% \scalebox{0.9}{
{\small\begin{equation}
\label{eq:dominant} 
\hspace{-1.5em}\begin{aligned}    
&\frac L2\|\x_0- \x_\star\|^2+(1-q)\sum_{i=0}^{n-1} u_iQ_{i,\star}
+ (\eta-(1-q)(\onet u))Q_{n,\star}  -\eta(f_n-f_\star)
% --------------------------------------
% \\ &\geq \frac {L}{2(q s+1)^2}\left\|\x_0-\x_\star- \sum_{i=0}^n u_{i}\prod_{j=0}^{i}(1-q u_j)^{-1}(\g_i/L-q (\x_i-\x_\star))-\frac{q(\onet u)}{1-q}(\g_n/L-q (\x_n-\x_\star))\right\|^2.
% --------------------------------------
\\ &\geq  \frac { (1-q)L}{2(1+q \sigma(h))^2}\left\|\x_0-\x_\star- \sum_{i=0}^{n-1} u_{i}\left(\frac{\g_i}L-q(\x_i-\x_\star)\right)-\frac{\eta-(1-q)(\onet u)}{1-q}\left(\frac{\g_n}L-q(\x_n-\x_\star)\right)\right\|^2.
 \end{aligned} 
 \end{equation}}
% }
where $\eta\triangleq1+2\sigma(h)+q\sigma^2(h)$.
\end{definition}

\begin{definition}[\dual schedule]\label{def:g-dominant} A schedule $h = [h_0,\dotsc,h_{n-1}]^T\in \rk{n}$ is \dual if
\begin{equation}
    \label{eq:g-dominant}f_0 -f_n\geq (2\sigma(h)+q \sigma^2(h))\cdot \frac{1}{2L}\|\g_n\|^2.
    %-\frac{q L}{2(1-q)}\|\x_n-\x_\star-\g_n/L\|^2.
\end{equation}
\end{definition}

The three schedule classes depend on the curvature ratio $q$. For example, a schedule that is \dom for $q=1/2$ need not be \dom for $q=0$.

We use one-dimensional Huber-type functions with $x_\star=0$ to construct examples attaining the performance bounds. For $r\geq0$, define
\begin{equation}\label{eq:huber-f}
f(x)=\begin{cases}
    \frac L2 x^2,& \text{if }|x|<r,\\
    \frac\mu2 x^2+(L-\mu)r|x|-\frac{L-\mu}{2}r^2,&\text{otherwise}.
\end{cases}
\end{equation}
This function belongs to $\mathcal F_{\mu,L}(\R)$ because subtracting $\frac\mu2 x^2$ yields a convex, $(L-\mu)$-smooth Huber function. The origin is stationary, with $f_\star=f(0)=0$, and is a global minimizer when $\mu\geq0$.

The following lemma gives the gradient descent iterates for a suitable choice of $r$ and shows that they remain in $[r,+\infty)$. The proof is given in Appendix~\ref{sect:proofs-main-lemmas}.
\begin{lemma}\label{lemma:huber-gd}
Let $h\in\rk{n}$, $s\geq1+\sigma(h)$, and $x_0\geq0$. Consider the function $f$ in~\eqref{eq:huber-f} with $r=x_0/s$. Then gradient descent with schedule $h$ satisfies $x_k\geq r$ and
\begin{equation}\label{eq:huber-xk}
x_k=\frac{1-\frac{(1-q)\sigma^{(k)}}{s}}{1+q\sigma^{(k)}}x_0,\qquad k=0,1,\dotsc,n,
\end{equation}
where $\sigma^{(k)}=\sigma([h_0,\dotsc,h_{k-1}]^\top)$.
\end{lemma}

% These examples also give equality in any dimension by taking
% \[
% \widetilde f(\x)=f(x_1)+\frac L2\sum_{j=2}^d x_j^2,\qquad
% \x_0=(x_0,0,\dotsc,0)^\top,
% \]
% since the remaining coordinates stay zero under gradient descent.\todo{delete?}

The next two theorems give sharp bounds on the function-value gap and gradient norm in terms of $\sigma_q(h)$.
\begin{theorem}[Sharp function-value bound]\label{theorem:bound-of-fdom} Let $\x_\star$ be any stationary point of $f$. If $h=[h_0,\dotsc,h_{n-1}]^\top$ is \dom, then the iterates generated by gradient descent with schedule $h$ satisfy
\begin{equation}\label{eq:bound-fdom}
    f_n-f_\star\leq \frac{1}{1+2\sigma(h)+q\sigma^2(h)}\cdot \frac{L}{2}\|\x_0-\x_\star\|^2.
\end{equation}
Equality is attained by~\eqref{eq:huber-f} with $x_0>0$ and $r=\frac{x_0}{1+2\sigma(h)+q\sigma^2(h)}$.
\end{theorem}

\begin{proof}
Since $Q_{i,\star}\leq 0$ by~\cref{def:Qij}, $1-q\geq 0$ and $\eta-(1-q)\onet u\geq 0$ by~\cref{def:f-dominant}, and $\eta=1+2\sigma(h)+q\sigma^2(h)>0$ by~\cref{lemma:s-bound},~\eqref{eq:bound-fdom} follows from~\eqref{eq:dominant} immediately.

We next verify equality for~\eqref{eq:huber-f} with the stated value of $r$.
Denote $\sigma^{(k)}=\sigma([h_0,\dotsc,h_{k-1}]^\top)$. Since
\[
1+2\sigma(h)+q\sigma^2(h)=1+\sigma(h)+\sigma(h)(q\sigma(h)+1)\geq 1+\sigma(h),
\]
\cref{lemma:huber-gd} gives $x_k\geq r$ and
\begin{equation*}
x_k=\frac{1-\frac{(1-q)\sigma^{(k)}}{1+2\sigma(h)+q\sigma^2(h)}}{1+q\sigma^{(k)}}x_0,\qquad k=0,1,\dotsc,n.
\end{equation*}
Substituting this expression at $k=n$ into~\eqref{eq:huber-f} yields equality in~\eqref{eq:bound-fdom}.\end{proof}

\begin{theorem}[Sharp gradient-norm bound]\label{theorem:bound-of-dual} Suppose that $f_\star$ is the global infimum. If $h=[h_0,\dotsc,h_{n-1}]^\top$ is \dual, then the iterates generated by gradient descent with schedule $h$ satisfy
\begin{equation}\label{eq:bound-gdom}
    \frac{1}{2L}\|\g_n\|^2\leq \frac{1}{1+2\sigma(h)+q\sigma^2(h)}(f_0 - f_\star).
\end{equation}

For the equality cases, take $x_0>0$. If $q\geq0$, equality is attained by~\eqref{eq:huber-f} with $r=x_0/(1+\sigma(h))$.
If $q<0$, equality is attained by the following modified \hubertype{} with the same choice of $r$:
\begin{equation}\label{eq:huber-g-neg}
f(x)=\begin{cases}
    \frac L2 x^2,& \text{if }|x|<r,\\
    \frac\mu2 x^2+(L-\mu)r|x|-\frac{L-\mu}{2}r^2,& \text{if }r\leq |x|\leq x_0,\\
    (1+q\sigma(h))Lr|x|-\frac L2(1+2q\sigma(h)+q\sigma^2(h))r^2,& \text{otherwise}.
\end{cases}
\end{equation}

\end{theorem}

\begin{proof}
Since $f_\star$ is the global infimum,~\eqref{eq:fi*} shows $f_n-f_\star\geq \frac{1}{2L}\|\g_n\|^2$. Adding this inequality to~\eqref{eq:g-dominant} gives $f_0-f_\star\geq (1+2\sigma(h)+q\sigma^2(h))\cdot \frac{1}{2L}\|\g_n\|^2$. Since $1+2\sigma(h)+q\sigma^2(h)>0$ by~\cref{lemma:s-bound}, the desired bound~\eqref{eq:bound-gdom} follows.

For $q\geq0$, applying~\cref{lemma:huber-gd} to~\eqref{eq:huber-f} with the stated value of $r$ gives
\begin{equation*}
x_n=\frac{1-\frac{(1-q)\sigma(h)}{1+\sigma(h)}}{1+q\sigma(h)}x_0
=\frac{x_0}{1+\sigma(h)}=r.
\end{equation*}
Hence $f'(x_n)=Lr$ and
\[
\begin{aligned}
f_0-f_\star
&=\frac\mu2 x_0^2+(L-\mu)rx_0-\frac{L-\mu}{2}r^2=\frac L2\bigl(1+2\sigma(h)+q\sigma^2(h)\bigr)r^2
=\frac{1+2\sigma(h)+q\sigma^2(h)}{2L}|f'(x_n)|^2,
\end{aligned}
\]
which gives equality in~\eqref{eq:bound-gdom}.

For $q<0$, we use~\eqref{eq:huber-g-neg}. Since $x_0=(1+\sigma(h))r$ and $\mu=qL$, we have
\begin{equation}\label{eq:huber-g-neg-func}
\begin{aligned}
\frac\mu2 x_0^2+(L-\mu)rx_0-\frac{L-\mu}{2}r^2
&=\frac L2\bigl(1+2\sigma(h)+q\sigma^2(h)\bigr)r^2\\
&=(1+q\sigma(h))Lrx_0-\frac L2\bigl(1+2q\sigma(h)+q\sigma^2(h)\bigr)r^2,
\end{aligned}
\end{equation}
and
\begin{equation*}
\mu x_0+(L-\mu)r=(1+q\sigma(h))Lr>0.
\end{equation*}
These identities show that the middle and outer branches agree in value and derivative at $x=x_0$, and evenness gives the same conclusion at $x=-x_0$. Since $f$ agrees with~\eqref{eq:huber-f} on $[-x_0,x_0]$, it is continuously differentiable. Away from the joining points,
\[
f''(x)=\begin{cases}
L,& |x|<r,\\
\mu,& r<|x|<x_0,\\
0,& |x|>x_0.
\end{cases}
\]
Thus $\mu\leq f''(x)\leq L$ on each interval, so $f\in\mathcal F_{\mu,L}(\R)$. On each branch, $f$ is strictly increasing as a function of $|x|$. Hence $x_\star=0$ is a global minimizer with $f_\star=0$.

For $s=1+\sigma(h)$, the formula in~\cref{lemma:huber-gd} gives $r\leq x_k\leq x_0$ for $k=0,\dotsc,n$, using $1+q\sigma(h)>0$. Changing the function outside $[-x_0,x_0]$ therefore leaves the iterates unchanged. In particular, $x_n=r$ and $f'(x_n)=Lr$. Together with $f_0-f_\star=(1+2\sigma(h)+q\sigma^2(h))Lr^2/2$ from~\eqref{eq:huber-g-neg-func}, this gives equality in~\eqref{eq:bound-gdom}.
\end{proof}

The following proposition shows that every primitive schedule is also \dom and therefore satisfies the function-value bound in~\cref{theorem:bound-of-fdom}.
\begin{proposition}\label{prop:primitive-is-fdom}
Let $h=[h_0,\dotsc,h_{n-1}]^\top\in\rk{n}$ be a primitive schedule. Define $u=[u_0,\dotsc,u_{n-1}]^\top$
by $u_i=h_i\prod_{j=0}^i(1-q h_j)^{-1}$
for $i=0,\dotsc,n-1$.
Then $h$ is \dom with the vector $u$ in~\cref{def:f-dominant}.
\end{proposition}

\begin{proof}
Set $s=\sigma(h)$.
Since $h\in\rk{n}$ implies $h_i\geq 0$ and $q h_j<1$, we have $u\in\R_+^n$ by definition.
Also,
\[
\onet u=\sum_{i=0}^{n-1}h_i\prod_{j=0}^i(1-q h_j)^{-1}  \stackrel{\eqref{eq:telescope}}{=}s.
\] 
Hence $(1+2s+q s^2)-(1-q)(\onet u)= (1+2s+q s^2)-(1-q)s=(1+s)(q s+1)\geq 0$, where $q s+1>0$ follows from~\cref{lemma:s-bound}.

It remains to show that~\eqref{eq:dominant} holds. The vector inside the norm on its right-hand side satisfies
\begin{equation}\label{eq:primitive-is-dominant-rhs}\begin{aligned}&\ 
\x_0-\x_\star- \sum_{i=0}^{n-1} u_i\left(\frac{\g_i}L-q(\x_i-\x_\star)\right)-\frac{(1+s)(q s+1)}{1-q}\left(\frac{\g_n}L-q(\x_n-\x_\star)\right)
% \\ &= (q s+1) (\x_n-\x_\star)-\left(u_n+\frac{q (\onet u)}{1-q}\right)(\g_n/L-q  (\x_n-\x_\star))
% \\ &= (q s+1) (\x_n-\x_\star)-\left(1+s+q s^2+\frac{q (1+2s+q s^2)}{1-q}\right)(\g_n/L-q  (\x_n-\x_\star))
\\ &=(q s+1)(\x_n-\x_\star)- \frac{(q s+1)(s+1)}{1-q} \left(\frac{\g_n}L-q(\x_n-\x_\star)\right)
\\ &=\frac{q s+1 }{1-q}\left((q s+1)(\x_n-\x_\star) - (s+1)\frac{\g_n}{L}\right),
\end{aligned}
%+u_n\prod_{j=0}^i(1-q h_j)^{-1}(1-q u_n)^{-1}\frac{(1-q h_{n-1})x_{n-1}-x_n}{(1-q )h_{n-1}}.
\end{equation}
where the first equality follows from~\eqref{eq:sum-of-x} with $l=0$ and $r=n$.
For the left-hand side of~\eqref{eq:dominant}, we have
\begin{equation}\label{eq:primitive-is-dominant-lhs}
\begin{aligned}    
&\frac L2\|\x_0- \x_\star\|^2+(1-q)\sum_{i=0}^{n-1} u_{i} Q_{i,\star}+ (1+2s+q s^2-(1-q)s)Q_{n,\star}  -(1+2s+q s^2)(f_n-f_\star)\\ &\geq  
\frac L2\|\x_n-\x_\star\|^2+\frac{s(s+1)}{2L}\|\g_n\|^2+\frac{q s}{2L(1-q)}\|\g_n-L(\x_n-\x_\star)\|^2
\\ &\qquad +(1+s+q s^2)Q_{n,\star}-(1+s+q s^2)(f_n-f_\star)
% --------------------------------------
% \\ &=\frac L2\|\x_n-\x_\star\|^2+\frac{s(s+1)}{2L}\|\g_n\|^2+\frac{q Ls}{2(1-q)}\|\x_n-\x_\star - \g_n/L\|^2
% \\ &\qquad +(1-q)(1+s+q s^2)Q_{n,\star}+(1+s+q  s^2)(q Q_{n,\star}-f_n+f_\star)
% --------------------------------------
\\ & =\frac L2\cdot\frac{1}{1-q}\|(q s+1)(\x_n-\x_\star)-(s+1)\g_n/L\|^2,
% \\ &=\frac L2\|\x_n-\x_\star\|^2+\frac{s(s+1)}{2L}\|\g_n\|^2+\frac{q (1+2s+q s^2)}{2(1-q)}\|\x_n-\x_\star - \g_n/L\|^2
\end{aligned}
\end{equation}
where the inequality follows by adding $(1+s+q s^2)(Q_{n,\star}-(f_n-f_\star))$ to both sides of~\eqref{eq:primitive}, and the equality follows by substituting the definition of $Q_{n,\star}$ in~\eqref{eq:Qij} and completing the square.
Comparing \eqref{eq:primitive-is-dominant-lhs} with \eqref{eq:primitive-is-dominant-rhs}, we conclude that  \eqref{eq:dominant} holds and thus $h$ is \dom.\end{proof}

For the empty schedule, $\sigma([~])=0$, so both~\eqref{eq:primitive} and~\eqref{eq:g-dominant} hold with equality. Together with~\cref{prop:primitive-is-fdom}, this gives the following corollary.
\begin{corollary}\label{coro:zero-step} The empty schedule $[~]\in\R^0$ is primitive, \dom and \dual.
\end{corollary}

\subsection{Concatenation of Schedules}
\label{sect:concat}
This subsection introduces three concatenation operations, \ConPP, \ConPF{} and \ConGP. When $q=0$, these operations recover the concatenation rules in~\cite{zhang2024accelerated}, as discussed for smooth convex functions at the end of~\cref{sect:convergence}.
 
\begin{theorem}[Concatenation of primitive schedules]\label{theorem:conpp}
Suppose $h_a\in \rk{n}$ and $h_b\in \rk{m-n-1}$ are both primitive, where $m\geq n+1$. Let $s_a=\sigma(h_a)$ and $s_b=\sigma(h_b)$. Then $h_c=[h_a^\top,\a,h_b^\top]^\top\in\rk{m}$ is also primitive, where $\a=\varphi_q(s_a,s_b)$ and $\varphi_q$ is defined by 
\begin{equation*}
    \label{eq:phi}\tag{Def-$\varphi$}
\varphi_q(x,y)\!=\frac{-x-y-\frac{2q}{1-q}(x+1)(y+1)+\sqrt{(x+y+2+\frac{2q}{1-q}(x+1)(y+1))^2+4(x+1)(y+1)}}{2}.
\end{equation*}
\end{theorem}

\begin{theorem}[Concatenation of \dom schedules]\label{theorem:conpf}
Suppose $h_a\in \rk{n}$ is primitive and $h_b\in \rk{m-n-1}$ is \dom, where $m\geq n+1$. Let $s_a=\sigma(h_a)$ and $s_b=\sigma(h_b)$. Then $h_c=[h_a^\top, \b,h_b^\top]^\top\in\rk{m}$ is \dom, where $\b=\psi_q(s_a,s_b)$ is the unique root of the equation $\Psi_{s_a,s_b,q}(\b)=0$ in the interval $(1,s_a+2)$:
\begin{equation*}
\label{eq:Psi}\tag{Def-$\Psi$}%\footnotesize 
% \Psi_{x,y,q}(h)=\begin{cases}
% q^{-2}
% \left[\begin{aligned}
%   &  (1-q h)^3+(qx+1-2(1-q))(1-q h)^2\\ &-(2(qx+1)-(1-q))(1-q h)(qy+1)^2+(1-q)(qx+1)(qy+1)^2
% \end{aligned}\right],& \text{if}\ q \neq 0,\\
% 2h^2+(2y-3)h    -2xy-x-4y,\ \text{if}\ q=0.\end{cases}
\Psi_{x,y,q}(h)=(h-1)^2(2+q(x-h))-(qy^2+2y+1)(qx+x+2-(2qx+q+1)h).
\end{equation*}
% Moreover, $\Psi_{x,y,q}(h)\to\Psi_{x,y,0}(h)$ as $q\to0$, and the equation $\Psi_{x,y,0}(h)=0$ is quadratic.
\end{theorem}

\begin{theorem}[Concatenation of \dual schedules]\label{theorem:congp} Suppose $h_a\in \rk{n}$ is \dual and $h_b\in \rk{m-n-1}$ is primitive, where $m\geq n+1$. Let $s_a=\sigma(h_a)$ and $s_b=\sigma(h_b)$. Then $h_c=[h_a^\top, \b,h_b^\top]^\top\in\rk{m}$ is \dual, where $\b=\psi_q(s_b,s_a)$ is the unique root of the equation $\Psi_{s_b,s_a,q}(\b)=0$ in the interval $(1,s_b+2)$  and $\Psi$ is defined by \eqref{eq:Psi}.
\end{theorem}

The polynomial $\Psi_{x,y,q}$ in~\eqref{eq:Psi} is cubic when $q\neq 0$, and is quadratic when $q=0$. 
The existence and uniqueness of $\psi_q(x,y)$ as a root of $\Psi_{x,y,q}(h)=0$ in the interval $(1,x+2)$ follow from~\cref{lemma:roots-of-psi} in the appendix.
% When $q=0$, however, $\Psi_{x,y,q}$ degenerates to a quadratic equation,
% \[
% \Psi_{x,y,0}(h)=2h^2+(2y-3)h    -2xy-x-4y.
% \]

We denote the concatenation procedures in~\cref{theorem:conpp,theorem:conpf,theorem:congp} by $\ConPP_q$, $\ConPF_q$, and $\ConGP_q$, respectively, each taking $h_a$ and $h_b$ as inputs and returning $h_c$.
For the sake of notational simplicity, we will omit the subscript $q$ of  $\varphi_q,\psi_q,\ConPP_q,\ConPF_q,\ConGP_q$ if there is no ambiguity.

\subsection{Proofs of the Concatenation Theorems}

We use common notation for the certificate differences on subintervals of a gradient descent trajectory, with a fixed reference stationary point $\x_\star$.
For a segment $[h_\ell,\dotsc,h_{r-1}]^\top\in\rk{r-\ell}$, with $0\leq\ell\leq r$, set
\[
s_{\ell,r}\triangleq\sigma_q([h_\ell,\dotsc,h_{r-1}]^\top),\qquad
\omega_i^{\ell,r}\triangleq h_i\prod_{j=\ell}^i(1-q h_j)^{-1}
\quad(\ell\leq i<r).
\]

The primitive certificate difference is
\begin{equation}\label{eq:certificate-primitive}
\begin{aligned}
\mathcal C^{\mathrm p}_{\ell,r}\triangleq{}&
\frac L2\|\x_\ell-\x_\star\|^2
+(1-q)\sum_{i=\ell}^{r-1}\omega_i^{\ell,r}Q_{i,\star}+q s_{\ell,r}Q_{r,\star}-s_{\ell,r}(f_r-f_\star)\\
&-\frac L2\|\x_r-\x_\star\|^2
-\frac{s_{\ell,r}(s_{\ell,r}+1)}{2L}\|\g_r\|^2
-\frac{q s_{\ell,r}}{2L(1-q)}\|\g_r-L(\x_r-\x_\star)\|^2.
\end{aligned}
\end{equation}
For a weight vector $u=(u_i)_{i=\ell}^{r-1}$, set $\eta=1+2s_{\ell,r}+q s_{\ell,r}^2$ and define the \dom certificate difference by
\begin{equation}\label{eq:certificate-dominant}
\begin{aligned}
&\mathcal C^{\mathrm f}_{\ell,r}(u)\triangleq{}
\frac L2\|\x_\ell-\x_\star\|^2+(1-q)\sum_{i=\ell}^{r-1}u_iQ_{i,\star}
+\bigl(\eta-(1-q)\onet u\bigr)Q_{r,\star}-\eta(f_r-f_\star)\\
&-\frac{(1-q)L}{2(1+q s_{\ell,r})^2}
\left\|\x_\ell-\x_\star-\sum_{i=\ell}^{r-1}u_i\left(\frac{\g_i}L-q(\x_i-\x_\star)\right)
-\frac{\eta-(1-q)\onet u}{1-q}\left(\frac{\g_r}L-q(\x_r-\x_\star)\right)\right\|^2.
\end{aligned}
\end{equation}
By~\cref{def:primitive,def:f-dominant},
 $\mathcal C^{\mathrm p}_{\ell,r}\geq0$ when the segment is primitive,
and $\mathcal C^{\mathrm f}_{\ell,r}(u)\geq0$ when it is \dom with certificate $u$.
For $\ell=r$, all sums are empty and $s_{\ell,\ell}=0$.

For parameters $s,\alpha$ and integer endpoint indices $0\leq n<m$, define
\begin{equation}\label{eq:sd-coefficients}
\begin{aligned}
C_n^\alpha&\triangleq\frac{s+3\alpha-2\alpha^2-q\alpha^2(s+2-\alpha)}{q s+s+2-(2q s+q+1)\alpha},\qquad
C_m^\alpha&\triangleq\frac{2s^2+3s+\alpha+q s^2(s+2-\alpha)}{q s+s+2-(2q s+q+1)\alpha},
\end{aligned}
\end{equation}
whenever the common denominator is nonzero. The dependence on $q$ and $s$ is suppressed in this notation.

Before proving the theorems, we present the sufficient-decrease property associated with primitive schedules in the following lemma. This property was established for multiple steps in the smooth convex setting ($q = 0$) in~\cite[Lemma~3.5]{zhang2024accelerated}, and a single-step version with a lower curvature bound was established in~\cite[Lemma~4.3]{rotaru2024exact}. Here we extend this property to multiple steps for all $q<1$. %The algebraic constraints in~\cref{lemma:sd} are automatically satisfied for particular choices of $\alpha$, as discussed in~\cref{lemma:condition-of-sd}.

\begin{lemma}[Sufficient decrease]\label{lemma:sd} 
Let $m,n$ be integers with $0\leq n<m$, and let $h_a\in\R^{n}$ be an arbitrary schedule and $h_b\in\rk{m-n-1}$ be a primitive schedule. Set $s=\sigma(h_b)$ and suppose that $\alpha$ satisfies
\begin{equation}\label{eq:sd-conditions}
\alpha\geq 1,\qquad q\alpha<1,\qquad q s+s+2-(2q s+q+1)\alpha>0.
\end{equation}
Then the iterates generated by gradient descent with the schedule $h_c=[h_a^\top,\alpha,h_b^\top]^\top$ satisfy
\begin{equation*}
f_n-f_m\geq \frac{C_n^\alpha\|\g_n\|^2+C_m^\alpha\|\g_m\|^2}{2L},
\end{equation*}
where $C_n^\alpha,C_m^\alpha$ are given by~\eqref{eq:sd-coefficients} with $s=\sigma(h_b)$.
\end{lemma}

\begin{proof}
For $k\in\{n,m\}$, set
\[
T_k\triangleq-(1-q)\sum_{i=n+1}^{m-1}\omega_i^{n+1,m}Q_{i,k}\geq0.
\]
Taking $l=n+1$, $r=m$, and $u_i=\omega_i^{n+1,m}$ in~\eqref{eq:A4}, and using~\eqref{eq:sum-of-x} and~\eqref{eq:telescope}, we obtain
\[
T_k=\left\langle\g_k-L(\x_k-\x_\star),\x_{n+1}-\x_\star-(1+q s)(\x_m-\x_\star)\right\rangle
-(1-q)\bigl[sQ_{\star,k}+\sum_{i=n+1}^{m-1}\omega_i^{n+1,m}Q_{i,\star}\bigr].
\]

Set
$
B_{n,m}\triangleq\langle\g_n,\x_n-\x_m\rangle-\frac L2\|\x_n-\x_m\|^2$.
Using the definitions of $\mathcal C^{\mathrm p}_{n+1,m}$, $Q_{i,j}$,
the formula for $T_k$ above, and also $\x_{n+1}=\x_n-\a\g_n/L$ gives the following two identities:
\begin{equation}\label{eq:cert-sd-a}
\begin{aligned}
\mathcal C^{\mathrm p}_{n+1,m}+T_n&=(1-q)s(f_n-f_m)
+\frac{\alpha^2-2\alpha-s}{2L}\|\g_n\|^2-\frac{s(s+1)}{2L}\|\g_m\|^2 +(1+qs)B_{n,m},
\end{aligned} 
\end{equation}
and
\begin{equation}\label{eq:cert-sd-b}
\begin{aligned}
&\quad\  \mathcal C^{\mathrm p}_{n+1,m}+T_m
-(\alpha-1)Q_{n,m}-(1-q\a)Q_{m,n} \\ &=\bigl[2-(1+q)\alpha\bigr](f_n-f_m)
+\frac{\alpha(\alpha-1)}{2L}\|\g_n\|^2
-\frac{s(s+2)+\alpha}{2L}\|\g_m\|^2
-(1-q\a)B_{n,m}.
\end{aligned} 
\end{equation}

Computing $(1-q\a)\eqref{eq:cert-sd-a}+(1+qs)\eqref{eq:cert-sd-b}$ cancels the $B_{n,m}$ terms and yields the identity
\begin{equation}\label{eq:cert-sd}
\begin{aligned}
&\qquad (1-q\a)\bigl(\mathcal C^{\mathrm p}_{n+1,m}+T_n\bigr)
+(1+qs)\bigl[\mathcal C^{\mathrm p}_{n+1,m}+T_m-(\alpha-1)Q_{n,m}-(1-q\a)Q_{m,n}\bigr]\\
% \\ &=(q s+s+2-(2q s+q+1)\alpha)
% \left(f_n-f_m-\frac{C_n^\alpha\|\g_n\|^2+C_m^\alpha\|\g_m\|^2}{2L}\right)
% =&\bigl[q s+s+2-(2q s+q+1)\alpha\bigr](f_n-f_m)\\
&= \bigl[q s+s+2-(2q s+q+1)\alpha\bigr](f_n-f_m)\\
&-\frac{s+3\alpha-2\alpha^2-q\alpha^2(s+2-\alpha)}{2L}\|\g_n\|^2
-\frac{2s^2+3s+\alpha+q s^2(s+2-\alpha)}{2L}\|\g_m\|^2.\\
\end{aligned}
\end{equation}
% Write $r\triangleq1-q\alpha$, $t\triangleq1+q s$, and $D\triangleq q s+s+2-(2q s+q+1)\alpha$.
% Since $h_b$ is primitive, $\mathcal C^{\mathrm p}_{n+1,m}\geq0$. 
% We have
% \begin{equation}\label{eq:cert-sd}
% \begin{aligned}
% &D\left(f_n-f_m-\frac{C_n^\alpha\|\g_n\|^2+C_m^\alpha\|\g_m\|^2}{2L}\right)\\
% =&\bigl[q s+s+2-(2q s+q+1)\alpha\bigr](f_n-f_m)\\
% &-\frac{s+3\alpha-2\alpha^2-q\alpha^2(s+2-\alpha)}{2L}\|\g_n\|^2
% -\frac{2s^2+3s+\alpha+q s^2(s+2-\alpha)}{2L}\|\g_m\|^2\\
% =&\bigl[(1-q)sr+\bigl(2-(1+q)\alpha\bigr)t\bigr](f_n-f_m)\\
% &+\frac{r(\alpha^2-2\alpha-s)+t\alpha(\alpha-1)}{2L}\|\g_n\|^2
% -\frac{rs(s+1)+t\bigl(s(s+2)+\alpha\bigr)}{2L}\|\g_m\|^2\\
% =&r\bigl(\mathcal C^{\mathrm p}_{n+1,m}+T_n\bigr)
% +t\bigl[\mathcal C^{\mathrm p}_{n+1,m}+T_m-(\alpha-1)Q_{n,m}-rQ_{m,n}\bigr]\\
% \geq&0,
% \end{aligned}
% \end{equation}
Since $h_b$ is primitive, $\mathcal C^{\mathrm p}_{n+1,m}\geq0$. Also,
$T_k,-Q_{n,m},-Q_{m,n}\geq 0$, $1+qs>0$ by~\cref{lemma:s-bound},
$1-q\a>0$ and $\a-1\geq 0$ by assumption, the left-hand side of~\eqref{eq:cert-sd}
is nonnegative. Thus the right-hand side is nonnegative. Dividing it
by $q s+s+2-(2q s+q+1)\alpha$ and using the definitions of $C_n^\alpha,C_m^\alpha$
gives the desired inequality.
% the second uses $r=1-q\alpha$ and $t=1+q s$,
% the last equality follows from~\eqref{eq:certificate-primitive} and~\eqref{eq:Qij},
% the expression for $T_k$ above, and $\x_{n+1}=\x_n-\alpha\g_n/L$,
% and the inequality holds because $\mathcal C^{\mathrm p}_{n+1,m},T_n,T_m\geq0$, $\alpha\geq1$, $r,t>0$, and $Q_{n,m},Q_{m,n}\leq0$.
% Since $D>0$ by~\eqref{eq:sd-conditions}, we obtain~\eqref{eq:sd}.
\end{proof}

The following lemma verifies~\eqref{eq:sd-conditions} for the concatenation stepsizes.
It also ensures that the concatenated schedule belongs to $\rk{m}$, as defined in~\cref{def:rkn}.
Its proof is given in~\cref{sect:proofs-main-lemmas}.
\begin{lemma}\label{lemma:condition-of-sd}
Let $q<1$ and $s,t\geq0$ satisfy $1+q s>0$ and $1+q t>0$. If $\alpha=\varphi(s,t)$ or $\alpha=\psi(s,t)$, then we have
\[
\alpha>1,\qquad 1-q\alpha>0,\qquad q s+s+2-(2q s+q+1)\alpha>0.
\]
\end{lemma}

We next use \cref{lemma:sd} to prove \cref{theorem:congp}.

\begin{proof}[Proof of~\Cref{theorem:congp}] Let $s_a  = \sigma(h_a)$ and $s_b=\sigma (h_b)$. Let $C_n^\beta,C_m^\beta$ be given by~\eqref{eq:sd-coefficients} with $s=s_b$ and $\alpha=\beta$. When performing gradient descent with $h_c =[h_a^\top,\b,h_b^\top]^\top$, the leading $n$ steps satisfy
\begin{equation}\label{eq:congp-prefix}
f_0-f_n\geq \frac{2s_a+q s_a^2}{2L}\|\g_n\|^2,
\end{equation}
since $h_a$ is \dual by~\cref{def:g-dominant}. By~\cref{lemma:condition-of-sd}, the conditions~\eqref{eq:sd-conditions} hold with $s=s_b$ and $\alpha=\beta$. Thus,~\cref{lemma:sd} gives
\begin{equation}\label{eq:congp-suffix}
f_n-f_m\geq \frac{C_n^\beta\|\g_n\|^2+C_m^\beta\|\g_m\|^2}{2L}.
\end{equation}
Adding~\eqref{eq:congp-prefix} and~\eqref{eq:congp-suffix} yields
\begin{equation}
f_0 - f_m\geq \frac{C_n^\beta+2s_a+q s_a^2}{2L}\|\g_n\|^2+\frac{C_m^\beta}{2L}\|\g_m\|^2. \label{eq:congp-proof}\end{equation}
By~\eqref{eq:recursion-of-sigma}, $s_c\triangleq\sigma(h_c)=\frac{s_a+\b+s_b+q s_as_b}{1-q\b}$.
Using $\Psi_{s_b,s_a,q}(\beta)=0$ where $\Psi$ is defined in~\eqref{eq:Psi}, we obtain
\begin{equation}\label{eq:congp-sb}
s_b = \b- 2+\frac{2(\b-1)^2(qs_a+1)^2}{(qs_a+1)^2(qs_a^2+2s_a+1)-q(qs_a^2+2s_a+\b)^2}.
\end{equation}
The details of~\eqref{eq:congp-sb} are provided in~\cref{lemma:phi-psi-x}.
Substituting $s_c$, $s_b$ and the coefficients in~\eqref{eq:sd-coefficients}, we obtain
\[
C_n^\beta+2s_a+q s_a^2=0,
\qquad C_m^\beta=2s_c+q s_c^2.
\]
Thus,~\eqref{eq:congp-proof} reduces to $f_0-f_m\geq (2s_c+q s_c^2)\|\g_m\|^2/(2L)$, which is~\eqref{eq:g-dominant} for $h_c$. Hence $h_c$ is \dual.
 \end{proof}

We now turn to the proof of \cref{theorem:conpp}.

\begin{proof}[Proof of~\Cref{theorem:conpp}]
By~\cref{lemma:condition-of-sd} and the symmetry of $\varphi$,
the conditions~\eqref{eq:sd-conditions} hold with $s=s_b$.
Thus~\cref{lemma:sd} gives
\begin{equation}\label{eq:conpp-sd}
f_n-f_m-\frac{C_n^\alpha\|\g_n\|^2+C_m^\alpha\|\g_m\|^2}{2L}\geq0,
\end{equation}
where $C_n^\alpha,C_m^\alpha$ are given by~\eqref{eq:sd-coefficients} with $s=s_b$.
Set $s_c=\sigma(h_c)$ and abbreviate
\begin{equation*}
P_1\triangleq\sum_{i=0}^{n-1}\omega_i^{0,n}Q_{i,\star},\qquad
P_2\triangleq\sum_{i=n+1}^{m-1}\omega_i^{n+1,m}Q_{i,\star},\qquad
P_3\triangleq\sum_{i=0}^{m-1}\omega_i^{0,m}Q_{i,\star}.
\end{equation*}
With these notations and the definitions of $Q_{i,j}$, $\mathcal C^{\mathrm p}_{0,m},\mathcal C^{\mathrm p}_{0,n}$ and $\mathcal C^{\mathrm p}_{n+1,m}$ in~\eqref{eq:certificate-primitive}
are expressed as linear combinations of the function values and their corresponding $P_j$,
plus a quadratic form in the endpoint iterates and gradients.

By~\cref{def:primitive} and~\eqref{eq:certificate-primitive}, it remains to show $\mathcal C^{\mathrm p}_{0,m}\geq0$ to prove that $h_c$ is primitive.
Splitting the weighted sum $P_3$ at $i=n$ and using $h_n=\alpha$ gives
\[
P_3=P_1+\alpha\prod_{j=0}^n(1-q h_j)^{-1}Q_{n,\star}
+\prod_{j=0}^{n}(1-q h_j)^{-1}P_2
\stackrel{\eqref{eq:qs-add-1}}{=}
P_1+\frac{\alpha(1+q s_a)}{1-q\alpha}Q_{n,\star}
+\frac{1+q s_a}{1-q\alpha}P_2.
\]
We claim the identity
\begin{equation}\label{eq:cert-p}
\begin{aligned}
\mathcal C^{\mathrm p}_{0,m}
&=\mathcal C^{\mathrm p}_{0,n}
+\frac{1+q s_a}{1-q\alpha}\mathcal C^{\mathrm p}_{n+1,m}+\frac{(1-q)(s_a+\alpha)}{1-q\alpha}
\left(f_n-f_m-\frac{C_n^\alpha\|\g_n\|^2+C_m^\alpha\|\g_m\|^2}{2L}\right).
\end{aligned}
\end{equation}
Since $h_a$ and $h_b$ are primitive,
$\mathcal C^{\mathrm p}_{0,n},\mathcal C^{\mathrm p}_{n+1,m}\geq0$.
The sufficient-decrease inequality~\eqref{eq:conpp-sd} makes the last term nonnegative, and all coefficients are nonnegative by~\cref{lemma:s-bound,lemma:condition-of-sd}.
Hence~\eqref{eq:cert-p} gives $\mathcal C^{\mathrm p}_{0,m}\geq0$.
It remains to verify this identity. We use the substitutions
  \begin{align*}
      \x_{n+1}&=\x_n - \frac{\a}{L}\g_n,\\
s_c& = \sigma (h_c)\overset{\eqref{eq:recursion-of-sigma}}{=}\frac{s_a+\a+s_b+q s_as_b}{1-q\a},\\
      s_a &=\a-2+
     \frac{2(\a-1)^2(qs_b+1)}{ (q s_b+s_b+2)-(2q s_b+q +1)\a},
  \end{align*} 
where the expression for $s_a$ follows from~\eqref{eq:phi-x} with $x=s_a$ and $y=s_b$.
Using these substitutions, the difference between the two sides of~\eqref{eq:cert-p} can be represented as a linear form in $f_n,f_m,f_\star,P_1,P_2$ plus a quadratic form in $\g_n,\g_m,\x_0,\x_n,\x_m,\x_\star$. Each of the coefficients is in the domain $\mathbb Q(q,s_b,\a,L)$.
Expanding both sides and comparing coefficients shows that it is identically zero.
The script for verification of this identity is provided in~\cite{foreverhaibara2026stepsizeconcatenation}.
Thus,~\eqref{eq:cert-p} holds, completing the proof.\end{proof}

Finally, we prove~\cref{theorem:conpf} by constructing a vector $u_c\in\R_+^m$ and verifying that it satisfies~\eqref{eq:dominant}. We provide the explicit construction in~\cref{lemma:dom-u}. The inequalities claimed in~\cref{lemma:dom-u} are proved via purely algebraic computations,
and its proof is deferred to \cref{sect:proofs-main-lemmas}.

\begin{lemma}\label{lemma:dom-u}
    Under the assumptions and notation of~\cref{theorem:conpf}, let $u_b=[u_{n+1},\ldots,u_{m-1}]^\top$ be a certificate for the \dom property of $h_b$, and set $s_c=\sigma(h_c)$. Define
    $u_c\triangleq[u_0,\dotsc,u_{n-1}, \gamma,(v-w)u_b^\top]^\top$ 
    where
    \begin{equation*}
    \begin{aligned}
    u_i &\triangleq \omega_i^{0,n}=h_i\prod_{j=0}^{i}(1-q h_j)^{-1},\qquad (i=0,1,\dotsc,n-1),\\
        \gamma&\triangleq\frac{(s_a+2-q\b)(q s_a+1)}{(1-q)(1-q\b)},\qquad
         v\triangleq\frac{(s_a+1)^2}{(\b-1)^2},\qquad
         w\triangleq \frac{(s_a+1)(q s_a+2-q\b)}{ (1-q\b)(1+2s_b+q s_b^2)}.
    \end{aligned}
    \end{equation*}
Then $\gamma,v,w\geq0$, $v-w\geq0$, and $u_c\in\R_+^m$. Moreover,
    \[
    (1+2s_c+q s_c^2)-(1-q)(\onet u_c)\geq0.
    \]
\end{lemma}

\begin{proof}[Proof of \Cref{theorem:conpf}] Let \(u_b\) be a certificate for the \(f\)-dominant property of \(h_b\), and define \(u_c\) as in~\cref{lemma:dom-u}.
Set
\[
s_c\triangleq\sigma(h_c)
\stackrel{\eqref{eq:recursion-of-sigma}}{=}
\frac{s_a+s_b+\beta+q s_as_b}{1-q\beta},
\qquad
\eta_j\triangleq1+2s_j+q s_j^2\quad(j=b,c).
\]
For $h_a$ and $h_b$, respectively, define
\begin{align*}
P_1&\triangleq\sum_{i=0}^{n-1}\omega_i^{0,n}Q_{i,\star},&
R_1&\triangleq\sum_{i=0}^{n-1}\omega_i^{0,n}\left(\frac{\g_i}L-q(\x_i-\x_\star)\right),\\
P_2&\triangleq\sum_{i=n+1}^{m-1}u_iQ_{i,\star},&
R_2&\triangleq\sum_{i=n+1}^{m-1}u_i\left(\frac{\g_i}L-q(\x_i-\x_\star)\right).
\end{align*}
The corresponding quantities for $h_c$ are
\begin{equation*}
\begin{aligned}
P_3&\triangleq\sum_{i=0}^{m-1}(u_c)_iQ_{i,\star}=P_1+\gamma Q_{n,\star}+(v-w)P_2,\\
R_3&\triangleq\sum_{i=0}^{m-1}(u_c)_i\left(\frac{\g_i}L-q(\x_i-\x_\star)\right)=R_1+\gamma\left(\frac{\g_n}L-q(\x_n-\x_\star)\right)+(v-w)R_2.
\end{aligned}
\end{equation*}
With these notations and the definitions of $Q_{i,j}$, certificates $\mathcal C^{\mathrm f}_{0,m}(u_c),\mathcal C^{\mathrm p}_{0,n}$ and $\mathcal C^{\mathrm f}_{n+1,m}(u_b)$
from~\eqref{eq:certificate-primitive} and~\eqref{eq:certificate-dominant} are expressed as
linear combinations of the function values and their corresponding $P_j$,
plus a quadratic form in the endpoint iterates, gradients and $R_j$.

Using~\eqref{eq:A4}, we define $T$ by
\[
T\triangleq -(1-q)\sum_{i=n+1}^{m-1}u_iQ_{i,n}
=-(1-q)\bigl[(\onet u_b)Q_{\star,n}+P_2\bigr]+\langle\g_n-L(\x_n-\x_\star),R_2\rangle.
\]
Then $T\geq 0$ because $q<1$, $u_b\in\R_+^{m-n-1}$ and $Q_{i,n}\leq 0$.

% The construction of $u_c$ gives $\onet u_c=s_a+\gamma+(v-w)(\onet u_b)$.
By~\cref{lemma:dom-u}, $u_c\in\R_+^m$ and $\eta_c-(1-q)\onet u_c\geq0$.
Thus, by~\cref{def:f-dominant} and~\eqref{eq:certificate-dominant}, it remains to prove $\mathcal C^{\mathrm f}_{0,m}(u_c)\geq0$.
We claim the following identity:
\begin{equation}\label{eq:cert3}
\begin{aligned}
\mathcal C^{\mathrm f}_{0,m}(u_c)
={}&\mathcal C^{\mathrm p}_{0,n}+v\mathcal C^{\mathrm f}_{n+1,m}(u_b)
% -\underbrace{w\left[(1-q)\sum_{i=n+1}^{m-1}u_iQ_{i,n}
% +\bigl(\eta_b-(1-q)\onet u_b\bigr)Q_{m,n}\right]}_{\triangleq\,\mathcal I}.
+wT-w\bigl(\eta_b-(1-q)\onet u_b\bigr)Q_{m,n}.
\end{aligned}
\end{equation}

Here $\mathcal C^{\mathrm p}_{0,n},\mathcal C^{\mathrm f}_{n+1,m}(u_b)\geq0$ by the assumptions,
$v,w\geq0$ by~\cref{lemma:dom-u} and
$\eta_b-(1-q)\onet u_b\geq 0$ from~\cref{def:f-dominant}. % Also, $T\geq 0$ by $u_c\in\R_+^m$ and $Q_{i,\star}\leq 0$.
% and $\mathcal I\leq0$ because the bracketed sum is nonpositive.
Thus,~\eqref{eq:cert3} gives $\mathcal C^{\mathrm f}_{0,m}(u_c)\geq0$, as desired.
To verify~\eqref{eq:cert3}, we use the substitutions
\begin{subequations}\label{eq:cert3-params}
    \begin{align}
        \x_{n+1}& =\x_n-\frac{\b}{L}\g_n,
\\        \onet u_c&=s_a+\gamma+(v-w)(\onet u_b),\label{eq:confp-onet-uc}
    % \\ \gamma &=\frac{s_a+2-q\b}{1+q (s_a+1-\b)},
    % \\ s_a & =\frac{\beta(q\beta^2-2(q+1)\beta+3)-s_b(q \b+\b-2)(q s_b+2)}{q \b^2-1+(2q\b-q -1)s_b(q s_b+2)}\label{eq:sa},
    \\ s_a & =\b- 2+\frac{2(\b -1)^2(q s_b+1)^2}{(qs_b+1)^2(qs_b^2+2s_b+1)-q(qs_b^2+2s_b+\b)^2}\label{eq:sa},
    \\ R_1&\stackrel{\eqref{eq:sum-of-x}}{=}(\x_0-\x_\star)-(q s_a+1)(\x_n-\x_\star),
    \end{align}
\end{subequations}
where~\eqref{eq:confp-onet-uc} follows from the construction of $u_c$,
and~\eqref{eq:sa} follows by solving $\Psi_{s_a,s_b,q}(\beta)=0$ for $s_a$, as given by~\eqref{eq:psi-x} in~\cref{lemma:phi-psi-x}.
% Using~\eqref{eq:A4}, we first replace the weighted sum in $\mathcal I$ by
% \[
% (1-q)\sum_{i=n+1}^{m-1}u_iQ_{i,n}
% =(1-q)\bigl[(\onet u_b)Q_{\star,n}+P_2\bigr]-\langle\g_n-L(\x_n-\x_\star),R_2\rangle.
% \]
We then expand $\mathcal C^{\mathrm f}_{0,m}(u_c)$, $\mathcal C^{\mathrm p}_{0,n}$, and $\mathcal C^{\mathrm f}_{n+1,m}(u_b)$ using~\eqref{eq:certificate-primitive} and~\eqref{eq:certificate-dominant}.
The definitions above and~\eqref{eq:cert3-params} eliminate
$\eta_c,\allowbreak s_c,\allowbreak\onet u_c,\allowbreak T,\allowbreak P_3,\allowbreak R_3,\allowbreak\eta_b,\allowbreak v,\allowbreak w,\allowbreak\x_{n+1},\allowbreak\gamma,\allowbreak s_a,\allowbreak R_1$.
The difference between the two sides of~\eqref{eq:cert3} is then a linear form in $f_n,f_m,f_\star,P_1,P_2$ plus a quadratic form in
$\g_n,\allowbreak\g_m,\allowbreak\x_0,\allowbreak\x_n,\allowbreak\x_m,\allowbreak\x_\star,\allowbreak R_2$.
All coefficients lie in $\mathbb Q(q,s_b,\beta,L,\onet u_b)$ and vanish after expansion.
The script for verification of this identity is provided in~\cite{foreverhaibara2026stepsizeconcatenation}.
This verifies~\eqref{eq:cert3} and completes the proof.\end{proof}

\section{Optimized Schedules and Convergence Analysis}
\label{sect:convergence}
We now apply the concatenation theorems to construct optimized schedules. We first give the dynamic-programming construction, then establish sharp finite-step guarantees for the function-value gap and gradient norm, and finally derive growth estimates for $\sigma_q$ and the resulting iteration complexity.

\subsection{Optimized Schedules via Dynamic Programming}
\label{sect:dp-construction}
As an immediate consequence of~\cref{theorem:conpp,theorem:conpf,theorem:congp}, a schedule of $n$ steps can be formed by concatenating schedules of $k$ and $n-k-1$ steps. Since the bounds in~\cref{theorem:bound-of-fdom,theorem:bound-of-dual} improve as $\sigma_q(h)$ increases, and since $\varphi$, $\psi$, and $\sigma_q$ are monotone in the relevant arguments by~\cref{lemma:s-bound,lemma:monotonicity-of-phi-psi}, it is sufficient to retain, for each iteration budget, a schedule with the largest value of $\sigma_q$. This principle gives the optimized \ConPP{} and \ConPF{} schedules below. Reversing the latter yields the optimized \ConGP{} schedules.

\Needspace{8\baselineskip}
\begin{definition}\label{def:dynamic-programming}
We recursively define three families of schedules $\{h_\circ^{(n)}\}_{n=0}^\infty$, $\{h_\bullet^{(n)}\}_{n=0}^\infty$, and $\{h_\algC^{(n)}\}_{n=0}^\infty$ as follows.
\begin{enumerate}[label=(\roman*),ref=(\roman*),leftmargin=*,itemsep=0.3em]
\item \textbf{Optimized \ConPP{} schedules (\OConPP).} Set $h_\circ^{(0)}=[~]\in\rk{0}$. For $n\geq1$, define
\[
\begin{aligned}
\mathcal H_\circ^{(n)}&=\{\ConPP(h_\circ^{(k)},h_\circ^{(n-k-1)}):\ k\in\{0,1,\dotsc,n-1\}\},\qquad
h_\circ^{(n)}&\in \argmax_{h\in\mathcal H_\circ^{(n)}} \{\sigma_q(h)\}.
\end{aligned}
\]
\item \textbf{Optimized \ConPF{} schedules (\OConPF).} Set $h_\bullet^{(0)}=[~]\in\rk{0}$. For $n\geq1$, define
\[
\begin{aligned}
\mathcal H_\bullet^{(n)}&=\{\ConPF(h_\circ^{(k)},h_\bullet^{(n-k-1)}):\ k\in\{0,1,\dotsc,n-1\}\},\qquad
h_\bullet^{(n)}&\in \argmax_{h\in\mathcal H_\bullet^{(n)}} \{\sigma_q(h)\}.
\end{aligned}
\]
\item \textbf{Optimized \ConGP{} schedules (\OConGP).} For $n\geq0$, define $h_\algC^{(n)}=\rev {h_\algB^{(n)}}$, where $\rev h$ denotes the vector obtained by reversing the order of the entries. %For example, $\rev {[1,2]^\top}=[2,1]^\top$.
\end{enumerate}
\end{definition}

The dynamic-programming construction follows the approach of our earlier work~\cite{zhang2024accelerated}, using the curvature-dependent concatenation rules established above. Following the terminology of~\cite{grimmer2024composing}, we call these schedules ``optimized'' because they maximize $\sigma_q$ at each prescribed length within the corresponding class generated recursively from the empty schedule. We next verify that the three families of schedules are primitive, \dom, and \dual, respectively.

% The following are direct consequences of~\cref{theorem:conpp} and~\cref{theorem:congp}.

\begin{proposition}\label{prop:dp-are-primitive-dom} For every $n\in\N$, the schedules $h_\algA^{(n)}$ defined in~\cref{def:dynamic-programming} and their reversals $\rev h_\algA^{(n)}$ are primitive, the schedules $h_\algB^{(n)}$ are \dom, and the schedules $h_\algC^{(n)}$ are \dual.
Moreover, $\sigma(h_\algB^{(n)})=\sigma(h_\algC^{(n)})$.
\end{proposition}
\begin{proof} Since $h_\algA^{(0)}=h_\algB^{(0)}=[~]$ is primitive and \dom by~\cref{coro:zero-step}, each $h_\algA^{(n)}$ is primitive and $h_\algB^{(n)}$ is \dom 
by induction using~Theorems \ref{theorem:conpp} and~\ref{theorem:conpf}. For the remaining claims, we adapt the reversal argument in~\cite[Proposition 5.5]{zhang2024accelerated}.

We next show by induction that $\rev h_\algA^{(n)}$ is primitive for all $n\geq0$. The result holds for $n=0$. Assume it holds for all $n< N$ where $N\geq 1$. For $n=N$, since there exists an integer $k$ such that $0\leq k<N$  and $h_\algA^{(N)}=\ConPP(h_\algA^{(k)},h_\algA^{(N-k-1)})$, we have that
\begin{equation}\label{eq:induction-prim}
\begin{aligned}
\rev h_\algA^{(N)}
&=\rev\ConPP(h_\algA^{(k)},h_\algA^{(N-k-1)})\\
&=\left[(\rev h_\algA^{(N-k-1)})^\top,
\varphi\bigl(\sigma(h_\algA^{(k)}),\sigma(h_\algA^{(N-k-1)})\bigr),
(\rev h_\algA^{(k)})^\top\right]^\top\\
&=\ConPP(\rev h_\algA^{(N-k-1)},\rev h_\algA^{(k)}),
\end{aligned}
\end{equation}
where the last equality follows from $\sigma(h)=\sigma(\rev h)$ in~\eqref{eq:sigma} and the symmetry of $\varphi$ in~\eqref{eq:phi}.
By the induction hypothesis, $\rev h_\algA^{(N-k-1)}$ and $\rev h_\algA^{(k)}$ are primitive. Hence the last  $\ConPP$ operator in~\eqref{eq:induction-prim} is well-defined, and we conclude $\rev h_\algA^{(N)}$ is primitive by~\cref{theorem:conpp}.

Using this fact, we next prove that $h_\algC^{(n)}$ is \dual by induction.
When $n=0$, $h_\algC^{(0)}$ is \dual by~\cref{coro:zero-step}. Now we suppose $h_\algC^{(n)}$ is \dual for all $n<N$ where $N\geq 1$. Since  $\rev {h_\algC^{(N)}}=h_\algB^{(N)}\in\mathcal H_\algB^{(N)}$, there exists an integer $k$ such that $0\leq k<N$ and 
\begin{equation}\label{eq:induction-dual}
\begin{aligned}
h_\algC^{(N)}
&=\rev\ConPF(h_\algA^{(k)},h_\algB^{(N-k-1)})\\
&=\left[(h_\algC^{(N-k-1)})^\top,
\psi\bigl(\sigma(\rev h_\algA^{(k)}),\sigma(h_\algC^{(N-k-1)})\bigr),
(\rev h_\algA^{(k)})^\top\right]^\top\\
&=\ConGP(h_\algC^{(N-k-1)},\rev h_\algA^{(k)}),
\end{aligned}
\end{equation}
where the second equality follows from $h_\algC^{(N-k-1)}=\rev h_\algB^{(N-k-1)}$ in~\cref{def:dynamic-programming} and $\sigma(h)=\sigma(\rev h)$ in~\eqref{eq:sigma}.
Recalling that $h_\algC^{(N-k-1)}$ is \dual and $\rev h_\algA^{(k)}$ is primitive, the last $\ConGP$ operator in~\eqref{eq:induction-dual} is well-defined and $h_\algC^{(N)}$ is \dual by~\cref{theorem:congp}. The equality follows from $h_\algC^{(n)}=\rev h_\algB^{(n)}$ and~\eqref{eq:sigma}.\end{proof}
The identities~\eqref{eq:induction-prim} and~\eqref{eq:induction-dual}
hold for any choice of splits throughout the recursive construction.
They therefore establish a one-to-one correspondence under reversal
between the classes generated by $\ConPF$ and $\ConGP$,
with primitive components generated by $\ConPP$.
Since reversal preserves both the length and $\sigma$,
the maximizing schedules $h_\algB^{(n)}$ in the former class
yield maximizing schedules $h_\algC^{(n)}$ in the latter.
This justifies calling $h_\algC^{(n)}$ ``\OConGP{} schedules.''

We compute the schedules in~\cref{def:dynamic-programming} by dynamic programming. Algorithms~\ref{alg1} and~\ref{alg2} store the schedules $H_\algA[k]$, $H_\algB[k]$ and their $\sigma$-values $S_\algA[k]=\sigma(H_\algA[k])$ and $S_\algB[k]=\sigma(H_\algB[k])$. Using these stored values, the recursion~\eqref{eq:recursion-of-sigma} gives the $\sigma$-value of each candidate concatenation. For $n=1,\dotsc,N$, the algorithms compare all splits and construct a schedule of length $n$ with the largest $\sigma$-value. Evaluating all splits and constructing the schedule takes $O(n)$ time for each $n$, giving $O(N^2)$ total time and storage over all lengths up to $N$. This yields the following proposition.

\begin{proposition}
The schedules $h_\algA^{(n)}$ and $h_\algB^{(n)}$ for $n=0,1,\dotsc,N$ can be computed in $O(N^2)$ time using $O(N^2)$ space by Algorithms~\ref{alg1} and~\ref{alg2}.
\end{proposition}  

\begin{center}
\begin{minipage}{0.49\textwidth}
\begin{algorithm} [H]\small
\caption{\texttt{priDP}($N$,$q$)} \label{alg1}
\begin{algorithmic}[1]
\State   $H_\algA[0] \gets  [~]$
\State $S_\algA[0]\gets 0$
\For{$n \gets 1$ to $N$}
    \State $k_{\text{best}}\gets 0$, $s_{\text{max}} \gets 0$
    \For{$k \gets 0$ to $n-1$}
        \State $k'\gets n-k-1$
        \State $s\gets  \frac{ S_\algA [k]+\varphi(S_\algA[k],S_\algA[k'])+S_\algA[k']+qS_\algA[k]S_\algA[k']}{1-q\varphi(S_\algA[k],S_\algA[k'])} $\label{alg:line-7}
        % \State  $s\gets \frac{1}{q}\cdot\left({\frac{(q S_\algA [k]+1)(q S_\algA[k']+1)}{1-q\varphi(S_\algA[k],S_\algA[k'])}-1}\right)$\label{alg:line-7}
        \If{$s \geqslant s_{\text{max}}$}
            \State $k_{\text{best}} \gets k$, $s_{\text{max}} \gets s$
        \EndIf
    \EndFor
    \State $k\gets k_{\text{best}}$, $k'\gets n-k_{\text{best}}-1$
    \State $H_\algA[n] \gets  \ConPP({H_\algA[k], H_\algA[k']})$
    \State $S_\algA[n]\gets s_{\text{max}}$
\EndFor
\State \textbf{return} $H_\algA$, $S_\algA$
%\EndFunction
\end{algorithmic}
\end{algorithm}
\end{minipage}
\hfill
\begin{minipage}{0.499\textwidth}
\begin{algorithm}[H] \small
\caption{{\texttt{domDP}($N$,$q$)}}    \label{alg2}
\begin{algorithmic}[1] 
\State   $H_\algB[0] \gets [~]$, $S_\algB[0]\gets 0$ \State $H_\algA,S_\algA\gets\Call{priDP}{N,q}$ \label{line:invoke-alg1}
\For{$n \gets 1$ to $N$}
    \State $k_{\text{best}}\gets 0$, $s_{\text{max}} \gets 0$
    \For{$k \gets 0$ to $n-1$}        
        \State $k'\gets n-k-1$
        % \State $s\gets S_\algA[k]+\psi(S_\algA[k],S_\algB[k'])+S_\algB[k']$
        % \State $s\gets  \frac{1}{q}\cdot\left({\frac{(q S_\algA [k]+1)(q S_\algB[k']+1)}{1-q\psi(S_\algA[k],S_\algB[k'])}-1}\right)$
        \State $s\gets  \frac{ S_\algA [k]+\psi(S_\algA[k],S_\algB[k'])+S_\algB[k']+qS_\algA[k]S_\algB[k']}{1-q\psi(S_\algA[k],S_\algB[k'])} $
        \If{$s \geqslant s_{\text{max}}$}
            \State $k_{\text{best}} \gets k$, $s_{\text{max}} \gets s$
        \EndIf
    \EndFor 
    \State $k\gets k_{\text{best}}$, $k'\gets n-k_{\text{best}}-1$
    % \State $H_\algB[n] \gets  [{H_\algA[k]^\top,\psi(S_\algA[k],S_\algB[k']), H_\algB[k']^\top}]^\top$
    \State $H_\algB[n]\gets \ConPF({H_\algA[k],H_\algB[k']})$
    \State $S_\algB[n]\gets s_{\text{max}}$
\EndFor
\State \textbf{return} $H_\algB$, $S_\algB$
% \EndFunction
\end{algorithmic}
\end{algorithm}
\end{minipage}
\end{center}

\begin{remark}[Optimality of balanced splits]
    Algorithms~\ref{alg1} and~\ref{alg2} also generalize to the smooth convex case, $q=0$, where the $s$-update in  
    Line~\ref{alg:line-7} is interpreted as the value when $q=0$,
\[
\lim_{q\to 0}\frac{x+\a+y+q xy}{1-q\a}=x+\alpha+y
% \left.\frac{x+\a+y+q xy}{1-q\a}\right|_{q=0}=x+\alpha+y
\quad\text{($\alpha$ from~\eqref{eq:phi} or~\eqref{eq:Psi})}.
% \lim_{q\to 0_+}\frac{1}{q}\left(\frac{(q x+1)(q y+1)}{1-q \alpha}-1\right)=x+\alpha+y.
\]
Thus in the smooth convex case, the algorithms reduce to Algorithms 4.1 and 4.2 in~\cite{zhang2024accelerated}.
Liu et al.~\cite{liu2026optimalrecursivecompositiondyadic} proved our earlier conjecture on optimal splits for primitive concatenation in the smooth convex case~\cite[Conjecture~4.17]{zhang2024accelerated}. In particular, for $q=0$, the set $\mathcal H_\algA^{(n)}$ in~\cref{def:dynamic-programming} satisfies
\[
\ConPP_0(h_\algA^{\lfloor(n-1)/2\rfloor},h_\algA^{\lceil(n-1)/2\rceil})\in\argmax_{h\in\mathcal H_\algA^{(n)}}(\sigma_0(h)),\qquad\forall n\in\N_+.
\]
This balanced-split property does not hold in general in the smooth strongly convex case.
For example, when $q=1/4$ and $n=5$, direct evaluation of the recursion gives
\[
\begin{aligned}
\sigma_{1/4}\!\left(\ConPP_{1/4}(h_\algA^{(2)},h_\algA^{(2)})\right)&=20+12\sqrt 3<41,\\ 
\sigma_{1/4}\!\left(\ConPP_{1/4}(h_\algA^{(1)},h_\algA^{(3)})\right)&=\frac{19+9\sqrt 5+3\sqrt{102+46\sqrt 5}}{2}>41.
\end{aligned}
\]
Thus, the balanced split $(2,2)$ is suboptimal. For primitive lengths $n=2^k-1$, $k\geq1$, however,
the symmetric split remains uniquely optimal in the \OConPP{} recursion for every $0<q<1$, as shown in~\cref{prop:prim-asymp}.
\end{remark}

Combining~\cref{prop:dp-are-primitive-dom} with the performance guarantees in Section~\ref{sect:stepsize-shcedules} gives the following sharp finite-step bounds.

\begin{theorem}[Sharp finite-step bounds]\label{theorem:general-bound} Let $q=\mu/L<1$ and $n\in\N$. Gradient descent with the schedule $h=h_\algA^{(n)}$ or $h=h_\algB^{(n)}$ satisfies
\[
f_n - f_\star \leq \frac{1}{1+2\sigma(h)+q\sigma(h)^2 }\cdot \frac{L}{2}\|\x_0 - \x_\star\|^2. 
\] 
Gradient descent with the schedule $h=h_\algC^{(n)}$ satisfies
\[
\frac{1}{2L}\|\g_n\|^2\leq  \frac{1}{1+2\sigma(h)+q\sigma(h)^2 }\cdot(f_0-f_\star)
\]
if $f_\star$ is the global infimum. 
Each bound is sharp for the corresponding schedule, with equality attained by the functions in~\cref{theorem:bound-of-fdom,theorem:bound-of-dual}, respectively.
\end{theorem}
\begin{proof} Since~\Cref{prop:dp-are-primitive-dom} asserts that $h_\algA^{(n)}$, $h_\algB^{(n)}$, $h_\algC^{(n)}$ are primitive, \dom and \dual, respectively, the bounds are direct consequences of~\cref{theorem:bound-of-fdom} and~\cref{theorem:bound-of-dual} (using~\cref{prop:primitive-is-fdom} for $h_\algA^{(n)}$). Sharpness follows from the Huber-type examples in the proofs of
\cref{theorem:bound-of-fdom,theorem:bound-of-dual}.\end{proof} 

\subsection{Growth Bounds and Iteration Complexity}
\label{sect:growth-complexity}
We now analyze the case of smooth strongly convex functions where $0<q<1$.
\Cref{theorem:general-bound} reduces the performance analysis to estimating $\sigma_q(h)$. We first establish its growth for the \OConPP{} schedules, bound the exponent $\nu(q)$, and compare the cumulative quantities of the three families.

\begin{proposition}\label{prop:prim-asymp}
Let $q=\mu/L\in(0,1)$, and define the sequence $\{v_n\}_{n\in\N}$ by
\[
v_0=\frac{1+q}{1-q}>1,\qquad
v_{n+1}=\frac{v_n^2+1+(v_n-1)\sqrt{v_n^2+2v_n+5}}{2},\qquad n\in\N.
\]
Then the limit $\nu(q)\triangleq\lim_{n\to\infty}2^{-n}\log v_n$ exists, and for every $n\in\N$,
\begin{equation}
\label{eq:prim-asymp}
\sigma(h_\algA^{(n)})=\frac{(1-q)\exp\left\{(n+1)\cdot \nu(q)+R_n^\circ\right\}- (1+q)}{2q},
\end{equation}
where the remainder $R_n^\circ$ is bounded by
% \[
% R_n^\circ\triangleq
% \log\left(\frac{2q\sigma(h_\algA^{(n)})+1+q}{1-q}\right)
% -(n+1)\nu(q)
% \]
\begin{equation}\label{eq:prim-remainder-bound}
-\frac{1}{2}<-\frac{(1-q)^2}{2(1+q)}\leq R_n^\circ\triangleq
\log\left(\frac{2q\sigma(h_\algA^{(n)})+1+q}{1-q}\right)
-(n+1)\nu(q)\leq0
\end{equation}
and
\[
\lim_{l\to\infty}R_{2^{l}-1}^\circ=0.
\]
Moreover, for any $l\in\N_+$, $h_\algA^{(2^{l}-1)}=\ConPP(h_\algA^{(2^{l-1}-1)},h_\algA^{(2^{l-1}-1)})$, and
the unique optimal split of \OConPP{} in~\cref{def:dynamic-programming} at $n=2^{l}-1$ is $k=2^{l-1}-1$.
\end{proposition}
\begin{proof}
    Let $h_c=[h_a^\top,\a,h_b^\top]^\top$, $s_a=\sigma(h_a)$, $s_b=\sigma(h_b)$ and $s_c=\sigma(h_c)$. \cref{lemma:recursion-of-sigma} implies
\begin{equation}\label{eq:concat-rate}
s_{c}= %\frac{s_a+\a+s_b+qs_as_b}{1-q\a}=
 \frac{(q s_a+1)(1-q\a)^{-1}(q s_b+1)- 1}{q},
\end{equation}
where $\alpha=\varphi(s_a,s_b)$ and $\varphi(\cdot,\cdot)$ is defined in~\eqref{eq:phi}. 

Set $t_i=(q s_i+1)/(1-q)>1$ for $i=a,b,c$, where the inequality follows from $s_i\geq0$ and $0<q<1$.
Solving this definition for $s_i$ gives $s_i=((1-q)t_i-1)/q$.
Substituting these expressions into $\alpha=\varphi(s_a,s_b)$,
we obtain
\begin{equation}\label{eq:concat-rate-transform-alpha}
\begin{aligned}
   \frac{1-q}{1-q\alpha}
   & =\frac{1-q}{1-q\varphi\left(\frac{(1-q)t_a-1}{q},\frac{(1-q)t_b-1}{q}\right)}
  \\ &\overset{\mathclap{\eqref{eq:phi}}}{=}\frac{2}{2t_at_b-t_a-t_b+2-\sqrt{(2t_at_b-t_a-t_b) ^2+{4(t_a-1)(t_b-1)}}}
  \\ &=  \frac{2}{2t_at_b-t_a-t_b+2-\sqrt{(2t_at_b-t_a-t_b+2) ^2-4t_at_b}}  
  \\ &= \frac{2t_at_b-t_a-t_b+2+\sqrt{(2t_at_b-t_a-t_b+2) ^2-4t_at_b}}{2t_at_b}.
\end{aligned}
\end{equation}
Combining~\eqref{eq:concat-rate} and~\eqref{eq:concat-rate-transform-alpha}, we deduce that $t_c$ has the recursion
\begin{equation}\label{eq:concat-rate-v}
 \begin{aligned}
t_c& \overset{\mathclap{\eqref{eq:concat-rate}}}{=}\frac{(q s_a+1)(q s_b+1) (1-q\a)^{-1} }{1-q}=\frac{  (1-q )t_at_b  }{1-q\a}
\\ &\overset{\mathclap{\eqref{eq:concat-rate-transform-alpha}}}{=} \frac{2t_at_b-t_a-t_b+2+\sqrt{(2t_at_b-t_a-t_b+2) ^2-4t_at_b}}{2 }.
 \end{aligned}
\end{equation}

Next, set $w_i=2t_i-1>1$ for $i=a,b,c$.
Substituting $t_i=(w_i+1)/2$ into~\eqref{eq:concat-rate-v} gives
\begin{equation}
    \label{eq:concat-rate-eqw}
\begin{aligned}
w_c=\frac{2q\sigma_q(h_c)+1+q}{1-q}=F(w_a,w_b),\qquad
F(x,y)\triangleq\frac{xy+1+\sqrt{(xy+3)^2-4(x+1)(y+1)}}{2}.
\end{aligned}
\end{equation}

For the optimized schedules $h_\circ^{(n)}$ in~\cref{def:dynamic-programming}, define
\begin{equation}\label{eq:wn-definition}
w_n=\frac{2q\sigma_q(h_\circ^{(n)})+1+q}{1-q},\qquad n\geq0.
\end{equation}
For each $k=0,\dotsc,n-1$, take $h_a=h_\circ^{(k)}$ and $h_b=h_\circ^{(n-k-1)}$.
Then $w_a=w_k$ and $w_b=w_{n-k-1}$, so~\eqref{eq:concat-rate-eqw} gives
\[
\frac{2q\sigma_q\bigl(\ConPP(h_\circ^{(k)},h_\circ^{(n-k-1)})\bigr)+1+q}{1-q}
=F(w_k,w_{n-k-1}).
\]
Since the transformation in~\eqref{eq:wn-definition} is strictly increasing in $\sigma_q$, the maximization in~\cref{def:dynamic-programming} gives
\[
w_{n} = \max_{0\leq k\leq n-1}F(w_k,w_{n-k-1}),\qquad n\geq 1.
\]

By~\eqref{eq:wn-definition}, $w_0=v_0>1$.
Applying~\cref{lemma:bound-Rn} to $a_m=w_{m-1}$, $m\geq1$, with $E_{n+1}=R_n^\circ$, we obtain
\begin{equation*}
w_n=e^{(n+1)\nu(q)+R_n^\circ},\qquad -\frac{1}{w_0(w_0+1)}\leq R_n^\circ\leq0,\qquad \lim_{l\to\infty}R_{2^l-1}^\circ=0.
\end{equation*}
Substituting this expression into~\eqref{eq:wn-definition} gives~\eqref{eq:prim-asymp}.
Since $\sigma_q(h_\circ^{(0)})=0$, \eqref{eq:wn-definition} gives $w_0=(1+q)/(1-q)$, so the bound~\eqref{eq:prim-remainder-bound} follows from
$
\frac{1}{w_0(w_0+1)}=\frac{(1-q)^2}{2(1+q)}<\frac12.
$

The unique optimal split at $n=2^l-1$ follows from
the last assertion of~\cref{lemma:bound-Rn},
with $a_m=w_{m-1}$.
\end{proof}

% \begin{remark}[Remainder at lengths $2^l-1$]\label{rem:dyadic-remainder}
% For each fixed $q\in(0,1)$, the upper bound $R_n^\circ\leq0$ in~\eqref{eq:prim-remainder-bound} is asymptotically attained along $n=2^l-1$ as $l\to\infty$.
% To see this, let $w_n$ be defined by~\eqref{eq:wn-definition}.
% The optimal symmetric splits in~\cref{prop:prim-asymp} and $w_0=v_0$ give $w_{2^l-1}=v_l$ for every $l\in\N$.
% The sequence $\{v_l\}_{l\in\N}$ satisfies the recursion in~\cref{lemma:duplicate} with initial value $v_0=(1+q)/(1-q)>1$, and the limit $C$ there equals $\nu(q)$.
% Item~\ref{item:duplicate-sum} therefore implies that the series $\sum_{l=0}^\infty(\log v_l-2^l\nu(q))$ converges. Consequently,
% \[
% R_{2^l-1}^\circ\;\overset{\mathclap{\eqref{eq:wn-remainder}}}{=}\;\log w_{2^l-1}-2^l\nu(q)=\log v_l-2^l\nu(q)\longrightarrow0\qquad\text{as }l\to\infty.
% \]

% % Moved to the proposition above
% \end{remark}

\begin{figure}[H]
    \centering
\includegraphics[width=\linewidth]{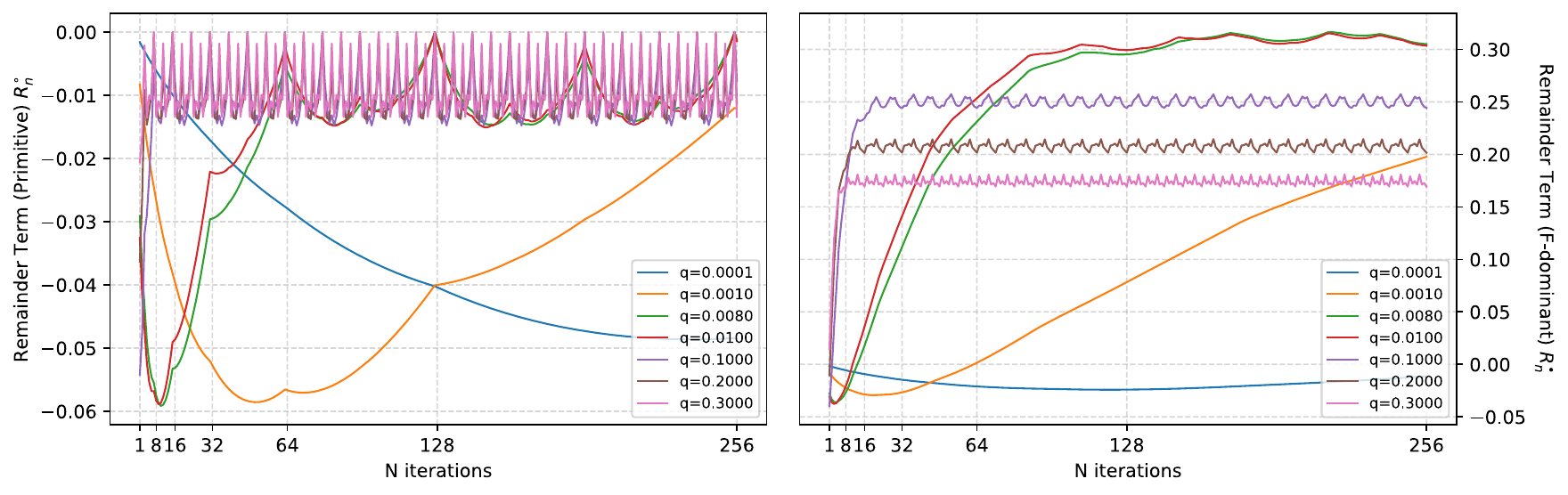}
\caption{The remainders $R_n^\circ$ for the \OConPP{} schedules $h_\circ^{(n)}$ (left) and $R_n^\bullet$ for the \OConPF{} schedules $h_\bullet^{(n)}$ (right).}
\label{fig:remainder-term}
\end{figure}

For the \OConPF{} schedules, define the analogous remainder (see~\cref{rem:improving-constants})
\begin{equation}\label{eq:fdom-remainder}
R_n^\bullet\triangleq
\log\left(\frac{2q\sigma(h_\algB^{(n)})+1+q}{1-q}\right)
-(n+1)\nu(q),\qquad n\geq0.
\end{equation}
Figure~\ref{fig:remainder-term} plots $R_n^\circ$ for the \OConPP{} schedules on the left and $R_n^\bullet$ for the \OConPF{} schedules on the right.
In the left panel, $R_n^\circ$ stays above $-0.06$ in all tested cases, suggesting that the lower bound in~\eqref{eq:prim-remainder-bound} may be improved.
A tighter bound would improve the constant factor in~\eqref{eq:prim-asymp} without changing the exponential  factor $e^{(n+1)\nu(q)}$.

The following proposition gives explicit bounds on $\nu(q)$
for deriving nonasymptotic iteration-complexity bounds. We recall $\rho=\log_{\sqrt 2+1}2\approx 0.7864$.
 \begin{proposition}\label{prop:nu-bound} Let $\nu(q)$ be defined as in~\cref{prop:prim-asymp}. Then, for every $q\in(0,1)$,
\[
\log\frac{1+q}{1-q}\leq\nu(q)\leq \log\frac{1+q}{1-q}+\log\frac{27}{25}.
\]
Moreover, if $q \in (0,\frac{3}{13})$, then
\[
0.7446     \left(\frac{2q}{1-q}\right)^{\rho}\leq \nu(q)\leq 0.9526\left(\frac{2q}{1-q}\right)^{\rho}.
\]

 \end{proposition}
\begin{proof}
The result follows from~\cref{lemma:duplicate}\ref{item:duplicate-rate}
with $a_0=v_0=(1+q)/(1-q)$ and $C=\nu(q)$.
\end{proof}

The following comparison extends the growth estimate in~\Cref{prop:prim-asymp} to the \OConPF{} and \OConGP{} schedules up to a constant factor.

\begin{proposition}\label{prop:prim-asymp-2}
  Suppose that $0<q<1$. Then, for every $n\in\N$,
\begin{equation}\label{eq:sigma-comparison}
\sigma (h_\algA^{(n)})\leq \sigma(h_\algB^{(n)})=\sigma(h_\algC^{(n)})\leq \frac{2}{1+q}\sigma(h_\circ^{(n)}).
\end{equation}
\end{proposition}
\begin{proof}
    The equality follows from~\cref{prop:dp-are-primitive-dom}. We next prove $\sigma (h_\algA^{(n)})\leq \sigma(h_\algB^{(n)})\leq \frac{2}{1+q}\sigma(h_\circ^{(n)})$  by induction. When $n=0$, $\sigma (h_\algA^{(n)})= \sigma(h_\algB^{(n)})=\sigma([~])=0$. Fix $N\geq1$ and suppose that both inequalities hold for all $0\leq n<N$.  
    
    Set $\tau=2/(1+q)$. By~\cref{lemma:monotonicity-of-phi-psi}, $\varphi$ and $\psi$ are increasing in each argument.

By the definition of $h_\algA^{(N)}$, choose an integer $0\leq k<N$ such that $h_\algA^{(N)}=\ConPP(h_\algA^{(k)},h_\algA^{(N-k-1)})$.
    Let $s_a=\sigma(h_\algA^{(k)})$, $s_b=\sigma(h_\algA^{(N-k-1)})$ and $s_c=\sigma(h_\algB^{(N-k-1)})$.
The induction hypothesis gives $s_b\leq s_c$. Therefore,
    \[ 
    \begin{aligned}
\sigma(h_\algA^{(N)}) &=\sigma(\ConPP(h_\algA^{(k)},h_\algA^{(N-k-1)}))\overset{\mathclap{\eqref{eq:recursion-of-sigma2}}}{=}
    \frac{(q s_a+1)(q s_b+1)(1-q \varphi(s_a,s_b))^{-1} - 1}{q}
    \\ &\leq     \frac{(q s_a+1)(q s_c+1)(1-q \psi(s_a,s_c))^{-1} - 1}{q}\\
    &\overset{\mathclap{\eqref{eq:recursion-of-sigma2}}}{=}\sigma(\ConPF(h_\algA^{(k)},h_\algB^{(N-k-1)}))\leq \sigma(h_\algB^{(N)}),
    \end{aligned}
    \] 
where the first inequality uses $s_b\leq s_c$, the monotonicity of $\varphi$, and $\varphi(s_a,s_c)<\psi(s_a,s_c)$ from~\cref{lemma:relation-phi-psi}.

    For the other inequality, choose an integer $0\leq l<N$ such that $h_\algB^{(N)}=\ConPF(h_\algA^{(l)},h_\algB^{(N-l-1)})$.
    Let $s_a=\sigma(h_\algA^{(l)})$, $s_b=\sigma(h_\algB^{(N-l-1)})$ and $s_c=\sigma(h_\algA^{(N-l-1)})$.  By the induction hypothesis, $s_b\leq \tau s_c$.
    By~\cref{lemma:condition-of-sd}, both $1-q\psi(s_a,s_b)$ and $1-q\psi(s_a,\tau s_c)$ are positive. Thus
    \[    \begin{aligned}
\sigma(h_\algB^{(N)}) &=\sigma(\ConPF(h_\algA^{(l)},h_\algB^{(N-l-1)}))\overset{\mathclap{\eqref{eq:recursion-of-sigma2}}}{=}
    \frac{(q s_a+1)(q s_b+1)(1-q \psi(s_a,s_b))^{-1} - 1}{q}
    \\ &\leq     \frac{(q s_a+1)(q (\tau s_c)+1)(1-q \psi(s_a,\tau s_c))^{-1} - 1}{q} \\ &\leq\frac{\tau\big((q s_a+1)(q s_c+1)(1-q \varphi(s_a,s_c))^{-1} - 1\big)}{q}\\
    &\overset{\mathclap{\eqref{eq:recursion-of-sigma2}}}{=}\tau\sigma(\ConPP(h_\algA^{(l)},h_\algA^{(N-l-1)}))\leq \tau\sigma(h_\algA^{(N)}),
    \end{aligned}  
    \]
where the first inequality uses $s_b\leq\tau s_c$ and the monotonicity of $\psi$ from~\cref{lemma:monotonicity-of-phi-psi}. The second inequality follows from~\cref{lemma:relation-phi-psi} with $x=s_a$ and $y=s_c$.\end{proof}

Combining these growth estimates with the sharp bounds in~\cref{theorem:general-bound} gives the following iteration complexity bounds for the \OConPF{} and \OConGP{} schedules.

\begin{theorem}[Iteration complexity bounds]
\label{theorem:iter-number-lower-bound}
Let $\kappa>1$ and $\epsilon\in(0,1)$.
Then gradient descent with $h_\algB^{(n)}$ and $h_\algC^{(n)}$ requires the same minimum number of iterations, denoted by $N(\kappa,\epsilon)$, to guarantee, respectively,
\[
f_n-f_\star\leq\epsilon\cdot\frac{L}{2}\|\x_0-\x_\star\|^2,
\qquad
\frac{1}{2L}\|\g_n\|^2\leq\epsilon\cdot(f_0-f_\star)
\]
for every $L$-smooth-$\mu$-strongly convex function and every initial point. Moreover, if $0<\epsilon\le1/(2\kappa)$, then

\begin{equation}\label{eq:simple-iteration-complexity}
\frac{1-\kappa^{-2}}{16}\,\kappa^{\rho}\log\frac{1}{\kappa\epsilon}
\le N(\kappa,\epsilon)
\le 1+3\kappa^{\rho}\log\frac{1}{\kappa\epsilon}.
\end{equation}
In particular, for each fixed $\kappa>1$ and $0<\epsilon\le1/(2\kappa^2)$,
\[
N(\kappa,\epsilon)=\Theta\!\left(\kappa^\rho\log\frac1\epsilon\right).
\]

\end{theorem}

\begin{proof}
By~\cref{prop:dp-are-primitive-dom,theorem:general-bound}, the two schedules have the same sharp normalized error bound at every iteration budget, so their minimum iteration counts coincide. It therefore suffices to consider $h_\algB^{(n)}$.

For every $\epsilon\in(0,1)$, we first establish the finer bounds
\begin{equation}\label{eq:iteration-budget-bounds}
\begin{aligned}
\left\lceil\underbrace{\frac{\frac12\log(q/\epsilon+1-q)+\log\frac{1+q}{1-q}}{\nu(q)}-1}_{=:A_-}\right\rceil
&\leq N(\kappa,\epsilon)\leq\left\lceil\underbrace{\frac{\frac12\log(q/\epsilon+1-q)+\log\frac{2}{1-q}+\frac{(1-q)^2}{2(1+q)}}{\nu(q)}-1}_{=:A_+}\right\rceil.
\end{aligned}
\end{equation}

Solving $1+2s+q s^2=1/\epsilon$ for $s>0$ gives
\begin{equation}\label{eq:accuracy-threshold}
s=\frac{\sqrt{1-q+q/\epsilon}-1}{q}.
\end{equation}
By sharpness, $h_\algB^{(n)}$ guarantees $\epsilon$-accuracy over the function class if and only if $\sigma(h_\algB^{(n)})\geq s$.
Set $c=(1-q)^2/(2(1+q))$. Combining~\eqref{eq:prim-asymp} and the bounds $-c\leq R_n^\circ\leq0$ in~\eqref{eq:prim-remainder-bound} with~\cref{prop:prim-asymp-2} gives
\begin{gather}
\frac{(1-q)e^{(n+1)\nu(q)-c}-(1+q)}{2q}
\leq\sigma(h_\algA^{(n)})\leq\sigma(h_\algB^{(n)}),\label{eq:dp-sigma-lower-bound}\\
\sigma(h_\algB^{(n)})\leq\frac{2}{1+q}\sigma(h_\algA^{(n)})
\leq\frac{(1-q)e^{(n+1)\nu(q)}-(1+q)}{q(1+q)}.\label{eq:dp-sigma-upper-bound}
\end{gather}
If $n\geq A_+$, then
\[
\begin{aligned}
(1-q)e^{(n+1)\nu(q)-c}
&\geq(1-q)e^{(A_++1)\nu(q)-c}\overset{\eqref{eq:iteration-budget-bounds}}{=}2\sqrt{q/\epsilon+1-q}\overset{\eqref{eq:accuracy-threshold}}{=}2(q s+1)\overset{q<1}{\geq}2q s+1+q.
\end{aligned}
\]
The lower bound~\eqref{eq:dp-sigma-lower-bound} therefore gives $\sigma(h_\algB^{(n)})\geq s$. Hence $N(\kappa,\epsilon)\leq\lceil A_+\rceil$.
If $n\leq A_-$, then
\[
(1-q)e^{(n+1)\nu(q)}\leq(1-q)e^{(A_-+1)\nu(q)}\overset{\eqref{eq:iteration-budget-bounds}}{=}(1+q)\sqrt{q/\epsilon+1-q}\overset{\eqref{eq:accuracy-threshold}}{=}(1+q)(q s+1).
\]
The upper bound~\eqref{eq:dp-sigma-upper-bound} therefore gives $\sigma(h_\algB^{(n)})\leq s$. By sharpness, there exists a function whose normalized error is at least $\epsilon$. For $n<A_-$, the inequality is strict, so sharpness yields a function whose normalized error exceeds $\epsilon$. Thus no integer $n<A_-$ guarantees $\epsilon$-accuracy, which gives $N(\kappa,\epsilon)\geq\lceil A_-\rceil$ and proves~\eqref{eq:iteration-budget-bounds}.

To prove~\eqref{eq:simple-iteration-complexity}, set $t=q/\epsilon\ge2$ and $B=\log((1+q)/(1-q))$. When $q<3/13$,~\cref{prop:nu-bound} gives $(5/4)q^{\rho}\le\nu(q)\le(5/2)q^{\rho}$. When $q\ge3/13$, the same proposition gives $B\leq \nu(q)\leq 2B$, while $q^{\rho}\le2q\le B\le2q/(1-q^2)$ and $q\le q^{\rho}$. Thus, for all $0<q<1$,
\begin{equation}\label{eq:nu-q-bound-q-pow-rho}
q^{\rho}\le\nu(q)\le\frac{4q^{\rho}}{1-q^2}.
\end{equation}
Using the upper bound in~\eqref{eq:iteration-budget-bounds}, $\lceil A_+\rceil\le A_+ +1$,
and $\log\frac{2}{1-q}=B+\log\frac{2}{1+q}$,
% and the estimates
% \[
% t+1-q\leq\frac{3t}{2}\qquad
% \log\frac{2}{1-q}\le B+\log2,\qquad
% \frac{(1-q)^2}{2(1+q)}\le\frac12,\qquad
% \frac{B}{\nu(q)}<1,
% \]
we obtain
\[
% N(\kappa,\epsilon)\le 1+\frac{\frac12\log t+\frac12\log\frac32+\log2+\frac12}{\nu(q)}
% \le1+\frac{3\log t}{q^{\rho}}.
N(\kappa,\epsilon)\leq 
 \frac{\frac12\log (t+1-q)+B+\log\frac{2}{1+q}+\frac{(1-q)^2}{2(1+q)}}{\nu(q)}
\leq \frac{\frac12\log(t+1)+\nu(q)+\log 2+\frac{1}{2}}{\nu(q)}
\leq 1+\frac{3\log t}{q^{\rho}},
\]
where the second inequality uses $B\leq \nu(q)$ and $q< 1$, while
the last inequality uses  $\nu(q)\geq q^\rho$ and $t\geq 2$.
This proves the upper bound in~\eqref{eq:simple-iteration-complexity}.

For the lower bound,~\eqref{eq:iteration-budget-bounds} gives $N(\kappa,\epsilon)\ge(\log t)/(2\nu(q))-1$.
Also, $N(\kappa,\epsilon)\ge1$ because the sharp bound for the empty schedule is $1>\epsilon$.
Adding up the two inequalities yields
\[
N(\kappa,\epsilon)\ge\frac{\log t}{4\nu(q)}\ge\frac{1-q^2}{16q^{\rho}}\log t,
\]
which proves the lower bound in~\eqref{eq:simple-iteration-complexity}.

For each fixed $\kappa>1$ and $0<\epsilon\le1/(2\kappa^2)$, the final assertion follows from~\eqref{eq:simple-iteration-complexity} and
$
\frac12\log\frac1\epsilon\le\log\frac1{\kappa\epsilon}\le\log\frac1\epsilon$.\end{proof}

\begin{figure}[H]
    \centering
\includegraphics[width=0.95\linewidth]{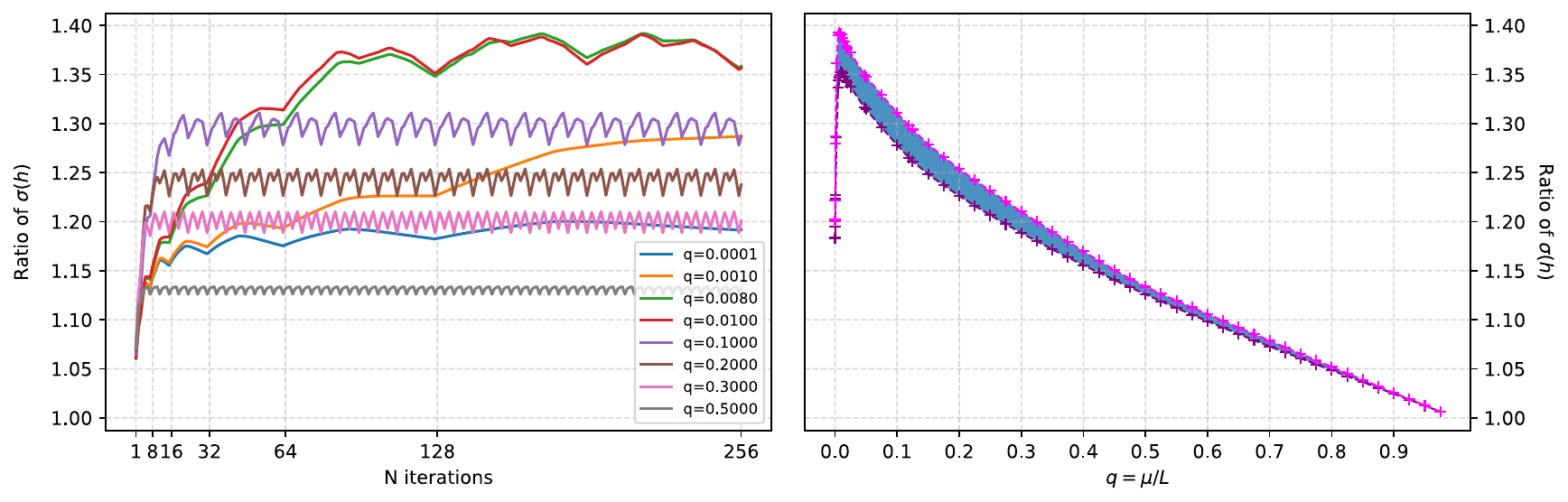} 
    \caption{The ratio $\sigma(h_\bullet)/\sigma(h_\circ)$ of the cumulative quantities of the \OConPF{} and \OConPP{} schedules, computed with high-precision arithmetic.
    Left: the ratio as a function of $n$ for several values of $q=\mu/L$.
    Right: its minimum and maximum over $n\in[128,256]$ as functions of $q$.}
    \label{fig:ratio-wrt-n-q}
\end{figure}

\begin{figure}[H]
    \centering
\includegraphics[width=0.95\linewidth]{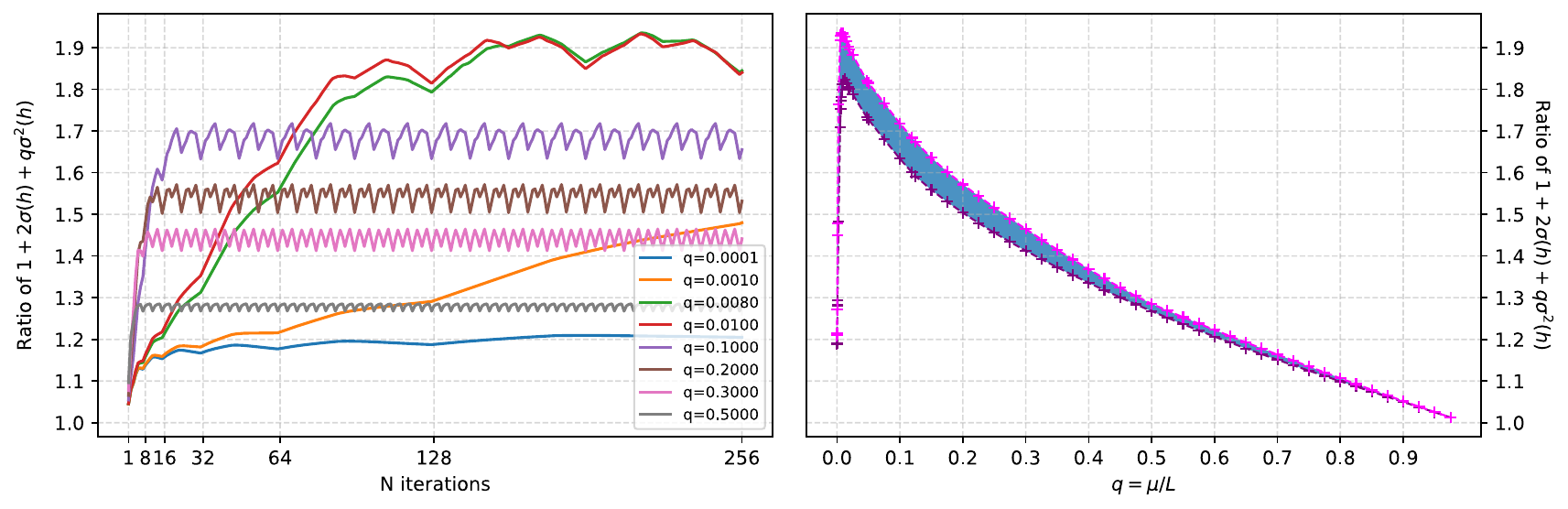} 
    \caption{Improvement of \OConPF{} over \OConPP{} in the function-value bound, measured by the ratio $(1+2\sigma(h_\bullet)+q\sigma^2(h_\bullet))/(1+2\sigma(h_\circ)+q\sigma^2(h_\circ))$ and computed with high-precision arithmetic. Left: the ratio as a function of $n$ for several values of $q=\mu/L$. Right: its minimum and maximum over $n\in[128,256]$ as functions of $q$.}
    \label{fig:bound-ratio-wrt-n-q}
\end{figure}

\begin{remark}[Improving the constants]\label{rem:improving-constants}
The \OConPP{} and \OConPF{} schedules have the same exponential growth rate, but the larger cumulative stepsizes of $h_\algB^{(n)}$ give smaller finite-step function-value bounds in~\cref{theorem:general-bound}. The ratio $\sigma(h_\algB^{(n)})/\sigma(h_\algA^{(n)})$ in Figure~\ref{fig:ratio-wrt-n-q} illustrates the constant-factor advantage of the \OConPF{} schedules for several values of $q$. This advantage translates into tighter function-value bounds, as shown in Figure~\ref{fig:bound-ratio-wrt-n-q}.

There may also be room to improve the constants in the iteration bounds of~\cref{theorem:iter-number-lower-bound}. Its proof estimates $\sigma(h_\algB^{(n)})$ using the bounds for $\sigma(h_\algA^{(n)})$ and the comparison~\eqref{eq:sigma-comparison}, which may lose some precision. A more direct approach is to estimate the remainder $R_n^\bullet$ in~\eqref{eq:fdom-remainder}, plotted in Figure~\ref{fig:remainder-term} (right). Indeed,
\[
\sigma(h_\bullet^{(n)})=\frac{(1-q)e^{(n+1)\nu(q)+R_n^\bullet}-(1+q)}{2q}.
\]
Thus, a sharper bound on $R_n^\bullet$ could improve the iteration bounds in~\eqref{eq:iteration-budget-bounds}.
\end{remark}

%\Needspace{0.38\textheight}

\begin{remark}[Near-optimal iteration complexity]\label{coro:near-optimal-complexity}
The complexity order is already close to optimal among predetermined nonnegative stepsize schedules. Liu and Ye~\cite[Corollary~1.2]{liu2026silverstronglyconvex} prove lower bounds for relative squared distance and relative function error. By strong convexity, these bounds transfer to the function-value gap and gradient-norm criteria in~\cref{theorem:iter-number-lower-bound} with tolerance $\delta=\kappa\epsilon$, giving a lower bound of $\Omega(\kappa^{\rho-o(1)}\log(1/(\kappa\epsilon)))$.
For sufficiently large $\kappa$ and $0<\epsilon\leq1/(2\kappa)$, this matches~\eqref{eq:simple-iteration-complexity} up to a subpolynomial factor in $\kappa$. In particular, for sufficiently large $\kappa$ and $0<\epsilon\leq1/(2\kappa^2)$, \eqref{eq:simple-iteration-complexity} gives $N(\kappa,\epsilon)=\Theta(\kappa^\rho\log(1/\epsilon))$, while the lower bound in~\cite{liu2026silverstronglyconvex} becomes $\Omega(\kappa^{\rho-o(1)}\log(1/\epsilon))$.
\end{remark}

\paragraph{Smooth convex functions.}
In the smooth convex case $q=0$, we have $\sigma_0(h)=\onet h$ and recover the three schedule classes in~\cite{zhang2024accelerated}, where \dom schedules are called dominant schedules. The bounds in~\cref{theorem:bound-of-fdom,theorem:bound-of-dual} recover~\cite[Theorems~2.6 and~5.2]{zhang2024accelerated}, and~\cref{prop:primitive-is-fdom} specializes to~\cite[Proposition~2.8]{zhang2024accelerated}. The connecting steps become
\begin{equation*}
\begin{aligned}
\varphi_0(x,y)&=\frac{-x-y+\sqrt{(x+y+2)^2+4(x+1)(y+1)}}{2},\\
\psi_0(x,y)&=\frac{3-2y+\sqrt{(2y+1)(2y+8x+9)}}{4}.\\
\end{aligned}
\end{equation*}
Thus,~\cref{theorem:conpp,theorem:conpf,theorem:congp} recover the concatenation rules in~\cite[Theorems~3.1, 3.2, and~5.3]{zhang2024accelerated}, with corresponding composition rules in~\cite{grimmer2024composing}.
Both $\sigma_0(h_\algA^{(n)})$ and $\sigma_0(h_\algB^{(n)})$ grow as $\Theta(n^{1/\rho})$, as shown in~\cite[Propositions~4.8, 4.9, 4.11, and~4.12]{zhang2024accelerated} and~\cite[Theorems~4 and~5]{grimmer2024composing}.
Liu et al.~\cite{liu2026optimalrecursivecompositiondyadic} further analyzed the asymptotic behavior of these recursions.

Since $\sigma_0(h)=\Theta(n^{1/\rho})$ for either schedule,~\eqref{eq:bound-fdom} with $q=0$ yields
\begin{equation}\label{eq:convex-function-value-rate}
f_n-f_\star\leq\frac{1}{1+2\sigma_0(h)}\cdot\frac{L}{2}\|\x_0-\x_\star\|^2
=\Theta(n^{-1/\rho})\cdot\frac{L}{2}\|\x_0-\x_\star\|^2.
\end{equation}
Ye and Liu~\cite[Theorem~1.1]{ye2026silver} proved a worst-case lower bound of $\Omega(n^{-1/\rho-o(1)})$ on the normalized function-value error for gradient descent with predetermined nonnegative stepsizes. Their lower bound matches the function-value convergence rate in~\eqref{eq:convex-function-value-rate} up to a subpolynomial factor in $n$.

\paragraph{Weakly convex functions.}
In the weakly convex case $q<0$, the bounds for \dom and \dual schedules alone do not guarantee arbitrarily small errors, because the factor
$\frac{1}{1+2 \sigma(h)+q\sigma^2(h)}$ in~\eqref{eq:bound-fdom} and~\eqref{eq:bound-gdom} is bounded away from zero. Indeed,
\begin{equation}
    \label{eq:limited-performance}
 {1+2 \sigma(h)+q\sigma^2(h)}
 =q(\sigma(h)+\frac{1}{q})^2+1-\frac{1}{q}\leq 1-\frac{1}{q}<+\infty.
\end{equation}
This is consistent with the possibility that gradient descent converges to a local minimizer in a nonconvex problem.

Restarting gives a stationarity guarantee when $f$ is bounded below. Suppose $f(\x)\ge f_\star$ for all $\x$, and fix a \dual schedule $h\in\R^n$ with $2\sigma(h)+q\sigma^2(h)>0$. Apply $h$ repeatedly, starting each run from the endpoint of the previous one. Denote the endpoint of the $k$-th run by $\x_{kn}$ and write $\g_{kn}=\nabla f(\x_{kn})$. Summing the decrease inequality in~\cref{def:g-dominant} over $N$ runs gives
\[
\frac{1}{2L}\sum_{k=1}^N\|\g_{kn}\|^2\leq \frac{f_0-f_{Nn}}{2\sigma(h)+q\sigma^2(h)}\leq \frac{f_0-f_\star}{2\sigma(h)+q\sigma^2(h)}.
\]
Hence the smallest squared gradient norm among the restart points is $O(1/N)$.

\section{Silver Schedules and Dynamic Concatenation}
\label{sect:connections-extensions}
In this section, we first extend the optimized primitive schedules to obtain distance and gradient-norm guarantees for arbitrary iteration budgets. We show that these schedules recover the silver schedule and its reversal, respectively, at power-of-two budgets.
%, with uniquely optimal symmetric splits in the primitive recursion.
We then develop dynamic schedules through repeated concatenation with a fixed primitive block, without prescribing the iteration budget in advance.

\subsection{The Silver Schedule}
\label{sect:silver-stepsizes}

We begin by extending primitive schedules to obtain distance and gradient-norm guarantees. The following proposition gives the two constructions and their corresponding bounds.

% \begin{proposition}\label{theorem:bound-x}
% \leavevmode
%     Let $q<1$, and let $h_a\in\rk{n}$ be a primitive schedule. Set $s=\sigma(h_a)$ and assume that $2qs+q+1>0$. Let $\x_n$ be the final iterate generated by gradient descent with $h_a$, and take one additional step
%     \[
%         \x_{n+1}=\x_n-\frac{q s+s+2}{2q s+q+1}\frac{\g_n}{L}.
%     \]
%     Then $\x_{n+1}$ satisfies
%     \begin{equation}
%         \label{eq:primitive-x-coro}
%     \|\x_{n+1}-\x_\star\|\leq  \frac{ 1-q }{2q s+q +1}\|\x_0-\x_\star\|.\end{equation}
%     For $0<q<1$, the bound~\eqref{eq:primitive-x-coro} is sharp, with equality attained by the one-dimensional quadratic $f(x)=\mu x^2/2$ for every nonzero initial point.
% \end{proposition}

% \begin{proposition}\label{theorem:bound-g}
% \leavevmode
%     Let $q<1$, and let $h _a\in\rk{n}$ be a primitive schedule. Set $s=\sigma(h_a)$ and assume that $2qs+q+1>0$. 
%    Gradient descent with $h=[\frac{q s+s+2}{2q s+q+1},h_a^\top]\in\rk{n+1}$ has
%     \begin{equation}
%         \label{eq:primitive-g-coro}
%     \|\g_{n+1} \|\leq  \frac{ 1-q }{2q s+q +1}\|\g_0 \|.\end{equation}
%     For $0<q<1$, the bound~\eqref{eq:primitive-g-coro} is sharp, with equality attained by the one-dimensional quadratic $f(x)=\mu x^2/2$ for every nonzero initial point.
% \end{proposition} 
\begin{proposition}\label{theorem:bound-x}
Let $q<1$, and let $h \in\rk{n}$ be a primitive schedule. Set $s=\sigma(h)$ and assume that $2qs+q+1>0$. Then the following assertions hold.
% Let $\x_n$ be the final iterate generated by gradient descent with $h $, and take one additional step
\begin{enumerate}[label=(\roman*),ref=(\roman*),leftmargin=*]
\item\label{item:bound-x}
Gradient descent with $h_a=[h^\top,\frac{qs+s+2}{2qs+q+1}]^\top\in\rk{n+1}$ satisfies
  \begin{equation}
        \label{eq:primitive-x-coro}
    \|\x_{n+1}-\x_\star\|\leq  \frac{ 1-q }{2q s+q +1}\|\x_0-\x_\star\|.
\end{equation}  
\item \label{item:bound-g}
Gradient descent with $h_b=[\frac{qs+s+2}{2qs+q+1},h^\top]^\top\in\rk{n+1}$ satisfies
\begin{equation}
        \label{eq:primitive-g-coro}
    \|\g_{n+1} \|\leq  \frac{ 1-q }{2q s+q +1}\|\g_0 \|.
\end{equation}
\end{enumerate} 
For $0<q<1$, both bounds~\eqref{eq:primitive-x-coro} and~\eqref{eq:primitive-g-coro} are sharp, with equality attained by the one-dimensional quadratic $f(x)=\mu x^2/2$ for every nonzero initial point.
\end{proposition}

\begin{proof}
\textit{\ref{item:bound-x}}
 Using the primitive certificate difference~\eqref{eq:certificate-primitive} with $(\ell,r)=(0,n)$, set
\begin{equation*}
   A_1\triangleq \frac L2\|\x_0-\x_\star\|^2
   +(1-q)\sum_{i=0}^{n-1}\omega_i^{0,n}Q_{i,\star}-\mathcal C^{\mathrm p}_{0,n}.
\end{equation*}
Since $\mathcal C^{\mathrm p}_{0,n}\geq0$, $\omega_i^{0,n}\geq0$, and $Q_{i,\star}\leq0$, we have $A_1\leq\frac L2\|\x_0-\x_\star\|^2$.

Direct expansion using~\eqref{eq:certificate-primitive} and \eqref{eq:Qij} gives
\begin{equation}
\label{eq:cert-x}
\begin{aligned}
   &A_1+\frac{(s+1-q )(2q s+q +1)+q^2(s+1)+(q s+1)}{1-q}Q_{n,\star}+\frac{2(s+1)(q s+1)}{1-q }Q_{\star,n}\\ &=\left(\frac{2q s+q +1}{1-q}\right)^2\cdot \frac{L}{2}\|\x_n - \frac{q s+s+2}{2q s+q +1}\g_n/L-\x_\star\|^2.
\end{aligned}%\tag{Certificate-X}
\end{equation}
By~\cref{lemma:s-bound}, $s\geq0$ and $q s+1>0$. Together with $q<1$ and the assumption $2qs+q+1>0$, these inequalities ensure that the coefficients of $Q_{n,\star}$ and $Q_{\star,n}$ in~\eqref{eq:cert-x} are nonnegative. Since $Q_{n,\star},Q_{\star,n}\leq0$, the right-hand side of~\eqref{eq:cert-x} is at most $A_1\leq\frac L2\|\x_0-\x_\star\|^2$. Taking square roots and using the definition of $\x_{n+1}$ gives~\eqref{eq:primitive-x-coro}.

\textit{\ref{item:bound-g}}
    Set $\a=\frac{qs+s+2}{2qs+q+1}$, so that $\a>1$ and $1-q\a>0$.
    Applying~\eqref{eq:cert-sd} to the schedule $h_b=[\a,h^\top]^\top$ and using the nonnegativity of its left-hand side gives
    \[
0+\frac{2(1- q)^2(s+1)^2( q s+1)^2}{(2 q s+q +1)^3}\|\g_0\|^2
-\frac{2(s+1)^2( q s+1)^2}{2 q s+ q +1}\|\g_{n+1}\|^2\geq 0,
    \]
    which therefore implies~\eqref{eq:primitive-g-coro}.

To prove sharpness for $0<q<1$, consider $f(x)=\mu x^2/2$ with $x_0\neq0$, as in~\cite[Remark~1.3 and Appendix~C]{silver1}.
By~\eqref{eq:sigma},
\[
x_{n+1}=x_0\left(1-q\cdot \frac{q s+s+2}{2qs +q +1}\right)\prod_{i=0}^{n-1}(1-q h_i)
=x_0\cdot \frac{(1-q)(qs+1)}{2qs+q+1}\cdot \frac{1}{1+qs}=\frac{1-q}{2qs+q+1}x_0
\]
holds for both cases, and the desired equalities follow.\end{proof}

In the above proposition,
the condition $2qs+q+1>0$ is automatic for $0\leq q<1$. For $q<0$, it is equivalent to $-1<q<0$ and $s<(1+q)/(-2q)$, so the bound requires $\mu>-L$ and a sufficiently small cumulative stepsize; the common factor in~\eqref{eq:primitive-x-coro} and~\eqref{eq:primitive-g-coro} then exceeds one, and these bounds control the possible growth of the distance to a stationary point and of the gradient norm.

For $0< q<1$ and each $N\in\N_+$, let $h_\circ^{(N-1)}$ be the \OConPP{} schedule and set $s_{N-1}=\sigma(h_\circ^{(N-1)})$. Set
\begin{equation}\label{eq:distance-beta}
\beta(s)=\frac{q s+s+2}{2q s+q+1},\qquad  \tau_N=\left(\frac{1-q}{2q s_{N-1}+q+1}\right)^2
\end{equation}
and define
\begin{equation}\label{eq:optimized-distance-schedule}
\hat h^{(N)}=[(h_\circ^{(N-1)})^\top,\beta(s_{N-1})]^\top,\qquad 
\check h^{(N)}=[\beta(s_{N-1}),(h_\circ^{(N-1)})^\top]^\top.
% \qquad r_N=\frac{1-q}{2q s_{N-1}+q+1}.
\end{equation}
 We call $\hat h^{(N)}$ and $\check h^{(N)}$ the \OConPPr{} schedule and \OConPPl{} schedule, respectively.
Here, \texttt{r} (right) and \texttt{l} (left) indicate adding a step at the end and at the beginning of the primitive schedule, respectively.

By~\cref{theorem:bound-x}, for any $d\geq1$, $f\in\mathcal F_{\mu,L}(\R^d)$, and initial point $\x_0\in\R^d$, gradient descent with $\hat h^{(N)}$
and $\check h^{(N)}$ satisfies, respectively,
\begin{equation}\label{eq:optimized-distance-bound}
\|\x_N-\x_\star\|^2\leq\tau_N\|\x_0-\x_\star\|^2,\qquad
\|\g_N\|^2\leq\tau_N\|\g_0\|^2.
\end{equation}

Previous work recovered the smooth convex silver schedule~\cite{silver2}
through concatenation~\cite{zhang2024accelerated} and composition~\cite{grimmer2024composing}.
We now relate our constructions to the silver schedule for smooth strongly convex optimization~\cite{silver1},
which has length $2^k$ and provides bounds on the relative squared distance
$\frac{\|\x_n-\x_\star\|^2}{\|\x_0-\x_\star\|^2}$.

For the strongly convex silver schedule, let $0<q<1$ and $y_1=z_1=q$. Following~\cite[Equations~(3.1)--(3.2)]{silver1} with $\kappa=1/q$, define
\begin{equation}\label{eq:silver-parameters}
z_{2^{k+1}}=z_{2^k}\left(1-z_{2^k}+\sqrt{1+(1-z_{2^k})^2}\right),\qquad
y_{2^{k+1}}=\frac{z_{2^k}^2}{z_{2^{k+1}}},\qquad k\geq0.
\end{equation}
Starting from $h_{\mathrm{sil}}^{(1)}=[2/(1+q)]^\top$, for each $k\geq0$ let $\widetilde h_{\mathrm{sil}}^{(2^k-1)}$ be the vector obtained by deleting the last entry of $h_{\mathrm{sil}}^{(2^k)}$, with $\widetilde h_{\mathrm{sil}}^{(0)}=[~]$, and define
\begin{equation}\label{eq:silver-schedule}
h_{\mathrm{sil}}^{(2^{k+1})}
=\left[(\widetilde h_{\mathrm{sil}}^{(2^k-1)})^\top,\frac{1+y_{2^{k+1}}/q}{1+y_{2^{k+1}}},(\widetilde h_{\mathrm{sil}}^{(2^k-1)})^\top,\frac{1+z_{2^{k+1}}/q}{1+z_{2^{k+1}}}\right]^\top.
\end{equation}

The following theorem identifies the silver schedule with the \OConPPr{} schedule and its reversal with the \OConPPl{} schedule at every power-of-two iteration budget.
Together with~\eqref{eq:optimized-distance-bound}, the latter identity yields a gradient-norm contraction bound for the reversed silver schedule, consistent with the discussion in~\cite[Section~1.1.4]{silver1}.
%This extends the uniqueness of the optimal symmetric split at primitive lengths $2^k-1$ from the smooth convex case~\cite[Corollary~4.4(i)]{liu2026optimalrecursivecompositiondyadic} to the smooth strongly convex case.

\begin{theorem}
\label{theorem:dyadic-silver-optimality}
Let $0<q<1$. For every $k\in\N$, 
% the \OConPP{} recursion in~\cref{def:dynamic-programming} satisfies
% \[
% \underset{j\in\{0,\dotsc,2^k-2\}}{\argmax}\;
% \sigma_q\!\left(\ConPP_q\!\left(h_\circ^{(j)},h_\circ^{(2^k-2-j)}\right)\right)
% =\{2^{k-1}-1\}.
% \]
% Consequently,
 the \OConPPr{} and \OConPPl{} schedules satisfy
\begin{equation*}
\hat h^{(2^k)}=h_{\mathrm{sil}}^{(2^k)},\qquad
\check h^{(2^k)}=\rev h_{\mathrm{sil}}^{(2^k)},\qquad k\in\N.
\end{equation*}
\end{theorem}

\begin{proof}
\emph{Recovery by symmetric concatenation.}
Starting from $h_\circ^{(0)}=[~]$, recall from~\cref{prop:prim-asymp} we have
\[
h_\circ^{(2^{k+1}-1)}=\ConPP(h_\circ^{(2^k-1)},h_\circ^{(2^k-1)}),\qquad k\in\N.
\]
Set $s_k=\sigma(h_\circ^{(2^k-1)})$ and let $y_{2^k},z_{2^k}$ be defined by~\eqref{eq:silver-parameters}. We first identify the last step of the silver schedule with $\beta(s_k)$. Set $t_k = \frac{q s_k+1}{1-q}>1$. Then $t_0=\frac{1}{1-q}$. Substituting $t_a=t_b=t_k$ and $t_c=t_{k+1}$ into~\eqref{eq:concat-rate-v} and using $t_k>1$ gives
\begin{equation}
    \label{eq:silver-recur}
  t_{k+1}=t_k^2-t_k+1+(t_k-1)\sqrt{t_k^2+1}.
\end{equation}
This rearranges to
\begin{equation}\label{eq:silver-recur-zb}
1-\frac{1}{t_{k+1}}=
\frac{t_k-1+(t_k-1)\sqrt{t_k^2+1}}{t_k^2}=
\left(1 - \frac{1}{t_k}\right)
\left( \frac{1}{t_k}
+\sqrt{1+\frac{1}{t_k^2}} 
\right).
\end{equation}
Thus the sequence $\{1-1/t_k\}$ satisfies the same recurrence as $\{z_{2^k}\}$. Since $1-1/t_0=q=z_1$, induction gives $z_{2^k}=1-1/t_k=q(s_k+1)/(q s_k+1)$ for all $k\in\N$, where the last equality follows from the definition of $t_k$ preceding~\eqref{eq:silver-recur}. Hence the last step of the silver schedule is
\begin{equation*}
\frac{1+z_{2^k}/q}{1+z_{2^k}}
=\frac{qs_k+s_k+2}{2q s_k+q+1}\overset{\eqref{eq:distance-beta}}{=}\beta(s_k).
\end{equation*}

It remains to show that $h_\circ^{(2^k-1)}=\widetilde h_{\mathrm{sil}}^{(2^k-1)}$. Using $z_{2^k}=1-1/t_k$, the connecting step $\a_{k+1}=\varphi(s_k,s_k)$ satisfies
\begin{equation*}
\frac{1-q}{1-q\a_{k+1}}
\overset{\eqref{eq:concat-rate-v}}{=}\frac{t_{k+1}}{t_k^2}
\overset{\eqref{eq:silver-recur}}{=}1+\frac{t_k-1}{1+\sqrt{t_k^2+1}}
\overset{\eqref{eq:silver-recur-zb}}{=}1+\frac{z_{2^k}^2}{z_{2^{k+1}}}
\overset{\eqref{eq:silver-parameters}}{=}1+y_{2^{k+1}}.
\end{equation*}
Hence $\a_{k+1}=(1+y_{2^{k+1}}/q)/(1+y_{2^{k+1}})$, as in~\eqref{eq:silver-schedule}.

By~\eqref{eq:silver-schedule}, the induction hypothesis $h_\circ^{(2^k-1)}=\widetilde h_{\mathrm{sil}}^{(2^k-1)}$ therefore implies $h_\circ^{(2^{k+1}-1)}=\widetilde h_{\mathrm{sil}}^{(2^{k+1}-1)}$.
Since $h_\circ^{(0)}=\widetilde h_{\mathrm{sil}}^{(0)}=[~]$, the identity holds for all $k\in\N$ by induction.
Appending $\beta(s_k)$ therefore yields
$\hat{h}^{(2^k)}=h_{\mathrm{sil}}^{(2^k)}$ for every $k\in\N$.
The symmetric concatenation above shows by induction that $\rev h_\circ^{(2^k-1)}=h_\circ^{(2^k-1)}$.
By~\eqref{eq:optimized-distance-schedule},
we immediately have
$\check h^{(2^k)}=\rev\hat h^{(2^k)}=\rev h_{\mathrm{sil}}^{(2^k)}$.
% \medskip
% \emph{Uniqueness of the optimal split.}
% Set $a_n=w_{n-1}$, with $w_n$ defined in~\eqref{eq:wn-definition}. By~\eqref{eq:concat-rate-eqw} and~\cref{def:dynamic-programming}, $\{a_n\}$ satisfies~\eqref{eq:max-split-recurrence} with $a_1=(1+q)/(1-q)>1$ and $F$ as in~\cref{lemma:F-lemma-four-points}. By~\cref{lemma:upper-bound-Rn}\ref{item:dyadic-optimal-split}, the unique maximizer of $F(a_m,a_{n-m})$ over $1\leq m<n$ at $n=2^k$ is $m=2^{k-1}$ for every $k\geq1$.
% For $n=2^k$, $k\geq1$, the unique split $m=n/2$ corresponds to two primitive blocks of length $m-1=2^{k-1}-1$. Applying this recursively gives $h_\circ^{(2^k-1)}=h^{(2^k-1)}$. By~\eqref{eq:silver-terminal-step}, appending $\beta(s_k)$ as in~\eqref{eq:optimized-distance-schedule} therefore yields $\hat h^{(2^k)}=h_{\mathrm{sil}}^{(2^k)}$. The case $k=0$ follows from $\beta(0)=2/(1+q)$.
\end{proof}

\begin{remark}
The following cases illustrate the potential optimality and nonuniqueness of our schedules. %\todo{purples of the remark?}
For $N=2$, the \OConPPr{} schedule is the unique optimal two-step schedule in the distance metric by~\cite[Theorem~2.2]{silver1}
as it matches the silver schedule.
For $N=3$ and $q=1/4$, the schedules are not unique. Different maximizers in~\cref{def:dynamic-programming} give
\begin{equation}\label{eq:optimized-silver-schedule-n-3}
\begin{aligned}
\hat h^{(3)}&=\left[\frac{4}{3},2\sqrt 3-2,\frac{14+6\sqrt 3}{11}\right]^\top\qquad\text{or}\qquad \left[2\sqrt 3-2,\frac{4}{3},\frac{14+6\sqrt 3}{11}\right]^\top,\\
\check h^{(3)}&=\left[\frac{14+6\sqrt 3}{11},\frac{4}{3},2\sqrt 3-2\right]^\top\qquad\text{or}\qquad \left[\frac{14+6\sqrt 3}{11},2\sqrt 3-2,\frac{4}{3}\right]^\top.
\end{aligned}
\end{equation}
The schedules in~\eqref{eq:optimized-silver-schedule-n-3} coincide with the optimal
schedules in the distance and  gradient-norm metrics derived in~\cite[Section 8.2.3]{altschuler2018greed}
for $N=3$, $q=1/4$, after converting the latter to our normalization.
\end{remark}

The following corollary gives the iteration complexity of the \OConPPr{} and \OConPPl{} schedules, recovering the distance complexity of the silver schedule at power-of-two budgets as established in~\cite[Theorem~1.2]{silver1}. Its proof is given in Appendix~\ref{sect:proofs-section5}.
\begin{corollary}[Iteration complexity]\label{coro:optimized-distance}
For the \OConPPr{} and \OConPPl{} schedules, let $N_2(\kappa,\epsilon)$ denote the minimum iteration budget needed to attain the respective guarantees in~\eqref{eq:optimized-distance-bound} with $\tau_N\leq\epsilon$.
If $\kappa>1$ and $0<\epsilon\leq e^{-1}$, then
\[
\frac{1-\kappa^{-2}}{8}\kappa^\rho\log\frac{1}{\epsilon}\leq N_2(\kappa,\epsilon)\leq  1+\kappa^\rho\log\frac{1}{\epsilon}. 
\]
As a consequence, if $\kappa>1$ is fixed and $0<\epsilon\leq e^{-1}$, then
\[
N_2(\kappa,\epsilon)=\Theta\!\left(\kappa^\rho\log\frac{1}{\epsilon}\right).
\]
\end{corollary}

% {For predetermined nonnegative stepsize schedules, the best-known iteration lower bound for reducing the relative squared distance is $\Omega(\kappa^{\rho-o(1)}\log(1/\epsilon))$~\cite[Corollary~1.2]{liu2026silverstronglyconvex}. For sufficiently large fixed $\kappa$ and $0<\epsilon\leq e^{-1}$, this lower bound matches the complexity in~\cref{coro:optimized-distance} up to a subpolynomial factor in $\kappa$.}

\subsection{Dynamic Schedules}
\label{sect:dynamic-schedules}
Repeated concatenation with a fixed primitive block gives schedules that can be extended without fixing the iteration budget in advance. The following proposition gives finite-step guarantees and convergence rates for repeated \ConGP{} and \ConPP{} concatenation. Its proof is given in Appendix~\ref{sect:proofs-section5}.

\begin{proposition}[Dynamic schedules]\label{prop:dynamic-congp}
Let $q<1$, let $h_b\in\rk{m}$ be primitive, and let $\bar h^{(0)}\in\rk{n_0}$, where $n_0,m\in\N$. For the constructions below, write $s_b=\sigma(h_b)$, $s_k=\sigma(\bar h^{(k)})$, $\theta_k=(1+2s_k+q s_k^2)^{-1}$, and $n_k=n_0+k(m+1)$.
For any $f\in\mathcal F_{\mu,L}$ and $\x_0$, the following bounds hold after $n_k$ steps for every $k\in\N$.
\begin{enumerate}[label=(\roman*),ref=(\roman*),leftmargin=*]
\item If $\bar h^{(0)}$ is \dual, define $\bar h^{(k+1)}=\ConGP(\bar h^{(k)},h_b)$ and $\gamma_k=\psi(s_b,s_k)$. Then each $\bar h^{(k)}$ is \dual, and if $f_\star=\inf f$,
\begin{equation}\label{eq:dynamic-congp-bound}
\frac{1}{2L}\|\g_{n_k}\|^2\leq\theta_k(f_0-f_\star).
\end{equation}
\item If $\bar h^{(0)}$ is primitive, define $\bar h^{(k+1)}=\ConPP(\bar h^{(k)},h_b)$ and $\gamma_k=\varphi(s_k,s_b)$. Then each $\bar h^{(k)}$ is primitive, and for any stationary point $\x_\star$, with $f_\star=f(\x_\star)$,
\begin{equation}\label{eq:dynamic-conpp-bound}
f_{n_k}-f_\star\leq\theta_k\cdot\frac L2\|\x_0-\x_\star\|^2.
\end{equation}
\end{enumerate}
For $0<q<1$, both bounds decay geometrically at $n_k$, with
\begin{equation}\label{eq:dynamic-linear-bound}
\theta_k\leq\theta_0\left(\frac{1-q\gamma_0}{1+q s_b}\right)^{2k},
\qquad
\lim_{k\to\infty}\frac{\theta_{k+1}}{\theta_k}
=\left(\frac{1-q}{2q s_b+q+1}\right)^2.
\end{equation}
When $q=0$, both constructions satisfy
\begin{equation}\label{eq:dynamic-congp-asymptotic}
% \lim_{k\to\infty}\frac{s_k}{n_k}=\frac{2(s_b+1)}{m+1}.
\lim_{k\to\infty}\theta_k  n_k= \frac{m+1}{4(s_b+1)}.
\end{equation}
\end{proposition}
At $q=0$,~\eqref{eq:dynamic-congp-asymptotic} gives
\[
\theta_k=\frac{m+1}{4(s_b+1)n_k}(1+o(1)),\qquad k\to\infty.
\]
Thus both bounds~\eqref{eq:dynamic-congp-bound} and~\eqref{eq:dynamic-conpp-bound} have the same asymptotic coefficient $(m+1)/(4(s_b+1))$.
The empty block gives $1/4$, so a block with $s_b>m$ improves this coefficient, as observed in~\cite[Remarks~4.19 and~5.11]{zhang2024accelerated}.
In particular, any nonempty \OConPP{} block satisfies $s_b>m$, since, by~\cref{lemma:roots-of-phi}(ii), all its stepsizes exceed $1$ when $q=0$.

The following corollary recovers the sequence in~\cite[Definition~2.16]{rotaru2024exact} by taking $\bar h^{(0)}=h_b=[~]$ in case~(i). Its proof is given in Appendix~\ref{sect:proofs-section5}.
\begin{corollary}\label{coro:dynamic-empty-block}
Let $q<1$. Define $h^{(0)}=[~]$ and $h^{(n+1)}=\ConGP(h^{(n)},[~])$ for $n\in\N$. Then $h^{(n)}=[h_0,\dotsc,h_{n-1}]^\top$ and $h_n\in[1,2/(1+\max\{q,0\}))$
for all $n\in\N$, and the steps satisfy
\begin{equation}\label{eq:weakly-convex-recur}
    \frac{h_{n+1}[(2-h_{n+1})(2-q h_{n+1})-1]}{2-h_{n+1}(1+q)}+\frac{h_{n }}{2-h_{n }(1+q)}=0
\end{equation}
for $n=-1,0,1,\dotsc$ if we define $h_{-1}=0$. If $f_\star$ is the global infimum, gradient descent with $h^{(n)}$ satisfies
\begin{equation}
\label{eq:weakly-convex-bound}     \frac{1}{2L}\|\g_n\|^2\leq \frac{1}{1+\frac{h_{n-1}}{2-h_{n-1}(1+q)}}(f_0-f_\star),\qquad \forall n\in\N.
\end{equation}
\end{corollary}

For $q=0$, the empty-block case of~\eqref{eq:dynamic-congp-asymptotic} yields
\[
\frac{1}{2L}\|\g_n\|^2\leq\frac{1+o(1)}{4n}(f_0-f_\star),\qquad n\to\infty.
\]
This recovers~\cite[Remark~5.11]{zhang2024accelerated}, and the finite-step bound~\eqref{eq:weakly-convex-bound} coincides with~\cite[Corollary~2.19]{rotaru2024exact}.

For $0<q<1$,~\eqref{eq:weakly-convex-bound} coincides with~\cite[Corollary~2.20]{rotaru2024exact}.
By~\cite[Proposition~2.17(v)]{rotaru2024exact}, we have $0<2/(1+q)-h_{n-1}\leq((1-q)/(1+q))^{2n}$. Substituting this estimate into~\eqref{eq:weakly-convex-bound} yields the geometric bound
\[
\frac{1}{2L}\|\g_n\|^2\leq\frac{(1+q)^2}{2}\left(\frac{1-q}{1+q}\right)^{2n}(f_0-f_\star),\qquad n\geq1.
\]

For $q<0$, the same telescoping argument as in
\cite[Appendix~A.4, Equation~(A.3)]{rotaru2024exact},
combined with Lemma~3.5 therein, also yields
\eqref{eq:weakly-convex-bound} for the untruncated sequence.
% \st{Since \mbox{$h_n\to2$} by \mbox{\cite[Proposition~2.17]{rotaru2024exact}},
% the coefficient in this bound tends to \mbox{$-q/(1-q)>0$},
% so it}
Consistent with~\eqref{eq:limited-performance}, this bound does not imply convergence to stationarity.
In contrast, \cite[Corollary~2.21]{rotaru2024exact}
establishes an $O(1/n)$ bound on the minimum squared
gradient norm for the capped sequence $\widetilde h_i=\min\{h_i,h_{\mathrm{cap}}(q)\}$, where $h_{\mathrm{cap}}(q)\in[1,2)$ is the normalized constant stepsize specified in~\cite[Proposition~2.15]{rotaru2024exact}.

\section{Numerical Comparisons and Optimality Conjectures}
\label{sect:compare}

In this section, we compare the worst-case performance of our proposed \OConPF{} schedules $h_\algB$ and \OConGP{} schedules $h_\algC$ with existing schedules and numerically optimized schedules. These comparisons illustrate the advantages of our constructions.

For $q>0$, we numerically searched for stepsize schedules that minimize the worst-case function-value error over smooth strongly convex functions. We used a branch-and-bound (BNB) method similar to that in~\cite{gupta2024branch}. Since the computational cost grows rapidly with the dimension $n$ of the stepsize sequence, we restricted the search to $n=1,\ldots,6$ and imposed a time limit of 600 seconds for each run, reporting the best schedule found.
Table~\ref{tab:bnb-vs-fdom-.1} compares the resulting schedules $h_{\text{BNB}}$ with our schedules $h_\algB$ for $\mu=0.1$ and $L=1$. Asterisks indicate runs that reached the time limit, for which global optimality was not certified. The two schedules agree closely in all reported cases, with $\|h_{\text{BNB}}-h_\algB\|_{\infty}<4\times 10^{-4}$.

For $n=1$ in~\cref{tab:bnb-vs-fdom-.1}, our schedule $h_\algB$ is exactly optimal in the function-value metric, and the BNB schedule recovers this optimum up to numerical accuracy.  Indeed, for $0<q<1$, \cref{def:dynamic-programming} gives
\begin{equation*}
h_\algB^{(1)}=\ConPF([~], [~])=[\psi(0,0)]^\top=[\frac{q+1-\sqrt{q^2-q+1}}{q}]^\top.
\end{equation*}
This coincides with the optimal one-step schedule derived in~\cite[Appendix A]{taylor2017thesis} for the worst-case normalized function-value error $\frac{f_1-f_\star}{\frac L2\|\x_0-\x_\star\|^2}$.

\begin{table}[!htb]
\centering
\setlength{\belowcaptionskip}{3pt}
\caption{BNB and \OConPF{} schedules and their function-value bounds for $\mu=0.1$ and $L=1$. Asterisks mark runs that reached the 600-second BNB time limit.}
\label{tab:bnb-vs-fdom-.1}
\small
{
\setlength{\tabcolsep}{4pt}
\renewcommand{\arraystretch}{1.1}
\begin{tabular*}{\textwidth}{@{\extracolsep{\fill}}cllr@{}}
\toprule
$n$ & Schedule & $h^\top$ & Bound \\
\midrule
\multirow{2}{*}{1} & $h_{\text{BNB}}$ & $[1.460615]$ & 0.212166 \\
 & $h_{\algB}$ & $[1.460608]$ & 0.212160 \\
\addlinespace
\multirow{2}{*}{2} & $h_{\text{BNB}}$ & $[1.383757, 1.786867]$ & 0.091183 \\
 & $h_{\algB}$ & $[1.383736, 1.786857]$ & 0.091171 \\
\addlinespace
\multirow{2}{*}{3*} & $h_{\text{BNB}}$ & $[1.383784, 2.219750, 1.460622]$ & 0.046495 \\
 & $h_{\algB}$ & $[1.383736, 2.219725, 1.460608]$ & 0.046478 \\
\addlinespace
\multirow{2}{*}{4*} & $h_{\text{BNB}}$ & $[1.550353, 1.383780, 2.648949, 1.460618]$ & 0.025736 \\
 & $h_{\algB}$ & $[1.550317, 1.383736, 2.648856, 1.460608]$ & 0.025723 \\
\addlinespace
\multirow{2}{*}{5*} & $h_{\text{BNB}}$ & $[1.383871, 1.892035, 1.383751, 3.000768, 1.460615]$ & 0.014666 \\
 & $h_{\algB}$ & $[1.383736, 1.892023, 1.383736, 3.000606, 1.460608]$ & 0.014651 \\
\addlinespace
\multirow{2}{*}{6*} & $h_{\text{BNB}}$ & $[1.383961, 1.892034, 1.383762, 3.346594, 1.383752, 1.786858]$ & 0.008674 \\
 & $h_{\algB}$ & $[1.383736, 1.892023, 1.383736, 3.346241, 1.383736, 1.786857]$ & 0.008659 \\
\bottomrule
\end{tabular*}
}
\end{table}
\FloatBarrier

In~\cref{tab:compare-other-schedules}, we compare the \OConPF{} schedule $h_\algB$ with the  original recursively defined silver schedules~\cite{silver1} using the metric $\frac{f_n-f_\star}{\frac{L}{2} \|\x_0-\x_\star\|^2}$ at iteration budgets $2^l$, where $l\in\N$.
 The result in~\cite{silver1} is stated in the distance metric $\frac{\|\x_n-\x_\star\|^2}{\|\x_0-\x_\star\|^2}$,
so we transfer the bound using $\frac{f_n-f_\star}{\frac{L}{2}\|\x_0-\x_\star\|^2}\leq \frac{\|\x_n-\x_\star\|^2}{\|\x_0-\x_\star\|^2}$. This bound is sharp for the silver stepsizes,
where the equality case is attained by the quadratic function  $f(x)=Lx^2/2$, as discussed in~\cite[Remark~1.3 and Appendix~C]{silver1}. 
% also use PEPit~\cite{pepit2022} to evaluate its worst-case performance in the function-value metric.
For the tested curvature ratios and iteration budgets, $h_\algB$ yields smaller worst-case function-value errors than the silver schedule.
We also compare the \OConGP{} schedule $h_\algC$ with the dynamic schedule~\cite[Theorem 2.18]{rotaru2024exact} using the gradient-norm metric $\frac{\|\g_n\|^2/(2L)}{f_0-f_\star}$. For the tested curvature ratios, $h_\algC$ matches the dynamic schedule at $N=1$ and yields smaller worst-case gradient-norm bounds for $N>1$.

\begin{table}[!htb]
    \centering
    \caption{Worst-case function-value bounds for \OConPF{} and the silver schedule (AP)~\cite{silver1}, and gradient-norm bounds for \OConGP{} and the dynamic schedule (RGP)~\cite{rotaru2024exact}. The identical bounds for \OConPF{} and \OConGP{} are shown in a single row.}
    \label{tab:compare-other-schedules}
    \small
    \begin{tabular*}{\textwidth}{@{\extracolsep{\fill}}cccccccc@{}}
        % \hline
        \toprule
        $q=\mu/L$ & $h$ & $N=1$ & $2$ & $4$ &   $8$ & $16$ & $32$  \\
        \midrule%hline
        \multirow{3}{*}{$10^{-1}$} & \texttt{O-ConPF/GP} & \textbf{0.212160} & \textbf{0.091171} & \textbf{0.025723} & \textbf{0.003077} & \textbf{0.000051} & $\mathbf{1.43\times 10^{-8}}$ \\ & AP & 0.669421 & 0.401033 & 0.138012 & 0.016978 & 0.000281 & $7.89\times 10^{-8}$ \\ & RGP & \textbf{0.212160} & 0.093205 & 0.029739 & 0.005010 & 0.000194 & $3.16\times 10^{-7}$ \\
        \midrule%hline
        \multirow{3}{*}{$10^{-3}$}  & \texttt{O-ConPF/GP} & \textbf{0.249625} & \textbf{0.131458} & \textbf{0.061872} & \textbf{0.027385} & \textbf{0.011482} & \textbf{0.004578} \\ & AP & 0.996008 & 0.990396 & 0.977011 & 0.945613 & 0.874839 & 0.729240 \\ & RGP & \textbf{0.249625} & 0.133542 & 0.066946 & 0.032679 & 0.015783 & 0.007529 \\
        \midrule%line
        \multirow{3}{*}{$10^{-4}$}& \texttt{O-ConPF/GP} & \textbf{0.249962} & \textbf{0.131849} & \textbf{0.062293} & \textbf{0.027820} & \textbf{0.011920} & \textbf{0.005009} \\ & AP & 0.999600 & 0.999035 & 0.997672 & 0.994391 & 0.986526 & 0.967855 \\ & RGP & \textbf{0.249962} & 0.133931 & 0.067365 & 0.033112 & 0.016221 & 0.007966 \\
        % \hline
        \bottomrule
    \end{tabular*}
\end{table}
\FloatBarrier

The BNB comparisons in~\cref{tab:bnb-vs-fdom-.1} show close
agreement with $h_\algB$ at the tested budgets, while
\cref{tab:compare-other-schedules} demonstrates the favorable
performance of $h_\algB$ and $h_\algC$ relative to existing
schedules in their respective metrics.
At $q=0$, $h_\algB$ and $h_\algC$ reduce to the schedules
in~\cite{zhang2024accelerated,grimmer2024composing}, which were
conjectured to be optimal among fixed-step memoryless gradient
descent schedules in their respective performance metrics.
These schedules match or surpass the best-known
schedules~\cite{taylor2017thesis,daccache2019performance,eloi2022worst,gupta2024branch,kamri2025numerical}.
These observations motivate
the following conjecture for $h_\algB$ and $h_\algC$ in their
respective metrics, extending the existing conjectures in~\cite{zhang2024accelerated,grimmer2024composing} to the
strongly convex setting.

\Needspace{8\baselineskip}
\begin{conjecture}
    Let $q\in[0,1)$ and $n\in\N$. The schedules $h_\algB^{(n)}$ and $h_\algC^{(n)}$ defined in~\cref{def:dynamic-programming} satisfy
    \[
  \begin{aligned}
  h_\algB^{(n)}&
\in \argmin_{h\in\R_+^n}
\sup_{\substack{
    d\in\N_+,\ L>0,\ f\in\mathcal F_{qL,L}(\R^d)\\
    \x_\star\in\argmin f,\ \x_0\notin\argmin f
}}
\frac{f_n-f_\star}
     {\frac{L}{2}\|\x_0-\x_\star\|^2},
    \\    h_\algC^{(n)}&\in \argmin_{h\in\R_+^n}
\sup_{\substack{
    d\in\N_+,\ L>0,\ f\in\mathcal F_{qL,L}(\R^d)\\
    \x_\star\in\argmin f,\ \x_0\notin\argmin f
}}\frac{\frac{1}{2L}\|\g_n\|^2}{f_0-f_\star}.
  \end{aligned}
    \]
\end{conjecture}

\section{Conclusion}
\label{sect:conclusion}
We have developed a curvature-dependent concatenation framework
for memoryless gradient descent. It yields the \OConPF{} and
\OConGP{} schedules with sharp worst-case guarantees and explicit bounds on iteration complexity, and also applies to smooth convex and smooth
weakly convex functions. We also propose the \OConPPr{} and \OConPPl{} schedules for distance and gradient norm, which extend the strongly convex silver schedule and its reversal, respectively, to arbitrary iteration budgets. Based on existing theoretical results and our numerical experiments, we conjecture that the \OConPF{} and \OConGP{} schedules are optimal in their respective metrics for every fixed $0\leq q<1$ and iteration budget.

\section*{Statements and Declarations}

\noindent\textbf{AI use disclosure}
\Cref{lemma:upper-bound-Rn,lemma:relation-phi-psi} were explored with the aid of ChatGPT. ChatGPT was also used in conducting the BNB experiment  in~\cref{sect:compare} and in revising or refining the writing.
The authors reviewed and revised the AI-assisted statements and proofs and take full responsibility for the manuscript.

\noindent\textbf{Competing interests} The authors have no conflict of interest to disclose.
\par

%\printbibliography
\bibliographystyle{abbrvnat}
 \bibliography{sn-bibliography}% common bib file
% \input sn-article.bbl

\begin{appendices}
\section{Proofs and Auxiliary Results for Section~\ref{sect:stepsize-shcedules}}\label{sect:proofs-section3}
\subsection{Lemmas of Algebraic Inequalities}\label{sect:algebraic-ineqs}
Recall that $\varphi_q$ and $\psi_q$ are defined in~\eqref{eq:phi} and~\eqref{eq:Psi}.
The next two lemmas characterize $\varphi_q$ and $\psi_q$ as roots of polynomials and establish the bounds needed to verify~\eqref{eq:sd-conditions}.

\begin{lemma}[Characterization of $\varphi_q$]\label{lemma:roots-of-phi}
Assume $q<1$, $x,y>-1$, $1+q x>0$, and $1+q y>0$. Then $\varphi_q(x,y)$ is the larger root of the quadratic polynomial
\begin{equation*}\label{eq:Phi}\tag{Def-$\Phi$}
    \Phi_{x,y,q}(h)=(1-q)h^2+(q(2xy+x+y+2)+x+y)h-((1+q)xy+2(x+y+1)).
\end{equation*}
Moreover, the following bounds hold.
\begin{enumerate}[label=(\roman*),ref=(\roman*),leftmargin=*]
    \item If $q\in (0,1)$, then $1< \a_1<\varphi_q(x,y)<\a_2<\min\{x+2,q^{-1}\}$, where
    \[\a_1=1+\frac{q(1-q)(x+1)(y+1)}{q^2(x+1)(y+1)+(qx+1)(qy+1)}\qquad \text{and}\qquad \a_2=\frac{qx+x+2}{2qx+q+1}.
    \]
    \item When $q\leq0$, $1<\varphi_q(x,y)<x+2$.
\end{enumerate} 
\end{lemma}
\begin{proof}
The fact that $\varphi_q(x,y)$ is the larger root of $\Phi_{x,y,q}(h)$ is trivial from the definition of~$\varphi_q$,~\eqref{eq:phi}.
% For $q<0$, $\Phi_{x,y,q}$ has positive leading coefficient and satisfies
% \[
% \Phi_{x,y,q}(1)=-(1-q)(x+1)(y+1)<0,\qquad
% \Phi_{x,y,q}(x+2)=2(x+1)^2(1+qy)>0.
% \]
% Hence its larger root $\varphi_q(x,y)$, given by~\eqref{eq:phi}, lies in $(1,x+2)$. It remains to consider $0\leq q<1$.
 
Note that
\begin{subequations}\label{eq:Phi-values}
\begin{align}
 % \lim_{h\to-\infty}  \Phi_{x,y,q}(h)&=+\infty \\ 
 \Phi_{x,y,q}(1)&=-(1-q)(x+1)(y+1)<0,\label{eq:Phi-at-one}
  \\ \Phi_{x,y,q}(\a_1)&=-\frac{(1-q)^3(x+1)(y+1)(qx+1)(qy+1)}{\bigl(q^2(x+1)(y+1)+(qx+1)(qy+1)\bigr)^2}<0,\label{eq:Phi-at-alpha1}
  \\ \Phi_{x,y,q}(\a_2)&= \frac{2(1-q)^2(x+1)^2(qx+1)}{(2qx+q+1)^2}>0,\qquad \text{if }2qx+q+1\neq 0,\label{eq:Phi-at-alpha2}
  \\ \Phi_{x,y,q}(x+2)&=2(x+1)^2(qy+1)>0.\label{eq:Phi-at-x-plus-two}
  %\\ \Phi_{x,y,q}(q^{-1})&=\frac{(1-q)(qx+1)(qy+1)}{q^2}>0,\qquad\text{if }q\neq0.
\end{align}
\end{subequations}

Since $\Phi_{x,y,q}$ is a quadratic polynomial with positive leading coefficient and $\Phi_{x,y,q}(1)<0$ by~\eqref{eq:Phi-at-one}, its smaller root lies below $1$.
The signs of $\Phi_{x,y,q}(1)$ and $\Phi_{x,y,q}(x+2)$ in~\eqref{eq:Phi-at-one} and~\eqref{eq:Phi-at-x-plus-two} thus give $1<\varphi_q(x,y)<x+2$.

Specifically, when $q\in(0,1)$, we have  $x+2-\a_2=\frac{2q (x+1)^2}{2qx+q+1 }>0$ and $\frac{1}{q}-\a_2=\frac{(1-q )(qx+1)}{q(2qx+q+1 )}>0$. It is also trivial that $\a_1>1$ and $\a_2= 1+\frac{(1-q)(x+1)}{2qx+q+1}>1$.
The signs at $\a_1,\a_2>1$ in~\eqref{eq:Phi-at-alpha1} and~\eqref{eq:Phi-at-alpha2} therefore give
\[
1<\a_1<\varphi_q(x,y)<\a_2<\min\{x+2,q^{-1}\}.
\]
\end{proof}

\begin{lemma}[Characterization of $\psi_q$]\label{lemma:roots-of-psi} Assume $q<1$, $x>-1$, $y\geq0$, $1+q x>0$, and $1+q y>0$.
    \begin{enumerate}[label=(\roman*),ref=(\roman*),leftmargin=*]
    \item\label{item:psi-1} 
    When $q\in (0,1)$, $\Psi_{x,y,q}(h) $ has three real roots $r_1<r_2<r_3$ and 
    \[
   r_1<1<r_2< \frac{qx+x+2}{2qx+q+1}<\min\{x+2,\frac{1}{q}\}\leq\max\{x+2,\frac{1}{q}\}<r_3.
    \]

    \item When $q =0$, $\Psi_{x,y,0}(h)$ is quadratic and has two real roots $r_1<r_2$, and
    \[
    r_1<1<r_2<x+2.\]
\item 
When $q<0$, $\Psi_{x,y,q}$ has three distinct real roots satisfying
\begin{equation*}
r_1<q^{-1}<r_2<1<r_3<x+2.
\end{equation*}

    \end{enumerate}
    Moreover, let $h_0$ be the unique root of $\Psi_{x,y,q}(h)=0$ in the interval $(1,x+2)$. Then $\Psi_{x,y,q}'(h_0)>0$.
\end{lemma}

\begin{proof}

% For $q<0$, the cubic $\Psi_{x,y,q}$ has positive leading coefficient $-q$, and direct evaluation gives
% \[
% \Psi_{x,y,q}(q^{-1})>0,\qquad \Psi_{x,y,q}(1)<0,\qquad \Psi_{x,y,q}(x+2)>0.
% \]
% The intermediate value theorem gives the three distinct roots in~\eqref{item:psi-2}. Thus $\psi_q(x,y)=r_3$ is the unique root in $(1,x+2)$, with $\Psi'_{x,y,q}(r_3)>0$. It remains to consider $0\leq q<1$.
 
Set $\a_2=\frac{qx+x+2}{2qx+q+1}=1+\frac{(1-q)(x+1)}{2qx+q+1}$ whenever $2qx+q+1\neq0$.  Note that
\begin{subequations}\label{eq:Psi-values}
\begin{align}
\Psi_{x,y,q}(1)&=-(1-q)(x+1)(qy^2+2y+1)<0,\label{eq:Psi-at-one}\\
\Psi_{x,y,q}(\a_2)&=\frac{2(1-q)^2(x+1)^2(qx+1)^2}{(2qx+q+1)^3}>0,\qquad \text{if }2qx+q+1>0,\label{eq:Psi-at-alpha2}\\
\Psi_{x,y,q}(x+2)&=2(x+1)^2(qy+1)^2>0,\label{eq:Psi-at-x-plus-two}\\
\Psi_{x,y,q}(q^{-1})&=\frac{(1-q)(qx+1)(qy+1)^2}{q^2}>0,\qquad \text{if }q\neq 0.\label{eq:Psi-at-q-inverse}
\end{align}
\end{subequations}

When $q\in(0,1)$, $\a_2>1$, $x+2-\a_2=\frac{2q (x+1)^2}{2qx+q+1}>0$ and $\frac{1}{q}-\a_2=\frac{(1-q )(qx+1)}{q(2qx+q+1 )}>0$ hold. Also, the leading coefficient of $\Psi_{x,y,q}$, $-q$, is negative when $q\in(0,1)$. Hence $\lim_{h\to-\infty}\Psi_{x,y,q}(h)=+\infty $ and $\lim_{h\to +\infty}\Psi_{x,y,q}(h)=-\infty$. Together with the signs in~\eqref{eq:Psi-values}, these limits prove~\cref{item:psi-1} by the intermediate value theorem.

When $q = 0$, $\Psi_{x,y,0}(h)=2h^2+(2y-3)h -2xy-x-4y$ and \eqref{eq:Psi-at-one} and~\eqref{eq:Psi-at-x-plus-two} give $r_1<1<r_2<x+2$.

When $q<0$, $\Psi_{x,y,q}$ has positive leading coefficient $-q>0$.
Hence $\lim_{h\to-\infty}\Psi_{x,y,q}(h)=-\infty $ and $\lim_{h\to +\infty}\Psi_{x,y,q}(h)=+\infty$.
Thus the signs at~\eqref{eq:Psi-at-q-inverse},~\eqref{eq:Psi-at-one} and~\eqref{eq:Psi-at-x-plus-two}
give $r_1<q^{-1}<r_2<1<r_3<x+2$ by the intermediate value theorem.

Since $h_0$ is a simple root at which $\Psi_{x,y,q}$
changes sign from negative to positive, we have
$\Psi_{x,y,q}'(h_0)>0$.\end{proof}

The following identities are used to verify the certificate decompositions in the concatenation proofs.
\begin{lemma}[Equivalent forms of $\varphi_q$ and $\psi_q$]\label{lemma:phi-psi-x} Assume $q<1$, $x,y\geq0$, $1+q x>0$, and $1+q y>0$. Let $\a=\varphi_q(x,y)$ and $\b=\psi_q(x,y)$. Then we have
\begin{subequations}
    \begin{align}\label{eq:phi-x}
        %  x &= \frac{(y+2-\a)\a(q-1)+2(\a-1)(y+1)}{ qy+y+2-(2qy+q+1)\a},\\
        x&=\a-2+\frac{2(\a-1)^2(qy+1)}{ qy+y+2-(2qy+q+1)\a},\\
        %  x &=\frac{\beta(q\beta^2-2(q+1)\beta+3)-y(q \b+\b-2)(q y+2)}{q \b^2-1+(2q\b-q -1)y(q y+2)}
        %  x &=\b - 2-\frac{2(\b-1)^2(qy+1)^2}{q \b^2-1+(2q\b-q -1)y(q y+2)}
                 x &=\b - 2+\frac{2(\b-1)^2(qy+1)^2}{(qy+1)^2(qy^2+2y+1)-q(qy^2+2y+\b)^2}
         \label{eq:psi-x}.
    \end{align}
\end{subequations}
Also, the denominators in~\eqref{eq:phi-x} and~\eqref{eq:psi-x} are nonzero.
\end{lemma}

\begin{proof}
Note that     $\Phi_{x,y,q}(\a) $ and $\Psi_{x,y,q}(\b) $ from~\eqref{eq:Phi} and~\eqref{eq:Psi} are both
linear in $x$.
%  For~\eqref{eq:phi-x}, 
Collecting the terms of $\Phi_{x,y,q}(\a)=0$ and $\Psi_{x,y,q}(\b)=0$ in $(x+2-\a)$ and $(x+2-\b)$,
respectively, gives
%  \begin{equation}
%     \label{eq:phi-x-proof}
% \bigl(qy+y+2-(2qy+q+1)\a\bigr)\cdot (x+2-\a)=2(\a-1)^2(qy+1).
%  \end{equation}
\begin{equation}\label{eq:phi-psi-x-proof}
    \begin{aligned}
        %  x &= \frac{(y+2-\a)\a(q-1)+2(\a-1)(y+1)}{ qy+y+2-(2qy+q+1)\a},\\
\bigl(qy+y+2-(2qy+q+1)\a\bigr)\cdot (x+2-\a)&=2(\a-1)^2(qy+1),\\
        %  x &=\frac{\beta(q\beta^2-2(q+1)\beta+3)-y(q \b+\b-2)(q y+2)}{q \b^2-1+(2q\b-q -1)y(q y+2)}
((qy+1)^2(qy^2+2y+1)-q(qy^2+2y+\b)^2)\cdot (x+2-\b)&=2(\b-1)^2(qy+1)^2.
        %  \label{eq:psi-x-proof2}.
    \end{aligned}
\end{equation}
Recalling $1<\a,\b<x+2$ by~\cref{lemma:roots-of-phi,lemma:roots-of-psi},~\eqref{eq:phi-psi-x-proof} thus implies
the denominators in~\eqref{eq:phi-x} and~\eqref{eq:psi-x} are both nonzero.
Solving for $x$ in~\eqref{eq:phi-psi-x-proof} gives the desired identities.
% The assumptions give $qx+1,qy+1>0$.
%  By~\cref{lemma:roots-of-phi}, $\alpha=\varphi_q(x,y)$ is the larger root of the quadratic polynomial $\Phi_{x,y,q}$ defined in~\eqref{eq:Phi}. Thus $\Phi_{x,y,q}(\a)=0$. Collecting the terms in $x$ gives
% \begin{equation}\label{eq:phi-x-proof}
% \Phi_{x,y,q}(\a)=-(qy+y+2-(2qy+q+1)\a)x+((y+2-\a)\a(q-1)+2(\a-1)(y+1)).
% \end{equation}
%  Since $\varphi_q(x,y)=\varphi_q(y,x)$, applying~\cref{lemma:condition-of-sd} with $x$ and $y$ interchanged yields $qy+y+2-(2qy+q+1)\a >0$. We can therefore solve~\eqref{eq:phi-x-proof} for $x$ to obtain~\eqref{eq:phi-x}.

%  For~\eqref{eq:psi-x}, collecting the terms of $\Psi_{x,y,q}(\b)=0$ in $(x+2-\b)$ gives
% \begin{equation}
%     \label{eq:psi-x-proof2}
% (q\b^2-1+(2q\b-q-1)y(q y+2))\cdot (x+2-\b)
% =-2(\b-1)^2(qy+1)^2.
% \end{equation}
% Recalling $1<\b<x+2$ by~\cref{lemma:roots-of-psi},~\eqref{eq:psi-x-proof2} thus 
% implies the denominator in~\eqref{eq:psi-x} is nonzero.
% Solving~\eqref{eq:psi-x-proof2} for $x$ gives~\eqref{eq:psi-x}.
\end{proof}

We next establish the monotonicity of the connecting stepsizes used in the dynamic-programming construction.
\begin{lemma}\label{lemma:monotonicity-of-phi-psi} 
    Let $q<1$, $x,y\geq0$, $1+q x>0$, and $1+q y>0$. Then 
    $\varphi_q(x,y)$ and $\psi_q(x,y)$ are monotonically increasing in both variables $x$ and $y$.
\end{lemma}
\begin{proof} 
%     For $\varphi_q$ with $q\geq0$, let $A=x+y+\frac{2q}{1-q}(x+1)(y+1)$ and $B=(A+2)^2+4(x+1)(y+1)$. Then $\varphi_q(x,y)=\frac{-A+\sqrt{B}}{2}$ and direct differentiation gives
%     \[
%     \begin{aligned}
%     \frac{\partial A}{\partial x}&=\frac{2qy+q+1}{1-q}>0,\qquad
%     \frac{\partial B}{\partial x}&=
%     \frac{2(1-q)(y+1)(2qy+3-q)+2(x+1)(2qy+q+1)^2}{(1-q)^2}>0.
%     \end{aligned}
%     \]
%     Also, direct computation shows that $\left(\frac{\partial B}{\partial x}\right)^2-4B\left(\frac{\partial A}{\partial x}\right)^2=\frac{32(qy+1)(y+1)^2}{1-q}> 0$. Since $B>0$ and both partial derivatives are positive, this gives $\frac{\partial B}{\partial x}>2\sqrt{B}\,\frac{\partial A}{\partial x}$. Hence
%     \[ 
% \frac{\partial}{\partial x}    \varphi_q(x,y)=-\frac{1}{2}\frac{\partial A}{\partial x}+\frac{1 }{4\sqrt{B}}\frac{\partial B}{\partial x}> 0.
%     \]
%     Note that $\varphi_q(x,y)$ is symmetric in $x$ and $y$, so it is monotonically increasing in both variables $x$ and $y$.

Set $\alpha=\varphi_q(x,y)$. Since $\alpha$ is the larger root of $\Phi_{x,y,q}$ defined in~\eqref{eq:Phi}, we have $\Phi'_{x,y,q}(\alpha)>0$. Moreover,
\[
\frac{\partial}{\partial x}\Phi_{x,y,q}(\alpha)
=(2qy+q+1)\alpha-(qy+y+2)<0,
\]
where the inequality follows from~\cref{lemma:condition-of-sd} with $x$ and $y$ interchanged.
Implicit differentiation gives $\frac{\partial}{\partial x}\varphi_q(x,y)>0$, and symmetry gives monotonicity in $y$.

For $\psi_q(x,y)$, define $\beta=\psi_q(x,y)$. By~\cref{lemma:roots-of-psi}, $\Psi'_{x,y,q}(\beta)>0$. Now we note that the following identities hold, which can be verified by expansion,
\begin{equation}\label{eq:diff-of-Psi}
\begin{aligned}    
\frac{\partial}{\partial x}\Psi_{x,y,q}(h)
&=\frac{q\Psi_{x,y,q}(h) - 2(1-q h)^2(qy^2+2y+1)}{(qx+1)+(1-qh)},\\
\frac{\partial}{\partial y}\Psi_{x,y,q}(h)
&=-2(qy+1)(qx+x+2-(2qx+q+1)h).
\end{aligned}
\end{equation}
By~\cref{lemma:condition-of-sd}, $1-q \beta >0$ and $qx+x+2-(2qx+q+1)\b>0$ hold. Since $\Psi_{x,y,q}(\beta)=0$,~\eqref{eq:diff-of-Psi} implies $\frac{\partial}{\partial x}\Psi_{x,y,q}(h)<0$ and $\frac{\partial}{\partial y}\Psi_{x,y,q}(h)<0$   when $h=\beta$.
The implicit function theorem therefore gives
\[
\frac{\partial\beta}{\partial x}
=-\frac{\partial_x\Psi_{x,y,q}(\beta)}{\Psi'_{x,y,q}(\beta)}>0,
\qquad
\frac{\partial\beta}{\partial y}
=-\frac{\partial_y\Psi_{x,y,q}(\beta)}{\Psi'_{x,y,q}(\beta)}>0.
\]
Thus, $\psi_q(x,y)$ is strictly increasing in both $x$ and $y$.
\end{proof}

The following comparison of the two connecting stepsizes is used to compare cumulative stepsizes in~\cref{prop:prim-asymp-2}.
\begin{lemma}\label{lemma:relation-phi-psi} Assume $q<1$, $x>-1$, $y\geq0$, $1+q x>0$, and $1+q y>0$. Then~$\varphi_q(x,y)< \psi_q(x,y)$. If $q\in (0,1)$, we also have
\begin{equation}\label{eq:compare-phi-psi-2}
\frac{(qx+1)(q \cdot\frac{2}{1+q}y+1)}{1-q\psi_q(x,\frac{2}{1+q}y)}-1< \frac{2}{1+q}\cdot\left(\frac{(qx+1)(qy+1)}{1-q\varphi_q(x, y)}-1\right).
\end{equation}
\end{lemma}

\begin{proof}
 Recall the polynomials $\Psi_{x,y,q}$ and $\Phi_{x,y,q}$ defined in~\eqref{eq:Psi} and~\eqref{eq:Phi}, respectively. The larger root of $\Phi_{x,y,q}$ is $\varphi_q(x,y)$. Note that
\begin{equation}\label{eq:pseudo-resultant-phi-psi}
\Psi_{x,y,q}(h)-(qy+1)\Phi_{x,y,q}(h)=-(1-q h)(x+2-h)((1-q)y+(h-1)).
\end{equation}
Taking $h=\varphi_q(x,y)$ in~\eqref{eq:pseudo-resultant-phi-psi} and using~\cref{lemma:roots-of-phi} give~$\Psi_{x,y,q}(\varphi_q(x,y))<0$. Moreover, by~\cref{lemma:roots-of-phi} and the proof of~\cref{lemma:roots-of-psi}, we have $1<\varphi_q(x,y)<x+2$ and $0<\Psi_{x,y,q}(x+2)$.  By continuity, $\Psi_{x,y,q}$ has a root in $(\varphi_q(x,y),x+2)$. Since $\psi_q(x,y)$ is its unique root in $(1,x+2)$ by~\cref{lemma:roots-of-psi}, we obtain $\varphi_q(x,y)<\psi_q(x,y)<x+2$.

 To prove~\eqref{eq:compare-phi-psi-2}, set $c_q=2/(1+q)$ and recall $\a_1,\a_2$ from~\cref{lemma:roots-of-phi}, with $\a_2$ independent of $y$. By~\cref{lemma:roots-of-psi,lemma:roots-of-phi}, we have $\psi_q(x,c_q y)<\a_2<q^{-1}$ and $\a_1<\varphi_q(x,y)<q^{-1}$. Hence
\[
\begin{aligned}
\frac{(qx+1)(q c_q y+1)}{1-q\psi_q(x,c_q y)}-1
&<\frac{(qx+1)(q c_q y+1)}{1-q\a_2}-1
\\ &
= c_q\left(\frac{(qx+1)(qy+1)}{1-q\a_1}-1\right)
% \\ 
< c_q\left(\frac{(qx+1)(qy+1)}{1-q\varphi_q(x, y)}-1\right).
\end{aligned}
\]
The equality follows by substituting the expressions for $\a_1$ and $\a_2$.\end{proof}

\subsection{\texorpdfstring{Proofs of \Cref{lemma:huber-gd,lemma:condition-of-sd,lemma:dom-u}}{Proofs of Lemmas \ref{lemma:huber-gd}, \ref{lemma:condition-of-sd}, and \ref{lemma:dom-u}}}\label{sect:proofs-main-lemmas}

\begin{proof}[Proof of~\Cref{lemma:huber-gd}]

% If $x_0=0$, then $r=0$ and $f(x)=\mu x^2/2$, so every iterate is zero and the claims follow.
% We therefore assume $x_0>0$.
Since $s\geq1$ and $\sigma^{(0)}=0$, we have $x_0\geq r$ and~\eqref{eq:huber-xk} holds for $k=0$.

Suppose $x_k\geq r$ and~\eqref{eq:huber-xk} holds for some $0\leq k<n$. Then $f'(x_k)=\mu x_k+(L-\mu)r$, so
\begin{equation}\label{eq:huber-f-xk-induction1}
x_{k+1}=(1-q h_k)x_k-(1-q)h_kr
=\left((1-q h_k)\frac{1-\frac{(1-q)\sigma^{(k)}}{s}}{1+q\sigma^{(k)}}-\frac{(1-q)h_k}{s}\right)x_0.
\end{equation}
By~\eqref{eq:qs-add-1}, $1+q\sigma^{(k+1)}=(1+q\sigma^{(k)})/(1-q h_k)$. Also,~\eqref{eq:telescope} gives
\begin{equation}\label{eq:sigma-add-h}
\sigma^{(k)}+h_k(1+q\sigma^{(k+1)})
=\sigma^{(k)}+h_k\prod_{i=0}^k(1-q h_i)^{-1}
=\sigma^{(k+1)}.
\end{equation}
Substituting these identities into~\eqref{eq:huber-f-xk-induction1} yields
\begin{equation}\label{eq:huber-f-xk-induction2}
\begin{aligned}
x_{k+1}
&=\frac{1-\frac{(1-q)(\sigma^{(k)}+h_k(1+q\sigma^{(k+1)}))}{s}}{1+q\sigma^{(k+1)}}x_0\overset{\eqref{eq:sigma-add-h}}{=}\frac{1-\frac{(1-q)\sigma^{(k+1)}}{s}}{1+q\sigma^{(k+1)}}x_0.
\end{aligned}
\end{equation}
By~\cref{lemma:s-bound}, padding the prefix with zeros gives $\sigma^{(k+1)}\leq\sigma(h)$, and $1+q\sigma^{(k+1)}>0$. Therefore,
\begin{equation}\label{eq:huber-f-xk-induction3}
\frac{1-\frac{(1-q)\sigma^{(k+1)}}{s}}{1+q\sigma^{(k+1)}}
=\frac1s+\frac{s-(1+\sigma^{(k+1)})}{s(1+q\sigma^{(k+1)})}
\geq\frac1s,
\end{equation}
where the inequality uses $s\geq1+\sigma(h)$. Combining~\eqref{eq:huber-f-xk-induction2} and~\eqref{eq:huber-f-xk-induction3} gives $x_{k+1}\geq x_0/s=r$, completing the induction.\end{proof}

\begin{proof}[Proof of~\Cref{lemma:condition-of-sd}]
Applying~\cref{lemma:roots-of-phi} when $\alpha=\varphi_q(s,t)$
and~\cref{lemma:roots-of-psi} when $\alpha=\psi_q(s,t)$,
we obtain $\alpha>1$ and $1-q\alpha>0$.
% For $q>0$, these lemmas also give
% $
% \alpha<\frac{q s+s+2}{2q s+q+1},
% $
% whose denominator is positive. This proves the last inequality. When $q=0$, it reduces to $\alpha<s+2$, which follows from the same lemmas.

% For $q<0$, set $D=2qs+q+1$. If $D\leq0$, then $qs+s+2-D\alpha>0$ follows from $qs+1>0$, $s+1>0$, and $\alpha>1$. If $D>0$, let $\alpha_2=(qs+s+2)/D>1$. By~\eqref{eq:Phi-at-alpha2} and~\eqref{eq:Psi-at-alpha2} with $(x,y)=(s,t)$, both polynomials are positive at $\alpha_2$. The root orderings in~\cref{lemma:roots-of-phi,lemma:roots-of-psi} then give $\varphi_q(s,t),\psi_q(s,t)<\alpha_2$, and hence $qs+s+2-D\alpha>0$.

Let $D=2qs+q+1$. When $D>0$, let $\a_2=\frac{qs+s+2}{D}$. Then 
\[\a_2=1+\frac{(1-q)(s+1)}{D}=s+2-\frac{2q(s+1)^2}{D}.\]
If $q\leq0$, this gives $\a_2\geq s+2>\a$ by~\cref{lemma:roots-of-phi,lemma:roots-of-psi}.
If $q>0$, the same lemmas give $\a<\a_2$ directly.
Thus $qs+s+2-D\a=D(\a_2-\a)>0$ in both cases.
When $D\leq 0$, the desired inequality $qs+s+2-D\a> 0$ automatically holds because $qs+s+2>0$.\end{proof}

\begin{proof}[Proof of~\Cref{lemma:dom-u}]
By~\cref{lemma:condition-of-sd}, $\b>1$ and $1-q\b>0$. Using~\eqref{eq:sa} and $s_a\geq0$, we obtain
\begin{equation*}
v - w=\frac{(s_a+1)^2}{(\b-1)^2}-\frac{(s_a+1)(q s_a+2-q\b)}{(1-q\b)(1+2s_b+q s_b^2)}=\frac{(s_a+1)(1+q s_a)}{(\beta-1)(1-q \b)}\geq 0.
\end{equation*}
Since $\gamma\geq0$, $u_i\geq0$ for $i=0,\dotsc,n-1$, and $u_b\geq0$, all entries of $u_c$ are nonnegative.

Since $u_b$ is a certificate for the \dom property of $h_b$, we have $(1+2s_b+q s_b^2)-(1-q)(\onet u_b)\geq 0$. To prove the corresponding inequality for $u_c$, we compute
\begin{equation}\label{eq:proof-lemma-dom-u-sc}
\begin{aligned}
    \onet u_c &= \sum_{j=0}^{n-1}u_j + \gamma +(v-w)(\onet u_b)
 \stackrel{\eqref{eq:telescope}}{=} s_a+\gamma +(v-w)(\onet u_b),\\ 
s_c &= \sigma(h_c) \stackrel{\eqref{eq:recursion-of-sigma}}{=\!=} \frac{s_a+\beta+s_b+q s_as_b}{1-q\beta}.
\end{aligned}
\end{equation}
Therefore,
\begin{equation*}
1+2s_c+q s_c^2-(1-q)(\onet u_c)\geq 
1+2s_c+q s_c^2-(1-q)(s_a+\gamma) - (v-w)(1+2s_b+q s_b^2)=0,
\end{equation*}
where the last equality follows by substituting $s_c$ from~\eqref{eq:proof-lemma-dom-u-sc} and then $s_a$ from~\eqref{eq:sa}.\end{proof}

\section{Auxiliary Results for Section~\ref{sect:convergence}}
\label{sect:lemma-convergence}
This appendix establishes the auxiliary estimates for the exponent and remainder term used in the proofs of~\cref{prop:prim-asymp,prop:nu-bound}.
We first analyze the scalar recurrence obtained by repeatedly concatenating two identical copies of a stepsize schedule.

\begin{lemma}\label{lemma:F-properties}
Define
    $G(x)=\frac{x^2+1+(x-1)\sqrt{x^2+2x+5}}{2}$. Let $x\geq 1$.
    Then we have
    \begin{enumerate}[label=(\roman*),ref=(\roman*),leftmargin=*]
        \item $\frac{G(x)}{x^2}\ge1 $ ,
        \item $\frac{G(x)}{x^2} \leq
         \min\left\{ 1+\frac{x-1}{x^2(x+1)},\ \frac{27}{25}\right\} $ and $\frac{G(x)}{x^2}$ increases over $[1,\frac{5}{3}]$, 
        \item $
1 +  (\sqrt2+1)(x-1)\leq   G(x) 
        \leq
        1+(\sqrt2+1)(x-1)+(x-1)^2$, and
        \item $\frac{xG'(x)}{G(x)}\geq \frac{50}{27}$.
    \end{enumerate}
\end{lemma}

\begin{proof}
Since $\sqrt{x^2+2x+5}\geq x+1$ for $x\geq1$, we have
\[
0\leq G(x)-x^2
=\frac{2(x-1)}{\sqrt{x^2+2x+5}+x+1}
\leq\frac{x-1}{x+1}.
\]
Dividing by $x^2$ gives~(i) and the first bound in~(ii). Differentiation also gives
\[
\frac{\mathrm d}{\mathrm dx}\left(\frac{G(x)}{x^2}\right)
=\frac{5-x-\sqrt{x^2+2x+5}}{x^3\sqrt{x^2+2x+5}}.
\]
The numerator decreases strictly and is zero at $x=5/3$, so $G(x)/x^2$ increases on $[1,5/3]$ and decreases thereafter. Its maximum is $G(5/3)/(5/3)^2=27/25$, completing~(ii).

For~(iii), we use
$
2\sqrt2\leq\sqrt{x^2+2x+5}\leq x-1+2\sqrt2,
$
where the upper bound follows from
\[
(x-1+2\sqrt2)^2-(x^2+2x+5)=4(\sqrt2-1)(x-1)\geq0.
\]
Substituting these bounds on \(\sqrt{x^2+2x+5}\) and using \(x+1\ge2\), we obtain
\[
\begin{aligned}
(1+\sqrt2)(x-1)
&\leq G(x)-1
=\frac{x-1}{2}\bigl(x+1+\sqrt{x^2+2x+5}\bigr)\leq(1+\sqrt2)(x-1)+(x-1)^2.
\end{aligned}
\]

For~(iv), differentiation gives
\[G'(x)=x+\frac{x^2+x+2}{\sqrt{x^2+2x+5}}
=x+\frac{\sqrt{x^2(x^2+2x+5)+4(x+1)}}{\sqrt{x^2+2x+5}}\geq 
2x.
\]
This combined with~(ii) gives~(iv).\end{proof}

\begin{lemma}\label{lemma:duplicate}
    Let $a_0 >1$ and $a_{n+1}=G(a_n)$ where $G(x)=\frac{x^2+1+(x-1)\sqrt{x^2+2x+5}}{2}$. Then the following statements hold:
    \begin{enumerate}[label=(\roman*),ref=(\roman*),leftmargin=*]
    \item\label{item:duplicate-limit} The limit $C\triangleq\lim_{n\to\infty}2^{-n}\log a_n$ exists and is finite.
    \item\label{item:duplicate-upper} $\log a_n-C\cdot 2^n\leq 0$,
        \item\label{item:duplicate-sum}    $-\frac{1}{2}<-\frac{1}{a_0(a_0+1)}
    \leq \sum_{n=0}^\infty( \log a_n-C\cdot 2^n)\leq 0$, and
\item\label{item:duplicate-rate}  $\log a_0\leq C\leq \log a_0+\log \frac{27}{25}$. Moreover, $0.7446(a_0-1)^{\rho}\leq C\leq 0.9526(a_0-1)^{\rho}$ holds if $1<a_0\leq 8/5$.
    \end{enumerate}
\end{lemma}
\begin{proof}
\textit{\ref{item:duplicate-limit}}
% By~\cref{lemma:F-properties}, $1\leq \frac{G(x)}{x^2}\leq 1+\min\{(\sqrt 2-1)(x-1),\frac{1}{x^2}\}$ for $x>0$.
The recurrence gives
\begin{equation}
    \frac{\log a_{k+1}}{2^{k+1}}-\frac{\log a_{k}}{2^{k}}
  =\frac{1}{2^{k+1}}\log \frac{G(a_{k})}{a_{k}^2}.\label{eq:diff-log-an-div-2pow}
\end{equation}
By~\cref{lemma:F-properties}, $0\leq\log(G(a_k)/a_k^2)\leq\log(27/25)$. Equation~\eqref{eq:diff-log-an-div-2pow} also shows that $2^{-n}\log a_n$ is nondecreasing. Summing~\eqref{eq:diff-log-an-div-2pow} from $k=0$ to $N-1$ gives
\[
0\leq\frac{\log a_N}{2^N}-\log a_0
\leq\sum_{k=0}^{N-1}\frac{1}{2^{k+1}}\log\frac{27}{25}
=(1-2^{-N})\log\frac{27}{25}.
\]
Thus the sequence is bounded above and converges to a finite limit $C$.

\textit{\ref{item:duplicate-upper}}
Summing~\eqref{eq:diff-log-an-div-2pow} from $k=n$ to $N-1$, letting $N\to\infty$, and multiplying by $-2^n$, we obtain
\begin{equation}
    \label{eq:err-an-c2n}
\log a_n - C\cdot 2^n=-2^n\sum_{k=n}^\infty \frac{1}{2^{k+1}}\log \frac{G(a_k)}{a_k^2}
=-\sum_{k=n}^\infty \frac{1}{2^{k+1-n}}\log \frac{G(a_k)}{a_k^2},
\end{equation} 
which implies $\log a_n - C\cdot 2^n\leq 0$ directly.

\textit{\ref{item:duplicate-sum}}
Since $\log(G(a_k)/a_k^2)\geq0$ by~\cref{lemma:F-properties}(i), Tonelli's theorem permits interchanging the sums of these nonnegative terms. Summing~\eqref{eq:err-an-c2n} over $n\geq0$ therefore gives
\begin{equation}
\begin{aligned}\sum_{n=0}^\infty (\log a_n - C\cdot 2^n)=
-\sum_{k=0}^\infty\log \frac{G(a_k)}{a_k^2} \sum_{n=0}^k\frac{1}{2^{k+1-n}}
\geq - \sum_{k=0}^\infty\log \frac{G(a_k)}{a_k^2}.
% \geq 
\label{eq:sum-log-minus-C}
  \end{aligned}
\end{equation}
By~\cref{lemma:F-properties}(ii), and the inequality $\log(1+u)\leq u$ for $u\geq0$, we have, for $x\geq1$,
\[
\log\frac{G(x)}{x^2}
\leq\log\left(1+\frac{x-1}{x^2(x+1)}\right)
\leq\frac{x-1}{x^2(x+1)}
=\frac{1}{x+1}\left(\frac{1}{x}-\frac{1}{x^2}\right).
\]
Applying this bound with $x=a_k$ and summing over $k\geq0$ gives

\begin{equation}
\sum_{k=0}^\infty\log \frac{G(a_k)}{a_k^2}
\leq  
\sum_{k=0}^\infty\frac{1}{a_k+1}\left(\frac{1}{a_k}-\frac{1}{a_k^2}\right)\leq 
\sum_{k=0}^\infty\frac{1}{a_0+1}\left(\frac{1}{a_k}-\frac{1}{a_{k+1}}\right)=\frac{1}{a_0(a_0+1)},
\label{eq:sum-log-minus-C2}
\end{equation}
where we used $a_k\geq a_0>1$ and $a_{k+1}=G(a_k)\geq a_k^2$, both consequences of~\cref{lemma:F-properties}(i).
Combining~\eqref{eq:sum-log-minus-C} and~\eqref{eq:sum-log-minus-C2} proves Item~\ref{item:duplicate-sum}.

\textit{\ref{item:duplicate-rate}}
Rearranging~\eqref{eq:err-an-c2n}, we obtain
\begin{equation}
\label{eq:C-series-expansion}
C=2^{-n}\left(\log a_n+\sum_{k=n}^\infty \frac{1}{2^{k+1-n}}\log \frac{G(a_k)}{a_k^2}\right),\qquad \forall n\in\N.
\end{equation}
By~\cref{lemma:F-properties}, $1\leq \frac{G(z)}{z^2}\leq \frac{27}{25}$ holds and thus~\eqref{eq:C-series-expansion} immediately gives $\log a_0\leq C\leq \log a_0+\log \frac{27}{25}$ by specifying $n=0$.

For brevity, write
\[
r:=\sqrt2+1,\qquad c:=1-\frac1{\sqrt2},\qquad
\ell:=\frac{25+3\sqrt{65}}{40}.
\]

Now assume $1<a_0\leq8/5$. By~\cref{lemma:F-properties}(i), $a_n\geq a_0^{2^n}\to\infty$. Hence there is a largest integer $N$ such that $a_N\leq8/5$. Set $x=a_N$. We first prove the lower bound. By~\cref{lemma:F-properties}(iii), we have $a_{n+1}-1=G(a_n)-1\geq r(a_n-1)$ and thus $x-1\geq   r^{N}(a_0-1)$, i.e., $N\leq\log_{r}\frac{x-1}{a_0-1}$. 

Since $G$ is strictly increasing on $[1,\infty)$, the choice of $N$ and $G(\ell)=\frac85<a_{N+1}=G(x)$ imply $\ell<x\leq\frac85$. Then we obtain 
\[\log \frac{G(x)}{x^2}\geq \log \frac{G(\ell)}{\ell^2}=2\log (10\sqrt {10}-6\sqrt{26})\] 
by the monotonicity of $\frac{G(z)}{z^2}$ (\cref{lemma:F-properties}(ii)). Using~\eqref{eq:C-series-expansion} with $n=N$, retaining the $k=N$ term and discarding the remaining nonnegative terms, and using $N\leq\log_r\frac{x-1}{a_0-1}$, we obtain
\[
\begin{aligned}
C&=2^{-N}\left(\log x+\sum_{k=N}^\infty\frac{1}{2^{k+1-N}}\log \frac{G(a_k)}{a_k^2}\right)\geq \frac{\log x+\log(10\sqrt {10}-6\sqrt{26})}{2^{\log _{r}(x-1)}}(a_0-1)^{\rho}.
% \geq0.7446(a_0-1)^\rho.
\end{aligned}
\]
Differentiating the quotient shows that it increases and then decreases on $[\ell,\frac85]$. Its minimum is therefore the smaller endpoint value, attained at $x=\frac85$ and lying in $(0.7446064,0.7446065)$,
which is greater than $0.7446$. This proves the lower bound for $C$.

It remains to prove the upper bound for $C$. From~\cref{lemma:F-properties} we have  $a_{n+1}-1\leq r(a_n-1)+(a_n-1)^2$ and thus $\frac{a_{n+1}-1}{r^{n+1}}\leq \frac{a_n-1}{r^n}+\frac{(a_n-1)^2}{r^{n+1}}$ for $0\leq n\leq N-1$.   
Telescoping summation   yields
\begin{equation}
   \label{eq:bound-seq-lower-telescope1}
   \frac{x-1}{r^{N }} \leq  a_0-1+\sum_{k=0}^{N-1} \frac{(a_k-1)^2}{r^{k+1}}\leq a_0-1+\sum_{k=0}^{N-1} \frac{r^{2(k-N)}(x-1)^2}{r^{k+1}},
\end{equation}
where the last inequality follows again from $a_k-1\leq r^{k-N}(x-1)$.
Multiplying~\eqref{eq:bound-seq-lower-telescope1} by $r^{N }$ and summing up a geometric series gives
\begin{equation*}
   x-1\leq r^{N }(a_0 - 1)+\sum_{k=0}^{N-1}r^{k-N-1}(x-1)^2\leq r^{N }(a_0-1)+c(x-1)^2 .\end{equation*}
Since $1<x\leq8/5$ and $0<c<1$, we have $x-1-c(x-1)^2=(x-1)(1-c(x-1))>0$. Taking logarithms therefore gives
\begin{equation}\label{eq:duplicate-N-lower-bound}
N\geq\frac{\log\left(\frac{x-1-c(x-1)^2}{a_0-1}\right)}{\log r}.
\end{equation}

Recall $\ell\leq x\leq\frac85$. By monotonicity of $G(z)/z^2$ we have $\log \frac{G(x)}{x^2}\leq \log\frac{G(8/5)}{(8/5)^2}$. Also,~\cref{lemma:F-properties}(ii) guarantees that $\frac{G(z)}{z^2}\leq \frac{27}{25}$. Therefore we obtain
\[
\begin{aligned}
    C&\overset{\eqref{eq:C-series-expansion}}{=}2^{-N}\left(\log x+\sum_{k=N}^\infty \frac{1}{2^{k+1-N}}\log\frac{G(a_k)}{a_k^2}\right)
\leq 
2^{-N}\left(\log x+\frac{\log \frac{G(8/5)}{(8/5)^ 2}}{2}+\sum_{k=N+1}^\infty\frac{\log\frac{27}{25}}{2^{k+1-N}} \right)
\\ &=2^{-N}\left(\log x+\frac{1}{2}\log \frac{27(89+3\sqrt{269})}{3200}\right)
\overset{\eqref{eq:duplicate-N-lower-bound}}{\leq} \frac{\log x+\frac{1}{2}\log \frac{27(89+3\sqrt{269})}{3200}}{2^{ \frac{\log(x - 1 - c(x-1)^2)}{\log r}}}(a_0 - 1)^{\rho}.
\end{aligned}
\]
Differentiation shows that the quotient in the last line first decreases and then increases on $[\ell,\frac85]$, so its maximum is attained at an endpoint. Comparing the endpoint values gives the maximum at $x=\ell$, with value in $(0.9525906,0.9525907)$, which is less than $0.9526$. Hence the proof is complete. \end{proof}

We next consider the bivariate map $F$ in~\eqref{eq:concat-rate-eqw}, for which $F(x,x)=G(x)$. The following two inequalities are used to control maximization over all splits.
\begin{lemma}\label{lemma:F-lemma-four-points} Let $x\geq 1, y\geq 1$ and define
     $F(x,y) = \frac{xy+1+\sqrt{(xy+3)^2-4(x+1)(y+1)}}{2}$. Then we have
 \begin{enumerate}[label=(\roman*),ref=(\roman*),leftmargin=*]
     \item\label{item:F-lemma-cauchy1}
     $F(x,x)F(y,y)\geq F(x,y)^2$, and
     \item\label{item:F-lemma-cauchy2} $F\left(F(x,x),F(y,y)\right)\geq F(F(x,y), F(x,y))$.
 \end{enumerate}
The second inequality is strict when $x,y>1$ and $x\ne y$.

 \begin{proof}
 Since $x,y\geq 1$, there exist $r,s\geq \sqrt 2$ such that $x=r-\frac{1}{r+1}$ and $y = s - \frac{1}{s+1}$. Plugging them into $F(x,x)$ and $F(y,y)$ yields $F(x,x)=r^2-1$ and $F(y,y)=s^2-1$.  
Noting that $F(x,y)\geq (xy+1)/2$ and it satisfies the quadratic equation
\begin{equation}\label{eq:F-lemma-defeq}
F(x,y)^2-(xy+1)F(x,y)+x+y-1=0,
\end{equation}
we claim that
\begin{equation}
    \label{eq:F-lemma-plug-in-gm}
F(x,x)F(y,y)-(xy+1)\sqrt{F(x,x)F(y,y)}+x+y-1\geq 0.
\end{equation}
Plugging in $x=r-\frac{1}{r+1}$, $y=s-\frac{1}{s+1}$, $F(x,x)=r^2-1$ and $F(y,y)=s^2-1$, we have
\begin{equation}
\begin{aligned}
    (F(x,x)F(y,y)+x+y-1)^2 - (xy+1)^2F(x,x)F(y,y)=\frac{\left(r-s\right)^{2}P(r,s)}{\left(r+1\right)^{2}\left(s+1\right)^{2}}
    \label{eq:F-lemma-plug-in-rs},
\end{aligned}
\end{equation}
where there is the identity
\begin{equation}\label{eq:Prs-sos}\begin{aligned}
P(r,s)%&=r^{4}s^{4}-3r^{4}s^{2}+2r^{4}-3r^{2}s^{4}+10r^{2}s^{2}+2r^{2}s-5r^{2}\\
% &\qquad+2rs^{2}+6rs+4r+2s^{4}-5s^{2}+4s+8
% \\ &=(r^2-2)(s^2-2)(r+s+1)^2+(r^2s^2-2r^2-2s^2-rs-r-s+2)^2\geq 0.
=(r^2-2)(s^2-2)(r^2-1)(s^2-1)+(rs+r+s+2)^2\geq 0.
\end{aligned}
\end{equation}
The nonnegativity of~\eqref{eq:Prs-sos} uses the fact that $r,s\geq \sqrt 2$.
As a result, the left-hand side of~\eqref{eq:F-lemma-plug-in-rs} is nonnegative and hence~\eqref{eq:F-lemma-plug-in-gm} is proved. Moreover,
\[
\sqrt{F(x,x)F(y,y)}\geq \sqrt{\frac{x^2+1}{2}\cdot \frac{y^2+1}{2}}\geq \frac{xy+1}{2}.\]
The polynomial $Q(z)=z^2-(xy+1)z+x+y-1$ is increasing on $[(xy+1)/2,\infty)$ and has $F(x,y)$ as its larger root by~\eqref{eq:F-lemma-defeq}. Thus~\eqref{eq:F-lemma-plug-in-gm} gives $\sqrt{F(x,x)F(y,y)}\geq F(x,y)$, proving~\cref{item:F-lemma-cauchy1}.

%---------------------------------------- item 2 ---------------------------------------- 

We next prove~\cref{item:F-lemma-cauchy2} of the lemma. We only need to deal with the case  $x,y>1$ so that the cases $x=1$ or $y=1$ follow from continuity. 
We denote $u=F(x,x) >1$, $v=F(y,y)> 1$ and $w=F(x,y)> 1$ for notational simplicity. Applying~\eqref{eq:F-lemma-defeq} to $(x,x)$, $(y,y)$, and $(x,y)$, we obtain $e_1=e_2=e_3=0$, where
\begin{align}
    e_1  &\triangleq u^2-(x^2+1)u+2x-1,\notag\\
    e_2  &\triangleq v^2-(y^2+1)v+2y-1,\notag\\
    e_3  &\triangleq w^2-(xy+1)w+x+y-1=0.\label{eq:F-lemma-e3}
\end{align}
Also, \(F(u,v)\) and \(F(w,w)\) are the larger roots of the equations \(e_4(z)=0\) and \(e_5(z)=0\), respectively, where we define
\[
\begin{aligned}
    e_4(z)&\triangleq z^2-(uv+1)z+u+v-1,\\
    e_5(z)&\triangleq z^2-(w^2+1)z+2w-1.\\
\end{aligned}
\]
First we claim that the resultant of $e_4$ and $e_5$ is nonpositive.
% \begingroup
% \newcommand{\deletedmath}[1]{\text{\st{\mbox{$\displaystyle #1$}}}}
% \st{In fact, the resultant \mbox{$\mathcal R\triangleq \text{Resultant}_z(e_4,e_5)$} is given by}
% \[
%     \begin{aligned}
%         &\deletedmath{\mathcal R=2u^{2}v^{2}w-u^{2}v^{2}-u^{2}vw^{2}-u^{2}v+u^{2}-uv^{2}w^{2}-uv^{2}-2uvw^{3}+2uvw^{2}}\\
%         &\deletedmath{\qquad+2uvw+2uv+uw^{4}+uw^{2}-4uw+v^{2}+vw^{4}+vw^{2}-4vw-w^{4}-2w^{3}+4w^{2},}
%     \end{aligned}
% \]
% \st{and there exists three polynomials of degree \mbox{$12$}, \mbox{$P_1,P_2,P_3\in \mathbb Q[u,v,w,x,y,z]$}, such that}
% \begin{equation}\label{eq:F-lemma-resultant-sos}
% \deletedmath{P_4^2+4(u-1)(v-1)(w^2-1)^2\cdot\mathcal R=P_1e_1+P_2e_2+P_3e_3,}
% \end{equation}
% \st{where}
% \[
% \begin{aligned}
%     &\deletedmath{P_4=u^{2}v^{2}w^{2}-u^{2}vw^{2}-u^{2}v-u^{2}w^{2}+2u^{2}w-uv^{2}w^{2}-uv^{2}-uvw^{4}}\\
%     &\deletedmath{\qquad-2uvw^{3}+4uvw^{2}+2uvw+uv+2uw^{4}-uw^{2}-4uw+u}\\
%     &\deletedmath{\qquad-v^{2}w^{2}+2v^{2}w+2vw^{4}-vw^{2}-4vw+v-2w^{4}-2w^{3}+7w^{2}-2w.}
% \end{aligned}
% \]
% \st{The exact values of \mbox{$P_1,P_2,P_3$} are presented in the supplementary materials.}\st{\mbox{\textbf{TODO:} add}} \st{We remark that the left hand side of \mbox{\eqref{eq:F-lemma-resultant-sos}} is actually (up to a constant factor) an element in the Gr\"{o}bner basis of \mbox{$e_1,e_2,e_3$}.}
% \st{Since \mbox{$e_1=e_2=e_3=0$}, \mbox{\eqref{eq:F-lemma-resultant-sos}} immediately implies \mbox{$\mathcal R\leq 0$} as \mbox{$u,v,w>1$}.}
% \endgroup

Set $A=uv-w^2$ and $B=u+v-2w$. Since $e_4(z)-e_5(z)=B-Az$, their resultant is
\[
\mathcal R\triangleq\operatorname{Resultant}_z(e_4,e_5)
=B^2-(w^2+1)AB+(2w-1)A^2.
\]
With $r,s\geq\sqrt 2$ as in the proof of~\cref{item:F-lemma-cauchy1}, we substitute $u=r^2-1$ and $v=s^2-1$ into the expression for $\mathcal R$ above. Using~\eqref{eq:F-lemma-e3} with $x=r-1/(r+1)$ and $y=s-1/(s+1)$
and defining $E_3=(r+1)(s+1)e_3=0$, we have
the polynomial identity
\begin{equation}\label{eq:F-lemma-resultant-square}
\mathcal R=E_3\cdot\bigl[E_3+2(r+s)(w-1)(rs-w-1)\bigr]-(r^2-2)(s^2-2)(w-1)^2(rs-w-1)^2\leq 0.
% =-\frac{(u-1)(v-1)}{(r+1)^2(s+1)^2}
% (K-Hw)^2\leq0.
\end{equation}

Let $\alpha,\beta$ be the roots of $e_4(z)=0$, and let $\gamma,\delta$ be the roots of $e_5(z)=0$, ordered so that $\alpha<\beta$ and $\gamma<\delta$. Since $e_4(1)=-(u-1)(v-1)<0$ and $e_5(1)=-(w-1)^2<0$, we have
\begin{equation}\label{eq:F-lemma-root-order}
\alpha<1<\beta=F(u,v),\qquad \gamma<1<\delta=F(w,w).
\end{equation}

Since $e_4$ and $e_5$ are monic, the root-product formula for the resultant and~\eqref{eq:F-lemma-resultant-square} give
\[
\mathcal R=e_5(\alpha)e_5(\beta)
=(\alpha-\gamma)(\alpha-\delta)(\beta-\gamma)(\beta-\delta)\leq0.
\]
Because $\alpha<1<\delta$ and $\gamma<1<\beta$, we have $(\alpha-\delta)(\beta-\gamma)<0$. Dividing by this negative factor yields $(\alpha-\gamma)(\beta-\delta)\geq0$.
Moreover, by Vieta's formulas and~\cref{item:F-lemma-cauchy1},
\[
\begin{aligned}
    (\alpha-\gamma)+(\beta-\delta)
    &=(uv+1)-(w^2+1)=F(x,x)F(y,y)-F(x,y)^2\geq0.
\end{aligned}
\]
Since their product is also nonnegative, both $\alpha-\gamma$ and $\beta-\delta$ are nonnegative. Using~\eqref{eq:F-lemma-root-order} and the definitions of $u,v,w$, we obtain the desired result from $\beta\geq\delta$.

To prove strictness, suppose that $x,y>1$ and $x\ne y$, so that $r,s>\sqrt2$ and $r\ne s$. We now show $\mathcal R<0$ strictly.
If $\mathcal R=0$,~\eqref{eq:F-lemma-resultant-square} gives either $w=1$ or $w=rs-1$. Substituting
this into~\eqref{eq:F-lemma-e3} gives either $e_3 =-(x-1)(y-1)$ or $e_3 =\frac{rs(r-s)^2}{(r+1)(s+1)}$,
both contradictory to $e_3 =0$. Thus $\mathcal R<0$.
Since $\beta\geq\delta$ has already been established,
and $\beta=\delta$ would imply $\mathcal R=0$, we obtain $\beta>\delta$, proving strictness.
\end{proof}
\end{lemma}

We next consider a general bivariate map $F$ and study growth, optimal splits, and remainder bounds for a new sequence $\{a_n\}_{n\geq1}$ defined by
\begin{equation}\label{eq:max-split-recurrence}
a_1>1,\qquad a_n=\max_{1\leq k<n}F(a_k,a_{n-k}),\qquad n\geq2.
\end{equation}

\begin{lemma}\label{lemma:upper-bound-Rn}
Let $\{a_n\}_{n=1}^\infty$ be defined by~\eqref{eq:max-split-recurrence}. For $x,y\geq 1$, suppose that $F$ is monotonically nondecreasing entrywise,  $xy\leq F(x,y)\leq Mxy$ holds  where $M\geq 1$ is  a constant, and $F(F(x,y),F(x,y))\leq F(F(x,x),F(y,y))$.
Put $D(x)=F(x,x)$, and define $D_0(x)=x$ and $D_r(x)=D(D_{r-1}(x))$ for $r\geq1$.
Then the following assertions hold.
\begin{enumerate}[label=(\roman*),ref=(\roman*),leftmargin=*]
\item\label{item:growth-rate-limit} The asymptotic growth rate satisfies
\[
\lim_{n\to\infty}\frac{1}{n}\log a_n =  \sup_{n\geq 1}\frac{1}{n}\log a_n=\lim_{r\to\infty} 2^{-r}\log D_r(a_1).
\]

\item\label{item:dyadic-optimal-split} Suppose in addition that
%$F$ is the function in~\cref{lemma:F-lemma-four-points},
there exists a constant $c>1$ such that
\begin{equation}
\label{eq:dyadic-additional-assumptions}
\log D(x)-\log D(y)\geq c\left(\log x-\log y\right),\qquad \forall x>y\geq 1.
\end{equation}
Then $a_{2^k}=D_k(a_1)$ for every $k\in\N$.
\item\label{item:dyadic-optimal-split-uniqueness}
Suppose~\eqref{eq:dyadic-additional-assumptions} holds and in addition
that $F(F(x,y),F(x,y))< F(F(x,x),F(y,y))$ is strict when $x,y>1$ and $x\neq y$.
Then, for every $k\geq1$, the unique maximizing index in~\eqref{eq:max-split-recurrence} at $n=2^k$ is $j=2^{k-1}$.
\end{enumerate}
\end{lemma}
\begin{proof}
\emph{\ref{item:growth-rate-limit}}
By induction, $a_1^n\leq a_n\leq a_1^n M^{n-1}$ for all $n\geq1$, so $(\log a_n)/n$ is bounded between $\log a_1$ and $\log a_1+\log M$. The recursive definition and $F(x,y)\geq xy$ give $a_n\geq F(a_k,a_{n-k})\geq a_k a_{n-k}$, and hence
\[
\log a_n\geq\log a_k+\log a_{n-k},\qquad n\geq2,\quad 1\leq k\leq n-1.
\]
By Fekete's lemma, $\lim_{n\to\infty}(\log a_n)/n=\sup_{n\geq1}(\log a_n)/n$ exists and is finite.

We next claim that, for every $x\geq1$, the limit
\begin{equation}\label{eq:lambda-definition}
\Lambda(x)\triangleq\lim_{r\to\infty}2^{-r}\log D_r(x)
\end{equation}
exists and is finite. Fix $x\geq1$. Since $t^2\leq D(t)\leq Mt^2$ for $t\geq1$ and $D_r(x)\geq1$, the recursion gives
\begin{equation}\label{eq:lemma-upper-bound-lambda-n}
0\leq \frac{\log D_r(x)}{2^r} -\frac{\log D_{r-1}(x)}{2^{r-1}}\leq \frac{\log M}{2^r},\qquad r\geq2.
\end{equation}
Thus $2^{-r}\log D_r(x)$ is nondecreasing in $r$. Summing~\eqref{eq:lemma-upper-bound-lambda-n} and using $x^2\leq D_1(x)\leq Mx^2$, we obtain
\begin{equation}\label{eq:lemma-lambda-bound}
    \log x\leq 2^{-r}\log D_r(x)\leq\log x+(1-2^{-r})\log M\leq\log x+\log M,\qquad r\geq1.
\end{equation}
Thus the sequence is nondecreasing and bounded above, which proves the claim.

Using monotonicity of $F$ and applying the inequality $D(F(x,y))\leq F(D(x),D(y))$ recursively gives, for every $r\geq1$,
\begin{equation}\label{eq:lemma-upper-bound-rn-lambda-ineq0}
    D_r(F(x,y))\leq 
% F(D_r(x),D_r(y)).
D _{r-1}(F(D (x),D (y)))
%\leq D^{r-2}(F(D^2(x),D^2(y)))
\leq \dotsm\leq F(D_r(x),D_r(y)).
\end{equation}
Because $F(x,y)\leq Mxy$, for any $r\geq1$,~\eqref{eq:lemma-upper-bound-rn-lambda-ineq0} implies that
\begin{equation}\label{eq:lemma-upper-bound-rn-lambda-ineq}
   2^{-r}\log D_r(F(x,y))\leq 2^{-r}\log (MD_r(x)D_r(y))=\frac{\log M+\log D_r(x)+\log D_r(y)}{2^r}.
\end{equation}
Taking $r\to\infty$ in~\eqref{eq:lemma-upper-bound-rn-lambda-ineq} and using the definition of $\Lambda$ in~\eqref{eq:lambda-definition}, we obtain $\Lambda(F(x,y))\leq\Lambda(x)+\Lambda(y)$.

We can then show by induction that $\Lambda (a_n)\leq n\Lambda (a_1)$. When $n=1$, it holds trivially. If it holds for all $n<N$, then for $a_N$ there exists $1\leq k<N$ such that $a_N=F(a_k,a_{N-k})$ and thus
\[
\Lambda (a_N)\leq \Lambda (a_k)+\Lambda (a_{N-k})\leq k\Lambda (a_1)+(N-k)\Lambda(a_1) = N\Lambda (a_1).
\]

To complete the proof of~\ref{item:growth-rate-limit}, we combine the following bounds.
On the one hand, $D(x)\geq x^2$ implies $D_r(x)\geq x^{2^r}$. Hence $\Lambda(x)\geq \log x$ and thus
\begin{equation}\label{eq:lemma-upper-bound-rn-log-upper}
    \frac{\log a_n}{n}\leq \frac{\Lambda (a_n)}{n}\leq \Lambda (a_1).
\end{equation}
On the other hand, $a_{2^r}\geq F(a_{2^{r-1}},a_{2^{r-1}})=D(a_{2^{r-1}})$. By induction and the monotonicity of $D$, we have $a_{2^r}\geq D_r(a_1)$ and hence
\begin{equation}\label{eq:lemma-upper-bound-rn-log-lower}
\lim_{r\to\infty}\frac{\log a_{2^r}}{2^r}\geq\lim_{r\to\infty} \frac{\log D_r(a_1)} {2^r}=\Lambda(a_1).
\end{equation}
Since $\lim_{n\to\infty}n^{-1}\log a_n$ exists by the first paragraph of this proof, combining~\eqref{eq:lemma-upper-bound-rn-log-upper} and~\eqref{eq:lemma-upper-bound-rn-log-lower} yields $\lim_{n\to\infty}n^{-1}\log a_n=\Lambda(a_1)$, as desired.

\emph{\ref{item:dyadic-optimal-split}}
% Suppose now that $F$ is the function in~\cref{lemma:F-lemma-four-points}.
Assume~\eqref{eq:dyadic-additional-assumptions}. We have
\begin{equation}\label{eq:lambda-dyadic-properties}
\Lambda(D_r(x))=2^r\Lambda(x),\qquad
\Lambda(F(x,y))\leq\Lambda(x)+\Lambda(y),\qquad
\Lambda(a_n)\leq n\Lambda(a_1),
\end{equation}
where the identity follows from $D_m(D_r(x))=D_{m+r}(x)$ and~\eqref{eq:lambda-definition}, the bound on $\Lambda(F(x,y))$ follows by letting $r\to\infty$ in~\eqref{eq:lemma-upper-bound-rn-lambda-ineq}, and the last inequality was established by induction in the proof of~\ref{item:growth-rate-limit}.

We next prove that $\Lambda$ is strictly increasing. For any $x>y\geq1$ and $r\in\N_+$,
\begin{equation}
    \label{eq:lemma-lambda-strictly-increasing}
    \begin{aligned}
    2^r(\Lambda (x)-\Lambda(y))&\stackrel{\eqref{eq:lambda-dyadic-properties}}{=} \Lambda(D_r(x))-\Lambda(D_r(y))
  \stackrel{\eqref{eq:lemma-lambda-bound}}{\geq} \log D_r(x)-\log D_r(y)-\log M\\ &
    \geq c^r\left(\log x-\log y\right)-\log M,
    \end{aligned}
\end{equation}
where the last inequality uses~\eqref{eq:dyadic-additional-assumptions} repeatedly. Choosing sufficiently
large $r$ such that the right-hand side of~\eqref{eq:lemma-lambda-strictly-increasing} is positive,
we obtain $\Lambda(x)>\Lambda(y)$, which proves strict monotonicity.
%  Since $D(x)=F(x,x)$,
% iterating the bounds $D(x)\geq x^2$ and $D(x)+c_1\leq (x+c_1)^2$ gives
% \[
% x^{2^r}\leq D_r(x)\leq(x+c_1)^{2^r}-c_1,\qquad x\geq1,\quad r\in\N.
% \]
% % \cref{lemma:F-properties}(i)--(ii) gives $x^2\leq D(x)$ and $D(x)+1\leq x^2+3\leq(x+1)^2$ for $x\geq1$. Iterating these bounds gives
% % \[
% % x^{2^r}\leq D_r(x)\leq(x+1)^{2^r}-1,\qquad x\geq1,\quad r\in\N.
% % \]
% Taking logarithms, dividing by $2^r$, and letting $r\to\infty$ yields, by~\eqref{eq:lambda-definition},
% \begin{equation}\label{eq:lambda-log-bracket}
% \log x\leq\Lambda(x)\leq\log(x+c_1),\qquad x\geq1.
% \end{equation}
% Differentiating the expression for $D=G$ in~\cref{lemma:F-properties} gives
% \[
% D'(x)=x+\frac{\sqrt{x^2+2x+5}}2+\frac{x^2-1}{2\sqrt{x^2+2x+5}}\geq1+\sqrt2,\qquad x\geq1.
% \]
% By the mean value theorem and induction, for $x>y\geq1$, it follows that
% By assumption~\eqref{eq:dyadic-additional-assumptions}, for $x>y\geq 1$, there exists
% sufficiently large $r\in\N$ such that $D_r(x)-D_r(y)\geq c_2^r(x-y)>c_1$.
% Thus $\Lambda$ is strictly increasing, since
% \[
% \begin{aligned}
% 2^r\Lambda(x)&\;\overset{\eqref{eq:lambda-dyadic-properties}}{=}\;\Lambda(D_r(x))\;
% \overset{\eqref{eq:lambda-log-bracket}}{\geq}\;\log D_r(x)>\log(D_r(y)+c_1)\;
% \overset{\eqref{eq:lambda-log-bracket}}{\geq}\;\Lambda(D_r(y))\;
% \overset{\eqref{eq:lambda-dyadic-properties}}{=}\;2^r\Lambda(y).
% \end{aligned}
% \]

By~\eqref{eq:lambda-dyadic-properties},
\[
\Lambda(a_{2^k})\leq2^k\Lambda(a_1)=\Lambda(D_k(a_1)).
\]
Strict monotonicity gives $a_{2^k}\leq D_k(a_1)$. Together with $a_{2^k}\geq D_k(a_1)$, established immediately before~\eqref{eq:lemma-upper-bound-rn-log-lower}, this yields $a_{2^k}=D_k(a_1)$.

\emph{\ref{item:dyadic-optimal-split-uniqueness}}
% We now prove the uniqueness under the additional assumption.
For $x,y>1$ with $x\ne y$, 
% \cref{lemma:F-lemma-four-points}(ii) gives
% \[
% D(F(x,y))=F(F(x,y),F(x,y))<F(F(x,x),F(y,y))=F(D(x),D(y)).
% \]
the strict monotonicity of $\Lambda$ and~\eqref{eq:lambda-dyadic-properties} yield $\Lambda(F(x,y))<\Lambda(x)+\Lambda(y)$, since
\begin{equation}\label{eq:lambda-strict-split}
\begin{aligned}
2\Lambda(F(x,y))&=\Lambda(D(F(x,y)))<\Lambda(F(D(x),D(y)))\leq\Lambda(D(x))+\Lambda(D(y))=2\Lambda(x)+2\Lambda(y),
\end{aligned}
\end{equation}
where the first strict inequality follows from the assumption.
Since $a_{n+1}\geq F(a_n,a_1)\geq a_na_1>a_n$, any index $j\ne2^{k-1}$ gives $a_j\ne a_{2^k-j}$ and therefore
\[
\Lambda(F(a_j,a_{2^k-j}))\stackrel{\eqref{eq:lambda-strict-split}}{<}\Lambda(a_j)+\Lambda(a_{2^k-j})\leq2^k\Lambda(a_1)=\Lambda(a_{2^k}).
\] 
Hence \(F(a_j,a_{2^k-j})<a_{2^k}\). The symmetric split attains $F(a_{2^{k-1}},a_{2^{k-1}})=D_k(a_1)=a_{2^k}$ and is therefore the unique maximizer.
\end{proof}

\begin{lemma}  [Lower bound] \label{lemma:lower-bound-Rn}

Let $\{a_n\}_{n=1}^\infty$ be defined by~\eqref{eq:max-split-recurrence}, where $F$ is monotonically nondecreasing entrywise and $F(x,y)\geq xy$ for $x,y\geq1$.
Define $v_0=a_1$ and $v_{r+1}=F(v_r,v_r)$ for $r\geq0$, and assume that the limit $\nu\triangleq\lim_{r\to\infty}2^{-r}\log v_r$ exists and is finite.
If the remainders $E_n\triangleq\log a_n-n\nu$ satisfy $E_n\leq0$ for all $n\in\N_+$, then
\[
E_n \geq \sum_{k=0}^\infty E_{2^k}\geq\sum_{k=0}^\infty(\log v_k-\nu2^k),\qquad n\in\N_+.
\]
\end{lemma}

\begin{proof} Since $F(x,y)\geq xy$, we have $a_n\geq a_ka_{n-k}$ for $1\leq k\leq n-1$. So $\log a_n-n\nu \geq (\log a_k-k\nu )+(\log a_{n-k}-(n-k)\nu)$. Hence $E_n\geq E_k+E_{n-k}$ for $1\leq k\leq n-1$.
For any $n\in\N_+$, this implies
\begin{equation}\label{eq:super-add-binary-decomp}
E_n\geq \sum_{j=1}^l E_{2^{p_j}}\geq \sum_{k=0}^\infty E_{2^k},
\end{equation}
where $n=\sum_{j=1}^l2^{p_j}$, with $0\leq p_1<\cdots<p_l$ and $p_j\in\N$, is the binary decomposition of $n$, and the last inequality uses $E_{2^k}\leq0$ for all $k\geq0$.
% Substituting the values of $E_{2^k}$ into~\eqref{eq:super-add-binary-decomp} gives $E_n\geq\sum_{k=0}^\infty(\log a_{2^k} - \nu 2^k)$. Finally, 
Finally, using $a_1=v_0$ as the base case, we obtain by induction
\[
a_{2^n}\geq F(a_{2^{n-1}},a_{2^{n-1}})\geq F(v_{n-1},v_{n-1})=v_n,\qquad n\geq1,
\]
where the second inequality uses the induction hypothesis $a_{2^{n-1}}\geq v_{n-1}$ and the coordinatewise monotonicity of $F$. Hence $  E_{2^k}\geq  \log v_k - \nu 2^k $ and the proof is complete using~\eqref{eq:super-add-binary-decomp}.\end{proof}

We now apply these estimates to the recurrence for the \OConPP{} schedules.
\begin{lemma}\label{lemma:bound-Rn} Let $\{a_n\}_{n=1}^\infty$ be defined by~\eqref{eq:max-split-recurrence}, where $F(x,y)=\frac{xy+1+\sqrt{(xy+3)^2-4(x+1)(y+1)}}{2}$. Let $v_0=a_1$ and $v_{r+1}=G(v_r)$ for $r\geq0$, where $G$ is defined in~\cref{lemma:F-properties}. Set $\nu\triangleq\lim_{r\to\infty}2^{-r}\log v_r$, which exists and is finite by~\cref{lemma:duplicate}\ref{item:duplicate-limit}.
Then $a_n=\exp\{n\nu +E_n\}$ with $-\frac{1}{2}<-\frac{1}{a_1(a_1+1)}\leq E_n\leq 0$ for all $n\in\N_+$
and $\lim_{k\to\infty}E_{2^k}= 0$. 
Moreover, $a_{2^{k}}=G(a_{2^{k-1}})$ for $k\in\N_+$ and the unique maximizing index in~\eqref{eq:max-split-recurrence}
at $n=2^k$ is $j=2^{k-1}$.
\end{lemma}    
\begin{proof}
The identity $(xy+3)^2-4(x+1)(y+1)=(xy-1)^2+4(x-1)(y-1)$ shows that $F$ is well-defined and monotonically nondecreasing entrywise for $x,y\geq1$. Moreover, $F(x,x)=G(x)$, and
\[
xy\leq F(x,y)\leq xy+2\leq3xy,\qquad x,y\geq1.
\]
Together with~\cref{lemma:F-lemma-four-points}(ii), these properties verify the assumptions of~\cref{lemma:upper-bound-Rn} with $M=3$. With $D_r$ defined in~\cref{lemma:upper-bound-Rn}, the identity $F(x,x)=G(x)$ gives $D_r(a_1)=v_r$ for $r\geq0$. Hence~\cref{lemma:upper-bound-Rn}\ref{item:growth-rate-limit} gives $\sup_{n\geq1}n^{-1}\log a_n=\nu$.
For each $n\in\N_+$, define $E_n\triangleq\log a_n-n\nu$. Since $n^{-1}\log a_n\leq\nu$, we have $E_n\leq0$.
Applying~\cref{lemma:lower-bound-Rn} and~\cref{lemma:duplicate}\ref{item:duplicate-sum}, with $v_0=a_1$, gives
\[
E_n \geq \sum_{k=0}^\infty E_{2^k}\geq\sum_{k=0}^\infty(\log v_k-\nu2^k)\geq-\frac{1}{a_1(a_1+1)}> -\frac12,\qquad n\in\N_+.
\]
Since $E_n\leq 0$, the above inequality also implies $E_{2^k}\to 0$ as $k\to \infty$.

We next verify~\eqref{eq:dyadic-additional-assumptions}. By Cauchy's mean value theorem and~\cref{lemma:F-properties}(iv),
for $x>y\geq 1$ there exists $z$ such that $y<z<x$ and
\[
\frac{\log G(x)-\log G(y)}{\log x-\log y}=\frac{zG'(z)}{G(z)}\geq\frac{50}{27}.
\]
Hence~\eqref{eq:dyadic-additional-assumptions} holds with $c=50/27$. This yields  $a_{2^{k}}=G(a_{2^{k-1}})$ for $k\in\N_+$ by~\cref{lemma:upper-bound-Rn}.
The uniqueness of the maximizing index follows from the strictness of~\cref{lemma:F-lemma-four-points} and~\cref{lemma:upper-bound-Rn}\ref{item:dyadic-optimal-split-uniqueness}.\end{proof}

\section{Proofs for Section~\ref{sect:connections-extensions}}\label{sect:proofs-section5}

\begin{proof}[Proof of~\Cref{coro:optimized-distance}]
Both bounds in~\eqref{eq:optimized-distance-bound} are sharp by~\cref{theorem:bound-x}, with equality attained for $f(x)=\mu x^2/2$ and any $x_0\neq0$.
Hence both minimum iteration budgets equal $\min\{N\in\N_+:\tau_N\leq\epsilon\}$.
Combining~\eqref{eq:distance-beta}, \eqref{eq:prim-asymp}, and~\eqref{eq:prim-remainder-bound} gives
\[
\tau_N=e^{-2N\nu(q)-2R_{N-1}^\circ},\qquad
e^{-2N\nu(q)}\leq\tau_N\leq e^{1-2N\nu(q)}.
\]
The upper and lower estimates yield
\begin{equation*}
\left\lceil\frac{\log(1/\epsilon)}{2\nu(q)}\right\rceil
\leq N_2(\kappa,\epsilon)\leq
\left\lceil\frac{\log(1/\epsilon)+1}{2\nu(q)}\right\rceil.
\end{equation*}
Using $x\leq \lceil x\rceil<x+1$ and $q^\rho\leq\nu(q)\leq {4q^\rho}/({1-q^2})$ from~\eqref{eq:nu-q-bound-q-pow-rho},
we obtain
\[
\frac{1-q^2}{8q^\rho}\log \frac{1}{\epsilon}\leq N_2(\kappa,\epsilon)\leq 
1+\frac{1}{2q^\rho}\left(1+\log \frac{1}{\epsilon}\right)
\leq 1+\frac{1}{q^\rho}\log  \frac{1}{\epsilon},
\]
where the last inequality uses $\log (1/\epsilon)\geq 1$. Using $\kappa=1/q$ gives the desired complexity.
% Since $\log(1/\epsilon)\geq1$, it follows that $N_\epsilon=\Theta(1+\nu(q)^{-1}\log(1/\epsilon))$.
% By~\cref{prop:nu-bound}, $\nu(q)=\Theta(q^\rho)$ uniformly for $0<q<3/13$, which gives the claimed complexity.
\end{proof}

\begin{proof}[Proof of~\Cref{prop:dynamic-congp}]
By induction,~\cref{theorem:congp,theorem:conpp} preserve the respective schedule classes and append $m+1$ steps at each concatenation. The bounds follow from~\cref{theorem:bound-of-dual,theorem:bound-of-fdom}, using~\cref{prop:primitive-is-fdom} in case~(ii). For $q=0$, the limits are established in~\cite[Propositions~5.10 and~4.18]{zhang2024accelerated}, respectively.
We now prove the results for $0<q<1$.

Let $0<q<1$. By~\cref{lemma:recursion-of-sigma}, both constructions satisfy
\begin{equation}\label{eq:dynamic-congp-recur}
% s_{k+1}=\frac{s_k+s_b+\gamma_k+q s_ks_b}{1-q\gamma_k}.
s_{k+1}=\frac{(1+qs_k)(1+qs_b)(1-q\gamma_k)^{-1}-1}{q}.
\end{equation}
Set $t_k=1+q s_k$. Equation~\eqref{eq:dynamic-congp-recur} gives
\begin{equation}\label{eq:dynamic-t-recursion}
t_{k+1}=\frac{1+q s_b}{1-q\gamma_k}t_k.
\end{equation}
By~\cref{lemma:roots-of-phi,lemma:roots-of-psi}, $1<\gamma_k<1/q$, so $s_k$ is increasing. Hence $\gamma_k\geq\gamma_0$ by~\cref{lemma:monotonicity-of-phi-psi}. Setting $A=(1+q s_b)/(1-q\gamma_0)>1$, we obtain $t_{k+1}\geq A t_k$ and thus $s_k\to+\infty$. Since $\theta_k=q/(t_k^2-(1-q))$,
\[
\frac{\theta_{k+1}}{\theta_k}
=\frac{t_k^2-(1-q)}{t_{k+1}^2-(1-q)}
\leq A^{-2},
\]
which proves the first inequality in~\eqref{eq:dynamic-linear-bound} by induction.

By~\cref{lemma:roots-of-phi,lemma:roots-of-psi,lemma:relation-phi-psi} and the symmetry of $\varphi$,
both connecting steps satisfy
\begin{equation}\label{eq:dynamic-gamma-bound}
1+\frac{q(1-q)(s_k+1)(s_b+1)}{q^2(s_k+1)(s_b+1)+(qs_k+1)(qs_b+1)}
<\varphi(s_k,s_b)
\leq\gamma_k<\frac{qs_b+s_b+2}{ 2q s_b+q+1}.
\end{equation}
The left-hand side of~\eqref{eq:dynamic-gamma-bound} converges to the right-hand side when
$k\to\infty$ and $s_k\to+\infty$. Therefore it implies the limits
% Dividing $\Phi_{s_k,s_b,q}(\varphi_q(s_k,s_b))=0$, with $\Phi$ defined in~\eqref{eq:Phi}, by $s_k$ and
%  letting $k\to\infty$ shows that $\varphi_q(s_k,s_b)$ converges to the displayed upper bound. Therefore,
\[
\lim_{k\to\infty}\gamma_k=\frac{qs_b+s_b+2}{ 2q s_b+q+1},
\qquad
\lim_{k\to\infty}\frac{t_{k+1}}{t_k}\stackrel{\eqref{eq:dynamic-t-recursion}}{=}\frac{ 2q s_b+q+1}{1-q}.
\]
Together with $t_k\to+\infty$ and $\theta_k=q/(t_k^2-(1-q))$, this proves the limit in~\eqref{eq:dynamic-linear-bound}.
\end{proof}

\begin{proof}[Proof of~\Cref{coro:dynamic-empty-block}]
Set $s_n=\sigma(h^{(n)})$. By~\cref{theorem:congp}, $h_n=\psi(0,s_n)\in(1,2)$, and each concatenation appends one step. Suppose first that $q\neq0$. Expanding $\Psi_{0,s_n,q}(h_n)=0$ gives
\begin{equation}\label{eq:adaptive-sn}
(q s_n+1)^2=\frac{(1-q h_n)^2(2-h_n)}{2-h_n(1+q)}.
\end{equation}
The numerator and the left-hand side are positive, so $2-h_n(1+q)>0$. Together with $h_n<2$, this gives the stated stepsize range. Also,
\begin{equation}\label{eq:dynamic-recur-qs1}
q s_{n+1}+1\overset{\eqref{eq:qs-add-1}}{=}\frac{q s_n+1}{1-q h_n}.
\end{equation}
Hence
\begin{equation*}
\frac{(1-q h_n)^2(2-h_n)}{2-h_n(1+q)}
\frac{2-h_{n+1}(1+q)}{(1-q h_{n+1})^2(2-h_{n+1})}
\overset{\eqref{eq:adaptive-sn}}{=}\frac{(q s_n+1)^2}{(q s_{n+1}+1)^2}
\overset{\eqref{eq:dynamic-recur-qs1}}{=}(1-q h_n)^2.
\end{equation*}
Cancelling and rearranging yields
\begin{equation*}
\frac{q h_n}{2-h_n(1+q)}
=\frac{(1-q h_{n+1})^2(2-h_{n+1})}{2-h_{n+1}(1+q)}-1
=-\frac{q h_{n+1}[(2-h_{n+1})(2-q h_{n+1})-1]}{2-h_{n+1}(1+q)}.
\end{equation*}
Dividing by $q$ proves~\eqref{eq:weakly-convex-recur} for $n\geq0$. For $n=-1$, it follows from~\eqref{eq:adaptive-sn} at $n=0$, since $s_0=h_{-1}=0$.

Combining~\eqref{eq:adaptive-sn} and~\eqref{eq:dynamic-recur-qs1} gives $(q s_{n+1}+1)^2=(2-h_n)/(2-h_n(1+q))$. Thus, for $n\geq1$,
\[
1+2s_n+q s_n^2
=1+\frac{(q s_n+1)^2-1}{q}
=1+\frac{h_{n-1}}{2-h_{n-1}(1+q)}.
\]
The same identity holds at $n=0$ because $s_0=h_{-1}=0$.

For $q=0$, the stepsize range follows from $h_n\in(1,2)$. The equation $\Psi_{0,s_n,0}(h_n)=0$ and the sum $s_{n+1}=s_n+h_n$ give
\[
2s_n=\frac{h_n(2h_n-3)}{2-h_n},
\qquad
2s_{n+1}=\frac{h_n}{2-h_n}.
\]
Comparing these identities at consecutive indices proves~\eqref{eq:weakly-convex-recur} for $n\geq0$. Its $n=-1$ case follows from $s_0=h_{-1}=0$. They also give $1+2s_n=1+h_{n-1}/(2-h_{n-1})$ for $n\geq0$.

In both cases, substituting the identity for $1+2s_n+q s_n^2$ into~\eqref{eq:dynamic-congp-bound} proves~\eqref{eq:weakly-convex-bound}.\end{proof}

\end{appendices}
\end{document}